\documentclass[journal]{IEEEtran}

\usepackage{amsmath,amsthm, amssymb, bm}
\usepackage{graphicx}%\graphicspath{{figures/}}
\usepackage{multicol}
\usepackage{enumerate}
\usepackage[square,numbers]{natbib}
\usepackage[colorlinks,citecolor=blue,urlcolor=blue,linkcolor=blue]{hyperref}
\usepackage{hypernat}
\usepackage[frame]{crop}
\usepackage{tikz}
\usetikzlibrary{shapes.geometric,arrows,chains,matrix,positioning,scopes,calc,decorations.markings,decorations.pathreplacing}
\tikzstyle{mybox} = [draw=white, rectangle]
\usepackage{booktabs}
\usepackage{multirow}
\usepackage{datetime}
\usepackage{thmtools}
\usepackage{mathrsfs}
\usepackage{float}

\usepackage{placeins}
\usepackage{stfloats}
\usepackage[bf,footnotesize]{caption}
\usepackage[capitalize]{cleveref} 
\crefname{lemma}{lemma}{lemmas}

\usepackage{framed}
 {\endMakeFramed}
\definecolor{shadecolor}{gray}{0.85}

\declaretheorem[style=plain,numberwithin=section,qed=$\square$]{theorem}

\declaretheorem[style=plain,sibling=theorem,qed=$\square$]{corollary}

\declaretheorem[style=definition,qed=$\triangleleft$,sibling=theorem]{definition}
\declaretheorem[style=definition,qed=$\triangleleft$,sibling=theorem]{example}
\declaretheorem[style=definition,qed=$\triangleleft$,sibling=theorem]{remark}
\declaretheorem[style=definition,qed=$\triangleleft$,sibling=theorem]{fact}

\numberwithin{equation}{section}
\numberwithin{theorem}{section}

\def\[#1\]{\begin{align}#1\end{align}}
\def\Ind{1\!\!1}
\def\defn#1{{\rm\textbf{#1}}}
\def\dee{\mathrm{d}}

\def\Uniform{\mbox{\rm Uniform}}
\def\Bernoulli{\mbox{\rm Bernoulli}}
\def\ie{i.e.,\ }
\def\eg{e.g.,\ }
\def\iid{i.i.d.\ }
\def\dist{\sim}
\def\distiid{\dist_{\mbox{\tiny iid}}}
\def\given{\mid}
\def\simiid{\distiid}
\def\eqdist{\overset{\mbox{\tiny d}}{=}}
       
\def\dataspace{\mathbf{X}}        % sample space
\def\Beta{\mbox{\rm Beta}}

\def\defas{:=}
\def\st{:}

\def\mean{\mathbb{E}}
\def\EE{\mean}

\def\Nats{\mathbb{N}}
\def\Reals{\mathbb{R}}
\def\PosReals{\Reals_{>0}}
\def\NNReals{\Reals_+}
\def\Ints{\mathbb{Z}}
\def\NNInts{\Ints_+}

\def\oset#1#2{#1#2} % Kallenberg-style, no comma even

\def\nset#1{\lbrace #1 \rbrace}
\def\mset#1{\nset{#1}}%{\{\!\{#1\}\!\}}

\def\structspace{\xspace^{\infty}}

\def\graphon{w}
\def\Graphon{W}

\def\graphonspace{\mathbf{W}}
\def\quotientspace{\widehat{\graphonspace}}
\def\xspace{\mathbf{X}}

\def\pMeas{\mathbf{M}}
\def\tspace{\mathbf{T}}
\def\model{\mathcal{P}}
\def\Pr{\mathbb{P}}

\def\argdot{\,.\,}

\newcommand{\genkernel}{\mathbf{p}}
\newcommand{\kernelfamily}[3]{{\{ \kernelval{#1}{#2} : #2 \in #3\}}} % \family \genkernal \theta T

\newcommand{\kernelvalset}[3]{{#1_{#3}(#2)}} % \kernelofset <kernel> <set> <parameter>
\newcommand{\kernelval}[2]{#1_{#2}}

\begin{document}
\nonfrenchspacing

\title{Bayesian Models of Graphs, Arrays and Other Exchangeable Random Structures}

\author{
  Peter Orbanz 
  and 
  Daniel M.\ Roy 
  \thanks{\hspace{-.3cm}\textbf{Published as}: Orbanz, P. and Roy, D.M. (2015). 
    Bayesian Models of Graphs, Arrays and Other Exchangeable Random Structures, 
    \textit{IEEE Trans.\ Pattern Analysis and Machine Intelligence}, 
    Vol. 37, No. 2, pp. 437--461.}
}

\maketitle

\begin{abstract}
  The natural habitat of most Bayesian methods is data represented by exchangeable sequences of observations,
  for which de~Finetti's theorem provides the theoretical foundation.
  Dirichlet process clustering, Gaussian process regression, 
  and many other parametric and nonparametric Bayesian models fall
  within the remit of this framework; many problems arising in modern
  data analysis do not.  This article provides an
  introduction to Bayesian models of graphs, matrices, and other data
  that can be modeled by random structures.  We describe results in
  probability theory that generalize de~Finetti's theorem to such data
  and discuss their relevance to nonparametric Bayesian
  modeling.  With the basic ideas in place, we survey example models
  available in the literature; applications of such models include
  collaborative filtering, link prediction, and graph and network
  analysis.  We also highlight connections to recent developments in
  graph theory and probability, and sketch the more general
  mathematical foundation of Bayesian methods for other types of data
  beyond sequences and arrays.
\end{abstract}

\section{Introduction}

For data represented by exchangeable sequences,
Bayesian nonparametrics has developed into a flexible and powerful toolbox of models and algorithms. 
Its modeling primitives---Dirichlet processes, Gaussian processes, etc.---are widely applied and well-understood, and
can be used as components in hierarchical models \citep{Teh:Jordan:2010:1} or dependent models \citep{MacEachern:2000}
to address a wide variety of data analysis problems.
One of the main challenges for Bayesian statistics and machine learning is arguably 
to extend this toolbox to the analysis of data sets with additional structure, such as graph, network, and relational data.

In this article, we consider structured data---sequences,
graphs, trees, matrices, etc.---and ask:
\begin{center}
{  \em
  What is the appropriate class of statistical models
  for a given type of structured data?
}
\end{center}
Representation theorems for exchangeable random structures lead us to an answer,
and they do so in
a very precise way: They characterize
the class of possible Bayesian models for the given type of 
data, show how these models are
parametrized, and even provide basic convergence guarantees.
The probability literature provides such results for dozens
of exchangeable random structures, including
sequences, graphs, partitions, arrays, 
trees, etc. The purpose of this article is to explain
how to interpret these results 
and 
how to translate them into a statistical
modeling approach.

\begin{center}
  {\bf Overview}
\end{center}

\begin{table*}[b]
  \caption{Exchangeable random structures}
  \centering{
    \begin{tabular}{llll}
      \toprule
      Random structure & Theorem of & Ergodic distributions $\kernelval \genkernel \theta$ & Statistical application \\
      \midrule
      Exchangeable sequences & de~Finetti \citep{DeFinetti:1930,DeFinetti:1937} & product distributions & 
      most Bayesian models \citep[e.g.][]{Schervish:1995}\\
      {} & Hewitt and Savage \citep{Hewitt:Savage:1955} &  \\
      Processes with exchangeable increments & B\"uhlmann \citep{Buehlmann:1960} & L\'evy processes\\
      Exchangeable partitions & Kingman \citep{Kingman:1978:2} & ``paint-box'' distributions & clustering\\
      Exchangeable arrays & Aldous \citep{Aldous:1981} & sampling schemes \cref{eq:jechar}, 
      \cref{eq:sechar} &
      graph-, matrix- and array-valued \\
      {} & Hoover \citep{Hoover:1979} & & data (\eg \citep{Hoff:2007:1}); see \cref{sec:models}\\
      {} & Kallenberg \citep{Kallenberg:1999} & {}\\
      Block-exchangeable sequences & Diaconis and Freedman \citep{Diaconis:Freedman:1980:1} & Markov chains &
      e.g. infinite HMMs \citep{Beal:Ghahramani:Rasmussen:2002,Fortini:Petrone:2012}\\
      \bottomrule
    \end{tabular}
    \label{tab:structures}
  }
\end{table*}

Consider a data analysis problem in which observations are represented as edges in a graph.
As we observe more data, the graph grows. For statistical purposes, we might model the graph as
a sample from some probability distribution on graphs. 
Can we estimate the distribution or, at least, some of its properties? 
If we observed multiple graphs, all sampled independently from the same distribution, we would 
be in the standard setting of statistical inference.
For many problems in e.g.\ network analysis, that is
clearly not the approach we are looking for: There 
is just one graph, and so our sample size should relate to the size of the graph.
Of course, we could assume that the edges are independent and identically
distributed random variables, and estimate their distribution---that would indeed
be a way of performing inference within a single graph. The resulting model, however, 
is sensitive only to the number of edges in a graph
and so has only a single parameter ${p\in[0,1]}$, the probability that an edge is present. 
What are more expressive models for graph-valued data?

Compare the problem to a more familiar one, where data is represented by a random sequence
${X:=(X_1,X_2,\ldots)}$, whose elements take values in a sample space $\xspace$,
and $n$ observations are interpreted as the values of the initial
elements ${X_1,\ldots,X_n}$ of the sequence. 
In this case, a statistical model is a family 
${\model=\lbrace P_{\theta} |\theta\in\tspace\rbrace}$ of distributions on
$\xspace$, indexed by elements $\theta$ of some parameter space $\tspace$.
If the sequence $X$ is exchangeable, de Finetti's theorem
tells us that there is \emph{some} model $\model$ and \emph{some} distribution
$\nu$ on $\tspace$ such that the joint distribution of $X$ is
\begin{equation}
  \label{intro:de:Finetti}
  \mathbb{P}(X\in\argdot)
  =
  \int_{\tspace}P_{\theta}^{\infty}(\argdot)\nu(d\theta)\;.
\end{equation}
This means the sequence can be generated by first generating
a random parameter value ${\Theta\sim\nu}$, and then sampling
${X_1,X_2,\ldots|\Theta\simiid P_{\Theta}}$. In particular, the elements
of the sequence are conditionally \iid given $\Theta$. 

For the purposes of statistical inference, the (conditional) \iid structure
implies that we can regard the elements $X_i$ of the sequence as
repetitive samples from an unknown distribution $P_{\Theta}$, and
pool these observations to extract information about the value of $\Theta$.
This may be done using a Bayesian
approach, by making a modeling assumption on $\nu$ (\ie by defining a prior
distribution) and computing the posterior given data, or in a frequentist
way, by assuming $\theta$ is non-random and deriving a suitable estimator.

To generalize this idea to the graph case, we regard the infinite
sequence ${X:=(X_1,X_2,\ldots)}$ as an infinite random structure, of which
we observe a finite substructure (the initial segment ${X_1,\ldots,X_n}$).
From this perspective, de Finetti's theorem tells us how to break down a
random structure (the sequence) into components (the conditionally independent elements), which
in turn permits the definition of statistical inference procedures. What
if, instead of an infinite sequence, $X$ is an infinite graph, of which
we observe a finite subgraph? To mimic the sequence case,
we would need a definition of exchangeability applicable to graphs, and a suitable
representation theorem to substitute for de Finetti's.

Since generating a parameter ${\Theta\sim\nu}$ randomly determines a distribution
$P_{\Theta}$, we can think of $P_{\Theta}$ as a random probability measure with distribution
defined by $\nu$, and paraphrase de Finetti's theorem as follows:
\begin{center}
{  \em
  The joint distribution of any exchangeable sequence of random values
  in $\xspace$ is characterized
  by the distribution of a random probability measure on
  $\xspace$.
}
\end{center}

If we assume
a random graph to be exchangeable---where we put off the precise definition
for now---it is indeed also characterized by a representation theorem.
The implications for statistical models are perhaps
more surprising than in the case of sequences:
\begin{center}
{  \em
  The distribution of any exchangeable 
  graph is characterized by a distribution on the space of
    functions \\ from $[0,1]^2$ to $[0,1]$.
}
\end{center}
Hence, any specific function ${w:[0,1]^2\rightarrow[0,1]}$ defines a distribution $P_w$ on graphs 
(we will see in \cref{sec:2-arrays} how we can sample a graph from $P_w$).

For modeling purposes, this means that \emph{any} statistical model of
exchangeable graphs is a family of such distributions $P_w$, and $w$ can
be regarded as the model parameter.
Density estimation in exchangeable graphs can therefore be formulated as a
regression problem: It is equivalent to recovering the function $w$
from data.
Once again, we can choose a frequentist approach (define
an estimator for $w$) or a Bayesian approach (define a prior
distribution on a random function $W$); we can obtain nonparametric
models by choosing infinite-dimensional subspaces of functions, or
parametric models by keeping the dimension finite.

Since a graph can be regarded as a special type of matrix (the 
adjacency matrix), we can ask more 
generally for models of exchangeable matrices, and obtain a similar
result:
\begin{center}
  \em
  The distribution of any exchangeable two-dimensional, real-valued array is characterized by a 
  distribution on the space of functions from $[0,1]^3$ to $\mathbb{R}$.
\end{center}
There is a wide variety of random structures for which exchangeability
can be defined; Table \ref{tab:structures} lists some important
examples. Borrowing language from \citep{Aldous:2010:1}, we
collectively refer to such random objects as \defn{exchangeable random
structures}. 
This article explains representation theorems for exchangeable random
structures and their implications for Bayesian statistics and machine
learning.
The overarching theme is that key aspects of de Finetti's
theorem can be generalized to many types of data, and that these results
are directly applicable to the derivation and interpretation of
statistical models.

\begin{center}
  {\bf Contents}
\end{center}
\begin{trivlist}
\item[{\bf \cref{sec:structures}}:] reviews exchangeable random structures, their representation theorems, and
  the role of such theorems in Bayesian statistics.
\item[{\bf \cref{sec:2-arrays}}:] introduces the generalization of de~Finetti's theorem to models of
  graph- and matrix-valued data, the Aldous-Hoover theorem, and explains
  how Bayesian models of such data can be constructed.
\item[{\bf \cref{sec:models}}:] surveys models of graph- and relational data available in the machine learning and statistics literature.
  Using the Aldous-Hoover representation, models can be classified 
  and some close connections emerge between models which seem, at first glance, only loosely related.
\item[{\bf \cref{sec:graph:limits}}:] describes recent development in the mathematical theory of graph limits. The results of this theory
  refine the Aldous-Hoover representation of graphs and provide a precise understanding of how graphs converge and how
  random graph models are parametrized.
\item[{\bf \cref{sec:d-arrays}}:] explains the general Aldous-Hoover representation for higher-order arrays.
\item[{\bf \cref{sec:sparsity}}:] discusses sparse random structures and networks, why these models contradict exchangeability, and open questions
  arising from this contradiction.
\item[{\bf \cref{sec:references}}:] provides references for further reading.
\end{trivlist}

\section{Bayesian Models of Exchangeable Structures}
\label{sec:structures}

\def\pMeas{\mathbf{M}}

\begin{it}
  The fundamental Bayesian modeling paradigm based on exchangeable sequences
  can be extended to a very general approach, where data 
  is represented by a random structure. 
  Exchangeability properties are then used to deduce valid statistical
  models and useful parametrizations.
  This section sketches out the ideas underlying this approach, 
  before we focus on graphs, matrices, and arrays in \cref{sec:2-arrays}.
\end{it}

\newcommand{\RS}{X}
\subsection{Basic example: Exchangeable sequences}
\label{sec:de:finetti}

The simplest example of an exchangeable random structure is an
exchangeable sequence. We use the customary shorthand notation
${(x_i):=(x_1,x_2,\ldots)}$ for a sequence, and similarly
$(x_{ij})$ for a matrix, etc.
Suppose  
$(X_{i})$
is an infinite sequence of random variables in a sample space
$\xspace$.
We call $(X_i)$ \defn{exchangeable} if its joint distribution satisfies
\[
  \label{eq:exchangeability:joint:distributions}
  &\Pr(X_1 \in A_1, X_2\in A_2,\dotsc) \\&= 
  \Pr(X_{\pi(1)} \in A_1,X_{\pi(2)} \in A_2,\dotsc) \nonumber
\]
for every permutation $\pi$ of ${\Nats \defas \{1,2,\dotsc\}}$ 
and every collection of sets $A_1,A_2,\dotsc$.
Expressing distributional equalities this way is cumbersome,
and we can write 
\eqref{eq:exchangeability:joint:distributions} 
more concisely as 
\[
  (X_1,X_2,\dotsc)
  \eqdist
  (X_{\pi(1)},X_{\pi(2)},\dotsc)\;,
\]
or even $(X_i) \eqdist (X_{\pi(i)})$, 
where the notation $Y\eqdist Z$ means that the random variables $Y$
and $Z$ have the same distribution.
Informally, exchangeability means that the probability of observing a
particular sequence does not depend on the order of the elements in the sequence.

\newcommand{\cvar}{\vartheta}
\newcommand{\cvarval}{t}

If the elements of a sequence are exchangeable, de~Finetti's representation
theorem implies they are conditionally i.i.d. The conditional
independence structure is represented by a 
\defn{random probability measure}, a random variable with values in the
set $\pMeas(\dataspace)$ of probability distributions on $\dataspace$.

\begin{theorem}[de~Finetti]
\label{theorem:definetti}
  Let ${(X_1,X_2,\dotsc)}$ be an infinite sequence 
  of random variables with values in a space $\dataspace$.
  The sequence $X_1,X_2,\dotsc$ is exchangeable if and only if
  there is a random probability measure $\Theta$ on $\dataspace$
  such that the $X_i$ are conditionally i.i.d. given $\Theta$ and
  \[
  \label{eq:de:finetti}
  \Pr(X_1 \in A_1,X_2 \in A_2,\dotsc) 
  &=\int_{\pMeas(\dataspace)}\prod_{i=1}^{\infty}\theta(A_i) \,\nu(\dee \theta)
  \]
  where $\nu$ is the distribution of $\Theta$. 
\end{theorem}

The integral on
the right-hand side of \eqref{eq:de:finetti}
can be interpreted as a two-stage sampling procedure: 
\begin{enumerate}
\item Sample ${\Theta\sim\nu}$, \ie draw a probability distribution at
  random from the distribution $\nu$.
\item Conditioned on $\Theta$, sample the $X_n$ conditionally \iid as
  \begin{equation}
    X_1,X_2,\ldots|\Theta\simiid\Theta\;.
  \end{equation}
\end{enumerate}

The theorem says that \emph{any} exchangeable sequence can
be sampled by such a two-stage procedure; the distribution of the
sequence is determined by the choice of $\nu$.
The random measure $\Theta$ is called the \defn{directing
  random measure} 
of $X$. Its distribution
$\nu$ is called the \defn{mixing measure} or \defn{de~Finetti
  measure}.

Statistical inference is only possible if the distribution of the
data, or at least some of its properties, can be recovered from
observations. For \iid random variables, this is ensured by the
law of large numbers. 
The proof of de Finetti's theorem also implies a law of large
numbers for exchangeable sequences:
\begin{theorem}\label{theorem:definetti:convergence}
 If the sequence $(X_i)$ is exchangeable, the empirical distributions
    \[
    \label{eq:emp:measure}
    \hat{S}_n(\,.\,):=\frac{1}{n}\sum_{i=1}^n\delta_{X_i}(\,.\,)
    \]
    converge to $\Theta$, in the sense that
    \[
    \label{eq:lln}
    \hat{S}_n(A) \to \Theta(A) \quad \text{as} \quad n \to \infty
    \]
    holds with probability 1 for every set $A$.  
\end{theorem}

The two theorems have fundamental implications for Bayesian
modeling. If we assume the data can be represented by (some finite prefix of) 
an exchangeable sequence, this implies 
without any further assumptions:
\begin{itemize}

\item 
 Conditioned on a random probability measure $\Theta$ representing an unknown distribution $\theta$, every sample $X_1,X_2,\dotsc$ is \iid with  distribution $\Theta$.

\item 
   Every exchangeable sequence model is characterized by
  a unique distribution $\nu$ on $\pMeas(\xspace)$.

\item 
   A statistical model can be taken to be some subset of $\pMeas(\xspace)$
  rather than $\pMeas(\xspace^{\infty})$, which we would have to consider for a general random sequence.

\item Statistical inference is possible in principle:  
   With probability one, the empirical distributions $\hat S_n$ converge to the distribution $\Theta$ generating the data, according to \eqref{eq:lln}.

\end{itemize}
A modeling application might look like this:
We consider a specific data source or measurement
process, and assume that data generated by this source can be
represented as an exchangeable sequence. The definition of
exchangeability for an infinite sequence does not mean we have
to observe an infinite number of data points to invoke de Finetti's
theorem; rather, it expresses the assumption that samples of
\emph{any finite} size generated by the source would be exchangeable.
Hence, exchangeability is an assumption on the data source, rather
than the data.

According to de Finetti's theorem, the data can then be explained by  
the two-stage sampling procedure above, for \emph{some} distribution $\nu$ on 
$\pMeas(\xspace)$. 
A Bayesian model is specified by choosing a specific distribution
$\nu$, the prior distribution.
In this abstract formulation of the prior as a measure on
$\pMeas(\xspace)$, the prior also determines
the observation model, as the smallest set
${\model\subset\pMeas(\xspace)}$ on which $\nu$ concentrates all its
mass---since ${\Theta}$ then takes values in $\model$,
and the sequence $(X_i)$ is generated by a
distribution in $\model$ with probability 1.
If ${\xspace=\mathbb{R}}$, for example, we could choose $\nu$ to
concentrate on the set of all Gaussian distributions on $\mathbb{R}$,
and would obtain a Bayesian model with a Gaussian likelihood and prior
$\nu$.

Given observations ${X_1,\ldots,X_n}$, we then compute the posterior
distribution, by conditioning $\nu$ on the observations.
\cref{theorem:definetti:convergence} implies that,
if the empirical measure converges asymptotically to a specific
measure ${\theta\in\pMeas(\xspace)}$, the posterior
converges to a point mass at $\theta$. This result has to be
interpreted very cautiously, however: It only holds for a sequence
$(X_i)$ which was actually generated from the measure $\nu$ we use as
a prior. In other words, suppose someone generates $(X_i)$ from a
distribution $\nu_1$ on $\pMeas(\xspace)$ by the two-stage sampling
procedure above, without disclosing $\nu_1$ to us. In the sampling
procedure, the variable ${\Theta\sim\nu_1}$ assumes as its value a 
specific distribution $\theta_1$, from which the data is then
generated independently.
We model the observed sequence by choosing a prior $\nu_2$. 
The posterior under $\nu_2$ still converges to a point mass,
but there is \emph{no guarantee} that it
is a point mass at $\theta_1$, and \eqref{eq:lln} 
only holds if ${\nu_2=\nu_1}$.

Thus, there are several important questions that
exchangeability does not answer: 
\begin{itemize}
\item The de Finetti theorem says that there is \emph{some} prior
  which adequately represents the data, but provides no guidance
  regarding the choice of $\nu$: Any probability measure $\nu$ on
  $\pMeas(\xspace)$ is the prior for some exchangeable sequence.
\item \cref{theorem:definetti:convergence} 
  \emph{only} guarantees convergence for sequences of random variables
  generated from the prior $\nu$.
\item \cref{theorem:definetti:convergence} is a first-order result:
  It provides no information
  on how quickly the sequence converges. Results on convergence rates
  can only be obtained for more specific models; the set of all
  exchangeable distributions is too large and too complicated
  to obtain non-trivial statements.
\end{itemize}
Answers to these questions typically require further modeling assumptions.

\subsection{The general form of exchangeability results}
\def\Xinf{X_{\infty}}

Many problems in machine learning and modern statistics involve
data which is more naturally represented by
a random structure that is not a sequence:   often a graph, matrix, array, tree,
partition, etc.\ is a better fit.
If it is possible to define a suitable notion of
exchangeability, the main features of de Finetti's
theorem typically generalize. Although results differ in their
details, there is a general pattern, which we sketch in this
section before considering specific types of exchangeable structures.

The setup is as follows: 
The product space $\xspace^{\infty}$ of infinite
sequences is substituted by a suitable space 
$\xspace_{\infty}$ of more general, infinite structures. 
An infinite random structure
$\Xinf$ is a random variable with
values in $\xspace_{\infty}$.
Each element of $\xspace_{\infty}$ can be thought of
as a representation of an infinitely large data set or
``asymptotic'' sample. An actual, finite sample of size $n$
is modeled as a substructure $X_n$ of $\Xinf$, such as a the length-$n$ prefix of an infinite sequence or
a $n$-vertex subgraph of an infinite graph.

The first step in identifying a notion of exchangeability is to specify what it 
means to permute components of a structure
${x_{\infty}\in\xspace_{\infty}}$. If $x_{\infty}$ is
an infinite matrix, for example, a very useful
notion of exchangeability arises when one considers all permutations that exchange the ordering of rows/columns, rather than the ordering of individual entries. 
Exchangeability of a
\emph{random} structure $\Xinf$ then means that the distribution
of $\Xinf$ is invariant under the specified family of permutations.

Once a specific exchangeable random structure $\Xinf$ is defined,
the next step is to invoke a representation theorem that generalizes 
de Finetti's theorem to $\Xinf$. Probability theory provides such
theorems for a range of random structures; see \cref{tab:structures} 
for examples. 
A 
representation theorem can be interpreted as determining (1) a natural parameter space $\tspace$ for exchangeable models
on $\xspace_{\infty}$, and (2) a special family  
of distributions on $\xspace_{\infty}$, which are called the
\defn{ergodic distributions} or \defn{ergodic measures}.
Each element ${\theta\in\tspace}$ determines an ergodic distribution,
and we denote this distribution as ${\kernelval \genkernel \theta}$.
The set of ergodic distributions is
\begin{equation}
  \kernelfamily \genkernel \theta \tspace   \subset\pMeas(\xspace_{\infty})\;.
\end{equation}
The distribution of any exchangeable random
structure $\Xinf$ can then be represented as a mixture of these ergodic
distributions,
\begin{equation}
  \label{eq:integral:decomp}
  \Pr(X_{\infty}\in\argdot)
  =
  \int_{\tspace} \kernelvalset{\genkernel}{\argdot}{\theta} \,\nu(\dee \theta)\;.
\end{equation}
In the specific case of exchangeable sequences,
\eqref{eq:integral:decomp} is precisely the integral representation
\eqref{eq:de:finetti} in de Finetti's theorem, and
the ergodic measures are
the distributions of \iid sequences, that is,
\begin{equation}
  \tspace:=\pMeas(\xspace)
  \qquad\text{ and }\qquad
  \kernelvalset{\genkernel}{\argdot}{\theta}=\theta^{\infty}(\argdot)
  \;.
\end{equation}
For more general random structures, the 
ergodic measures are not usually product distributions, 
but they retain some key properties:
\begin{itemize}
\item They are particularly simple distributions
  on $\xspace_{\infty}$, and form a ``small'' subset of all
  exchangeable distributions.
\item They have a conditional independence property, in the sense
  that a random structure $\Xinf$ sampled from one of the ergodic
  distributions decomposes into conditionally 
  independent components.
  In de Finetti's theorem, these conditionally independent components
  are the elements $X_i$ of the sequence. 
\end{itemize}
As in the sequence case, the integral \eqref{eq:integral:decomp} in
the general case represents a 
two-stage sampling scheme:
\begin{equation}
  \begin{split}
    \label{eq:hierarchy}
    \Theta &\dist \nu\\
    X_{\infty} \given \Theta &\dist \kernelval{\genkernel}{\Theta}\;.
  \end{split}
\end{equation}
For Bayesian modeling, this means:
\begin{center}
  \emph{ A Bayesian model for an exchangeable random structure
    $\Xinf$ with representation \eqref{eq:integral:decomp} is 
      characterized by a prior distribution on $\tspace$.}
\end{center}
Suppose the prior $\nu$ concentrates on a subset
${\mathcal{T}\subset\tspace}$, that is, $\mathcal{T}$ is the
smallest subset to which the prior assigns probability 1. 
Then $\mathcal{T}$ defines a subset
\begin{equation}
  \model:= \kernelfamily \genkernel \theta {\mathcal T} 
\end{equation}
of ergodic measures. We thus have defined a Bayesian model on 
$\xspace_{\infty}$, with prior $\nu$ and observation model $\model$.
In summary:
\begin{itemize}
\item $\tspace$ is the natural parameter space for Bayesian models
  of $\Xinf$, and the prior
  distribution is a distribution on $\tspace$.
\item The observation model $\model$ is a subset of the ergodic measures.
  An exchangeability theorem characterizing the ergodic measures
  therefore also characterizes the possible observation models.
\item The representation \eqref{eq:integral:decomp} is typically 
  complemented by a convergence results: A specific function of the samples
  converges to $\Theta$ almost surely as ${n\to\infty}$,
  generalizing \cref{theorem:definetti:convergence}. In particular,
  the parameter space $\tspace$ can be interpreted as the set of 
  all possible limit objects.
\end{itemize}
If the set $\mathcal{T}$ on which the prior concentrates its
mass is a finite-dimensional subspace of $\tspace$, we call the resulting
Bayesian model \defn{parametric}. If $\mathcal{T}$ has infinite
dimension, the model is \defn{nonparametric}.

\subsection{Exchangeable partitions}
\label{sec:kingman}

An illustrative example of an exchangeable random structure is an exchangeable
partition. Bayesian nonparametric clustering models are based on such
exchangeable random partitions.
We again define the exchangeable structure as an infinite object: Suppose
${X_1,X_2,\dotsc}$ is a sequence of observations. To encode a
clustering solution, we have to specify which observations $X_i$ belong
to which cluster. To do so, it suffices to record which index $i$ belongs
to which cluster, and a clustering solution can hence be expressed as
a partition ${\pi=(b_1,b_2,\dotsc)}$ of the index set $\Nats$. 
Each of the sets $b_i$, called \defn{blocks},
is a finite or infinite subset of $\Nats$; every element of $\Nats$ is contained in exactly one block.
An \defn{exchangeable partition} is a random partition $X_{\infty}$ of $\Nats$ which is invariant under permutations of $\Nats$.
Intuitively, this means the probability of a partition depends only
on the relative sizes of its blocks, but not on which elements are in which block.
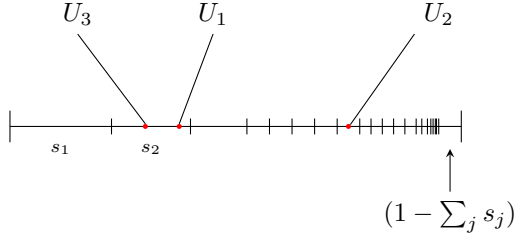
\begin{figure}
  \centering{
    \begin{tikzpicture}[scale=1.5,>=stealth]%[transform canvas={xshift=-0.5cm}]
      \draw[] (-0,0) -- (4.0,0) node {};
      \draw (-0,4pt) -- (-0,-4pt);
      \draw (4.0,4pt) -- (4.0,-4pt);
      \foreach \x/\xtext in {0.9/, 1.6/, 2.1/, 2.3/, 2.5/, 2.7/, 2.9/, 3.1/, 3.2/, 3.3/, 3.4/, 3.5/, 3.6/, 3.65/, 3.7/, 3.73/, 3.75/, 3.77/, 3.78/, 3.8/}
      \draw[shift={(\x,0)},font=\scriptsize] (0pt,2pt) -- (0pt,-2pt) node[below] {$\xtext$};
      \node[font=\scriptsize] at (0.45,-0.2) {$s_1$};
      \node[font=\scriptsize] at (1.25,-0.2) {$s_2$};
      \node[circle,fill,scale=0.2,color=red] (U3) at (1.2,0) {};
      \node[circle,fill,scale=0.2,color=red] (U1) at (1.5,0) {};
      \node[circle,fill,scale=0.2,color=red] (U2) at (3,0) {};
      \node (label3) at (0.6,1) {$U_3$};
      \node (label2) at (3.8,1) {$U_2$};
      \node (label1) at (1.8,1) {$U_1$};
      \draw (label1.south) -- (U1.north);
      \draw (label2.south west) -- (U2.north east);
      \draw (label3.south) -- (U3.north);
      \node (label0) at (3.9,-0.8) {$(1-\sum_j s_j)$};
      \draw[->] (label0.north)--($(label0)+(0,0.6)$);
    \end{tikzpicture}
  }
  \caption{Sampling from a paint-box distribution with parameter ${\mathbf{s}=(s_1,s_2,\dotsc)}$.
    Two numbers $i,j$ are assigned to the same block of the partition if the uniform variables $U_i$ and $U_j$ are 
    contained in the same interval.} 
  \label{fig:kingman}
\end{figure}

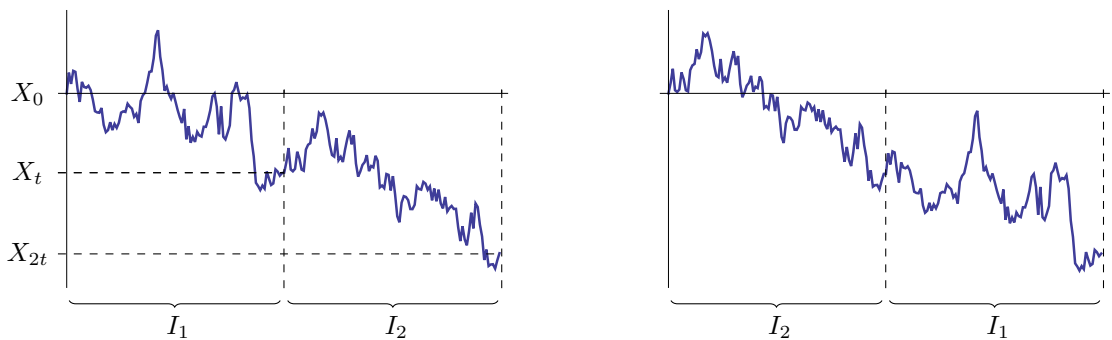
\begin{figure*}[b!]
  \caption{Exchangeable continuous-time processes: Shown on the left is a sample path of a process $X$ started
    at ${X_0=0}$. 
    We define two intervals of equal length, ${I_1:=(0,t]}$ and ${I_2:=(t,2t]}$. A permuted path
    is obtained by swapping both the intervals and the respective path segment of the process on each interval.
    The path segments are shifted vertically such that the path starts again at 0 and is "glued together"
    at the interval boundary. If $X$ is exchangeable, the permuted process defined by permuting paths in this
    manner has the same distribution as $X$. 
  }
  \begin{center}
    \begin{tikzpicture}[
        mirrorbrace/.style={
          decoration={brace, mirror},
          decorate},brace/.style={
          decoration={brace},
          decorate}
      ]
      \begin{scope}
        \node at (-.1,0.2) {\includegraphics[width=6cm]{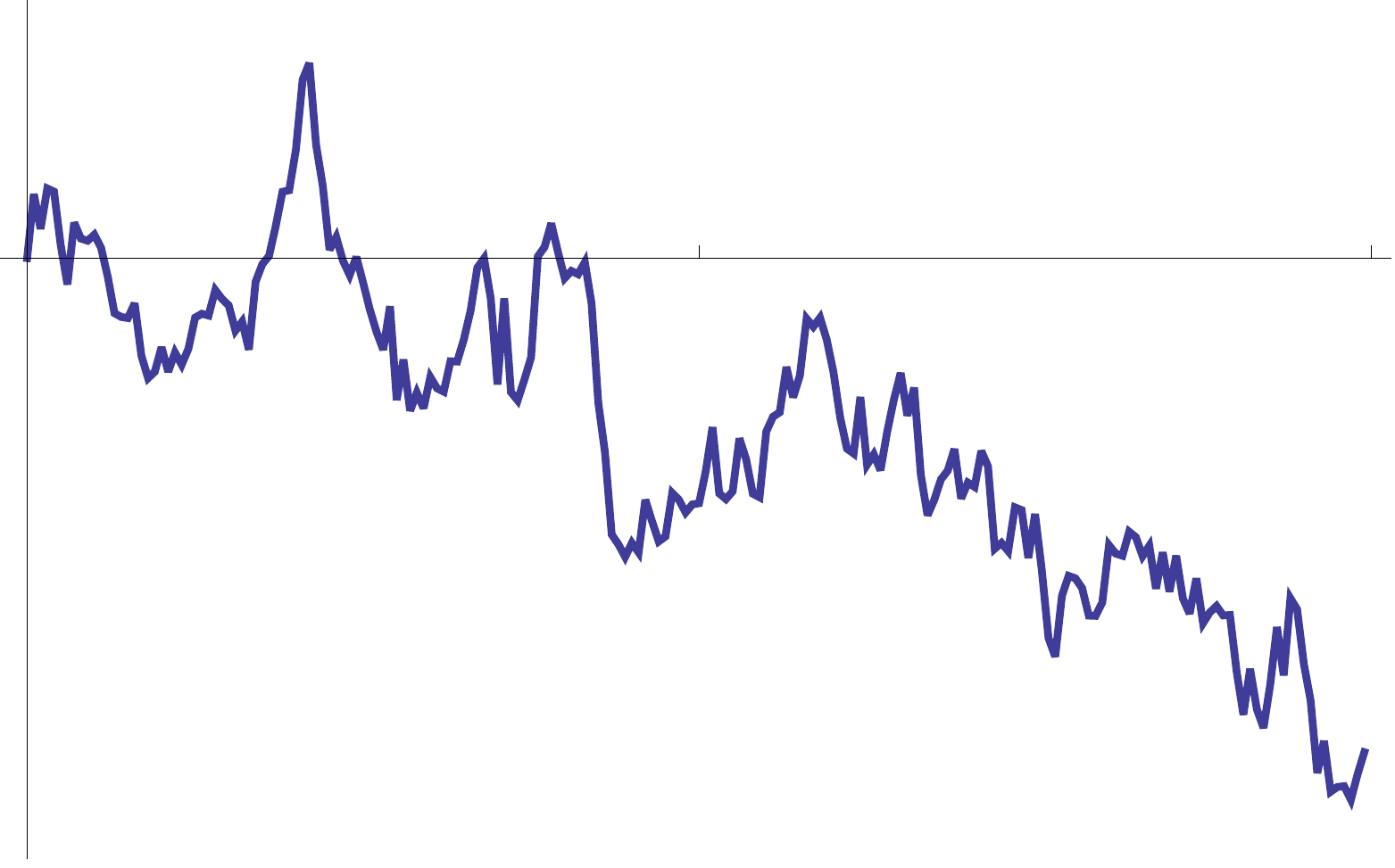}}; %OU2.pdf}};
        \draw[dashed] (-.095,1)--(-.095,-1.7);
        \draw[mirrorbrace] (-2.95,-1.8)--(-.145,-1.8);
        \draw[dashed] (2.8,1)--(2.8,-1.7);
        \draw[mirrorbrace] (-.045,-1.8)--(2.75,-1.8);
        \node at (-1.5,-2.15) {$I_1$};
        \node at (1.4,-2.15) {$I_2$};
        \draw[dashed] (-3.1,-.1)--(-.095,-.1);
        \node at (-3.5,0.95) {$X_0$};
        \draw[dashed] (-3.1,-.1)--(-.095,-.1);
        \node at (-3.5,-.1) {$X_{t}$};
        \draw[dashed] (-3.1,-1.18)--(2.8,-1.18);
        \node at (-3.5,-1.18) {$X_{2t}$};
      \end{scope}
      \begin{scope}[xshift=8cm]
        \node at (-.1,0.2) {\includegraphics[width=6cm]{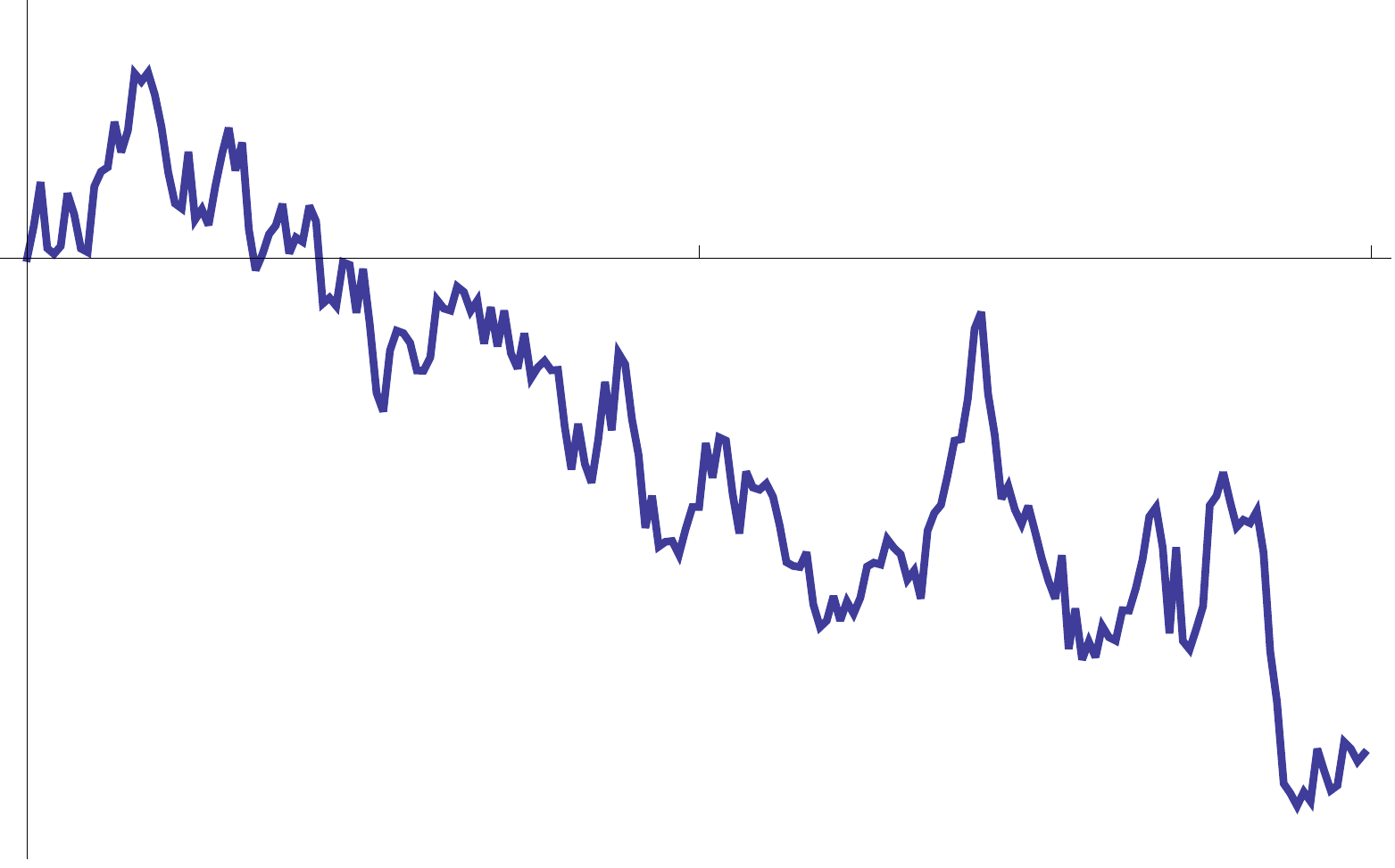}};
        \draw[dashed] (-.095,1)--(-.095,-1.7);
        \draw[mirrorbrace] (-2.95,-1.8)--(-.145,-1.8);
        \draw[dashed] (2.8,1)--(2.8,-1.7);
        \draw[mirrorbrace] (-.045,-1.8)--(2.75,-1.8);
        \node at (-1.5,-2.15) {$I_2$};
        \node at (1.4,-2.15) {$I_1$};
      \end{scope}
    \end{tikzpicture}
  \end{center}
  \label{fig:buehlmann}
\end{figure*}

Kingman \citep{Kingman:1978:2} showed that exchangeable random partitions can again be represented
in the form of \cref{eq:integral:decomp}. The parameter space $\tspace$
consists of all sequence 
${\theta:=(s_1,s_2,\dotsc)}$ of scalars ${s_i\in [0,1]}$
which satisfy 
\begin{equation}
  \label{eq:constraints:paintbox}
  s_1\geq s_2\geq\dots
  \quad\text{ and }\quad
  \sum_{i} s_i\leq 1\;.
\end{equation}
Let $\bar s_n \defas \sum_{i=1}^{n} s_i$.
Then $\theta$ defines a partition of
$[0,1]$ into intervals
\begin{equation}
  I_j:=\bigl[ \bar s_{j-1}, \bar s_j \bigr)
    \quad\text{ and }\quad
    \bar{I}:=\bigl(1- \bar s_\infty,1\bigr]\;,
\end{equation}
as shown in \cref{fig:kingman}. 
Each ergodic distribution $\kernelval \genkernel \theta$ is defined as
the distribution of the following random partition 
of $\mathbb{N}$:
\begin{enumerate}
\item Generate ${U_1,U_2,\dotsc\simiid\mbox{Uniform}[0,1]}$.
\item Assign ${n\in\mathbb{N}}$ to block $b_j$ if ${U_n\in I_j}$.
  Assign every remaining element (those $n$ with ${U_n\in\bar{I}}$) to its own block of size one.
\end{enumerate}
Kingman called this distribution a \defn{paint-box distribution}. 
\begin{theorem}[Kingman]
  \label{theorem:kingman}
  Let $X_{\infty}$ be random partition of $\mathbb{N}$.
  \begin{enumerate}
  \item $X_{\infty}$ is exchangeable if and only if
    \begin{equation}
      \label{eq:theorem:kingman}
      \Pr(X_{\infty}\in\argdot)=\int_{\tspace} \kernelvalset \genkernel \argdot \theta \, \nu(\dee \theta)\,,
    \end{equation}
    for some distribution $\nu$ on $\tspace$, 
    where ${\kernelval p \theta}$ is the paint-box distribution with parameter ${\theta\in\tspace}$.
  \item If $X_{\infty}$ is exchangeable, the scalars $s_i$ can be recovered asymptotically as
    limiting relative block sizes
    \begin{equation}
      \label{eq:kingman:asymptotic:frequencies}
      s_i=\lim_{n\to\infty}\frac{|b_i\cap\lbrace 1,\dotsc,n\rbrace|}{n}\;.
    \end{equation}
  \end{enumerate}
\end{theorem}

Part 1) is of course the counterpart to de Finetti's theorem, and part
2) corresponds to \cref{theorem:definetti:convergence}.
In \eqref{eq:kingman:asymptotic:frequencies}, we compute averages 
\emph{within} a single random structure, having observed
only a substructure of size $n$. Nonetheless, we can recover the 
parameter $\theta$ asymptotically from data. This is a direct
consequence of exchangeability, and would not generally be true
for an arbitrary random partition.

\begin{example}[Chinese restaurant process]
  A well-known example of a random partition is the Chinese restaurant
  process (CRP; see e.g.~\citep{MR2245368,Hjort:Holmes:Mueller:Walker:2010} for details). 
  The CRP is a one-parameter discrete-time stochastic process that induces a partition of $\Nats$.
  The parameter ${\alpha > 0}$ is called the \emph{concentration}; different
  values of $\alpha$ correspond to different distributions
  $\Pr(X_{\infty}\in\argdot)$ in \cref{eq:theorem:kingman}.
  If $X_{\infty}$ is generated by a CRP, the paint-box parameter $\Theta$ is essentially the
  sequence of weights generated by the ``stick-breaking'' construction of the
  Dirichlet process \citep{Hjort:Holmes:Mueller:Walker:2010}---with the difference that the elements of
  $\Theta$ are ordered by size, whereas stick-breaking weights are not.
  In other words, if $X_{\infty}$ in \eqref{eq:theorem:kingman} is a
  CRP, we can sample
  from $\nu$ by (1) sampling from a stick-breaking
  representation and (2) ordering the sticks by length. The lengths of
  the ordered sticks are precisely the scalars $s_i$ in the theorem.
\end{example}

\subsection{``Non-exchangeable'' data}
\label{sec:nonexchangeable}

Exchangeability seems at odds with many types of data; for example, a
sequence of stock prices over time 
will be poorly modeled by an exchangeable sequence. Nonetheless, a Bayesian model of
a time series will almost certainly imply an exchangeability assumption---the crucial
question is which components of the overall model are assumed to be exchangeable.
As the next example illustrates, these components need not be the variables representing
the observations.

\def\Levy{L\'evy }
\begin{example}[\Levy processes and B\"uhlmann's theorem]
  A widely-used class of models for real-valued time series in continuous time are L\'evy processes.
  The sample path of such a process is a random function $X$ on $\mathbb{R}_+$ that is piece-wise
  continuous (in more technical terms: right-continuous with left-hand limits).
  It is customary to denote the function value of $X$ at time ${t\in\mathbb{R}_+}$ as $X_t$.
  Recall the definition of a \Levy process: 
  If ${I=(t_1,t_2]}$ is an interval, then ${\Delta X_{I}:=X_{t_2}-X_{t_1}}$ is
  called an \defn{increment} of $X$. A process is called \defn{stationary} if the increments
  ${\Delta X_{I_1}}$ and ${\Delta X_{I_2}}$ are identically distributed whenever the intervals $I_1$ and $I_2$ 
  have the same length. We say that $X$ has \defn{independent increments} if 
  the random variables ${\Delta X_{I_1}}$ and ${\Delta X_{I_2}}$ are independent whenever $I_1$ and $I_2$ do
  not overlap. If $X$ is both stationary and has independent increments, it is a 
  \defn{\Levy process}.

  In other words, $X$ is a \Levy process if, for any disjoint ${I_1,I_2,\ldots}$ of equal length,
  the increments ${(\Delta X_{I_1},\Delta X_{I_2},\ldots)}$ form an \iid sequence.
  It is then natural to ask for an exchangeable process:
  We say that $X$ is an \defn{exchangeable continuous-time process} if,
  again for any disjoint ${I_1,I_2,\ldots}$ of equal length,
  the sequence ${(\Delta X_{I_1},\Delta X_{I_2},\ldots)}$ is exchangeable
  (see \cref{fig:buehlmann}).
  The representation theorem for such processes is due to Hans B\"uhlmann
  \citep[e.g.][Theorem 1.19]{Kallenberg:2005}:
  \begin{center}
    {\it
      A piece-wise continuous stochastic process on $\mathbb{R}_+$ is an
      exchangeable continuous-time process if and only if 
      it is a mixture of L\'evy processes.
    }
  \end{center}
  Hence, each ergodic measure $\kernelval \genkernel \theta$ is the distribution
  of a L\'evy process, and the measure $\nu$ is a distribution on parameters of L\'evy
  processes or---in the parlance of stochastic process theory---on L\'evy characteristics.
\end{example}

\begin{example}[Discrete times series and random walks]
  Another important type of exchangeability property 
  \citep{Diaconis:Freedman:1980:1,Zabell:1995:1} is defined for sequences
  $X_1,X_2,\dotsc$ taking values in a countable space $\xspace$. 
  Such a sequence is called \defn{Markov exchangeable} if the probability of observing an initial trajectory $x_1,\dotsc,x_n$ depends
  only on the initial state $x_1$ and, for every pair
  ${y,y'\in \xspace}$, on the number of transitions ${t_{y,y'}=\# \{ j < n \st x_j = y, x_{j+1} = y' \}}$.  In particular, the probability does not depend on when each transition occurs.
  Diaconis and Freedman \citep{Diaconis:Freedman:1980:1} showed the following:
  \begin{center}
    {\it
    If a (recurrent) process is Markov exchangeable, it is a mixture of Markov chains.
    }
  \end{center}
  (Recurrence means that each visited state is visited infinitely often if the process
  is run for an infinite number of steps.)
  Thus, each ergodic distribution $\kernelval \genkernel \theta$ is the distribution of a Markov
  chain, and a parameter value $\theta$ consists of a distribution on $\xspace$ (the distribution
  of the initial state) and a transition matrix.
  If a Markov exchangeable process is substituted for the Markov chain in a hidden Markov model,
  \ie if the Markov exchangeable variables are latent variables of the model, the resulting
  model can express much more general dependencies than Markov exchangeability. 
  The infinite hidden Markov model
  \citep{Beal:Ghahramani:Rasmussen:2002} is an example; see \citep{Fortini:Petrone:2012}.
  Recent work by \citet*{Bacallado:Favaro:Trippa:2013:1} constructs prior distributions on
  random walks that are Markov exchangeable and almost surely
  reversible.
\end{example}

A very general approach to modeling is to assume that an exchangeability
assumption holds marginally at each value of a covariate variable
$z$, \eg a time or a location in space:
Suppose $\structspace$
is a set of structures as described above, and $\mathbf{Z}$ is a space of covariate values.
A \defn{marginally exchangeable random structure} is a random measurable mapping
\begin{equation}
  \label{eq:marginally:exchangeable}
  \xi:\mathbf{Z}\to\xspace_{\infty}
\end{equation}
such that, for each $z\in\mathbf{Z}$, the random variable $\xi(z)$ is an exchangeable random structure
in $\xspace_{\infty}$.
\begin{example}[Dependent Dirichlet process]
  A popular example of a marginally exchangeable model
  is the dependent Dirichlet process (DDP) of \citet{MacEachern:2000}.
  In this case, for each ${z\in\mathbf{Z}}$, the random variable $\xi(z)$ is a random
  probability measure whose distribution is a Dirichlet process. 
  More
  formally, $\mathbf{Y}$ is some sample space, ${\xspace_{\infty}=\pMeas(\mathbf{Y})}$,
  and the DDP is a distribution on mappings ${\mathbf{Z}\to\pMeas(\mathbf{Y})}$;
  thus, the DDP is a random conditional probability.
  Since $\xi(z)$ is a Dirichlet process if $z$ is fixed, samples from $\xi(z)$ are exchangeable.
\end{example}
\cref{eq:marginally:exchangeable} is, 
of course, just another way of saying that $\xi$ is a $\xspace_{\infty}$-valued stochastic
process indexed by $\mathbf{Z}$, although we have made no specific requirements on the paths
of $\xi$. The interpretation as a path is more apparent in the next example.
\begin{example}[Coagulation- and fragmentation models]\hfill
  If $\xi$ is a coagulation or fragmentation process, 
  $\xspace_{\infty}$ is the set of partitions of $\Nats$ (as in Kingman's theorem), and
  $\mathbf{Z}=\NNReals$.
  For each $z\in\NNReals$, the random variable
  $\xi(z)$ is an exchangeable partition---hence, Kingman's theorem is applicable marginally
  in time. Over time, the random partitions become consecutively finer (fragmentation processes) or coarser (coagulation
  processes): At random times, a randomly selected block is split, or two randomly selected blocks merge.
  We refer to \citep{Bertoin:2006} for more details and to \citep{Teh:Blundell:Elliott:2011}
  for applications to Bayesian nonparametrics.
\end{example}

\subsection{Random functions vs random measures}
\label{sec:de:finetti:revisited}

\begin{figure}[t]
  \begin{center}
    \begin{tikzpicture}
      \draw[domain=0.001:1,samples=200,color=blue,thick] plot (6*\x,{1/5*(7+ln(\x/(1-\x)))}) ;
      \draw[->] (-0.2,0) -- (6.8,0) node[right] {}; 
      \draw[->] (0,-.2) -- (0,3.2) node[above] {$x$};
      \draw (6,-0.1)--(6,0.1) node[pos=0,label=below:{$1$}] {};
      \node[label=below:{$0$}] at (0,0) {}; 
      \node[label=left:{$a$}] at (0,0) {}; 
      \node[label=below:{$U_i$}] (U) at (3.4,0) {};
      \draw ($(U)-(0,0.1)$)--($(U)+(0,0.1)$);
      \draw[color=gray] (0,2.55)--(6,2.55);
      \draw[color=gray] (6,0)--(6,2.55);
      \draw (-0.1,2.55)--(0.1,2.55) node[pos=0,label=left:{$b$}] {};
      \node[label=left:{$X_i$}] (X) at (0,1.45) {};
      \draw ($(X)-(0.1,0)$)--($(X)+(0.1,0)$);
      \draw[dashed] (X)--($(X)+(3.4,0)$);
      \draw[dashed] (U)--($(U)+(0,1.45)$);
      \node[color=blue] at (4.5,1.3) {$F$};
    \end{tikzpicture}
  \end{center}
  \caption{de~Finetti's theorem expressed in terms of random functions: 
    If $F$ is the inverse CDF of the random measure $\Theta$ in the de Finetti representation,
    $X_i$ can be generated as $X_i:=F(U_i)$, where $U_i\sim\Uniform[0,1]$.}    
  \label{fig:de:finetti}
\end{figure}
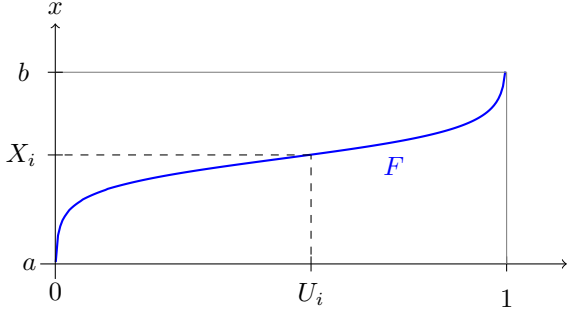

De Finetti's theorem can be equivalently formulated in terms of a random function, rather than a random measure,
and this formulation provides some useful intuition for \cref{sec:2-arrays}.
Roughly speaking, this random function is the inverse of the cumulative distribution function (CDF) 
of the random measure $\Theta$ in de Finetti's theorem;
see \cref{fig:de:finetti}.

\newcommand{\cdf}{c.d.f.}
More precisely, suppose that ${\xspace=[a,b]}$.
A measure $\mu$ on $[a,b]$ can be represented by its CDF,
defined as ${\psi(x) \defas \mu([a,x])}$. 
Hence, sampling the random measure $\Theta$ in de Finetti's theorem is equivalent to sampling
a random CDF $\Psi$. A CDF is not necessarily an invertible function, but it always admits a
so-called right-continuous inverse $\overline{\psi^{-1}}$, given by
\[
\overline{\psi^{-1}}(u) = \inf \, \lbrace x\in[a,b] \given \psi(x) \geq u\rbrace \;.
\]
This function inverts $\psi$ in the sense that $\psi\circ\overline{\psi^{-1}}(u)=u$ for all $u\in[0,1]$.
It is well-known that any scalar random variable $X_i$ with CDF $\psi$ can be generated as
\begin{equation}
  X_i\eqdist\overline{\psi^{-1}}(U_i)\qquad\text{ where } U_i\sim\Uniform[0,1]\;.
\end{equation}
In the special case $\xspace=[a,b]$, de Finetti's theorem therefore translates as follows:
If $X_1,X_2,\dotsc$ is an exchangeable sequence, then there is a random function 
$F \defas \overline{\Psi^{-1}}$ such that
\[
 (X_1,X_2,\dotsc) \eqdist (F(U_1),F(U_2),\dotsc)\;,
\]
where $U_1,U_2,\dotsc$ are \iid uniform variables.

It is much less obvious that the same should hold on an arbitrary sample space,
but that is indeed the case:
{
\begin{corollary}
  \label{thm:aldousdefinetti}
  Let $X_1,X_2,\dotsc$ be an infinite, exchangeable sequence 
  of random variables with values in a space $\dataspace$.
  Then there exists a random function $F$ from $[0,1]$ to $\dataspace$ such that,  
   if $U_1,U_2,\dotsc$ is an \iid sequence of uniform random variables,
\[
 (X_1,X_2,\dotsc) \eqdist (F(U_1),F(U_2),\dotsc).
\]
\end{corollary}
}
As we will see in the next section, this random function representation
generalizes to the more complicated case of array data, whereas the random measure
representation in \cref{eq:de:finetti} does not.

\section{Exchangeable Graphs, Matrices, and Arrays}
\label{sec:2-arrays}

\begin{it}
  Random arrays are a very general type of random
  structure, which include important special cases, such as
  random graphs and random matrices. The representation theorem 
  for exchangeable arrays is the Aldous-Hoover theorem.
  In this section, we focus on $2$-arrays, matrices and graphs.
  The general case for $d$-arrays is conceptually similar, but
  considerably more technical, and we postpone it
  until \cref{sec:d-arrays}.
\end{it}

A \defn{$d$-array} is a collection of elements
${x_{i_1,\ldots,i_d}\in\xspace}$ indexed by $d$ indices ${i_1,\ldots,i_d\in\Nats}$.
A sequence is a $1$-array. In this section, we assume the random
structure $\Xinf$ to be a random $2$-array 
\begin{equation}
  \Xinf 
  =
  (X_{ij})
  =
  \left(\begin{matrix} X_{11} & X_{12} & \ldots \\
    X_{21} & X_{22} & \ldots \\
    \vdots & \vdots & \ddots
  \end{matrix}\right)\;.
\end{equation}
A \defn{random matrix} is a random $2$-array, although the
term matrix usually implies that $\xspace$ has the algebraic
structure of a field. A \defn{random graph} is a random matrix
with ${\xspace=\lbrace 0,1\rbrace}$.
As in the sequence case, we assume $\Xinf$ is infinite in size, and
the statistical interpretation is that an observed, finite array
is a sub-array of $\Xinf$. In network analysis problems, for example,
an observed graph with $n$ vertices would be interpreted as
a random induced subgraph of an underlying graph $\Xinf$, which represents
an infinite population.

In this section, we are interested in the characterization of random arrays whose distributions are invariant to permutations reordering the rows and columns.
For a $2$-array, there are two natural ways to define exchangeability: we can ask that the distribution of the array be invariant only to the \emph{joint} (simultaneous) permutation of the rows and columns or also to  \emph{separate} permutations of rows and columns.

\begin{definition}
  A random 2-array $(X_{\oset i j})$ is called 
  \defn{jointly exchangeable} if
  \begin{equation}
    \label{eq:def:jointly:2D}
    (X_{\oset i j}) \eqdist (X_{\oset {\pi(i)} {\pi(j)} } )
  \end{equation}
  for every permutation $\pi$ of $\Nats$,
  and \defn{separately exchangeable} if
  \begin{equation}
    \label{eq:def:separately:2D}
    (X_{\oset i j }) \eqdist (X_{\oset {\pi(i)} {\pi'(j)}})
  \end{equation}
  for every pair of permutations $\pi,\pi'$ of $\Nats$.
\end{definition}
Invariance to all separate permutations of the rows and the columns is
an appropriate assumption if rows and columns
correspond with two \emph{distinct} sets of entities, such as in a collaborative filtering problem, 
where rows may correspond to users and columns to movies. 
On the other hand, if $(X_{\oset ij})$ is the adjacency matrix of a
random graph on the vertex set $\Nats$, we would require joint exchangeability, because
there is only a
single set of entities---the vertices of the graph---each  
of which corresponds both to a row and a column of the matrix.

\subsection{The Aldous-Hoover theorem}

The analogue of de~Finetti's theorem
for exchangeable arrays is the 
\defn{Aldous-Hoover theorem} \citep{Aldous:1981,Hoover:1979}.
It has two versions, for jointly and for separately exchangeable arrays.

\begin{theorem}[Aldous-Hoover]%[Jointly exchangeable arrays]
  \label{theorem:AH:2D:jointly}
  A random array $(X_{\oset i j })$ is jointly exchangeable if and only if it can be represented as follows:
  There is a random function $F:[0,1]^3\to\dataspace$ such that
  \begin{equation}
    \label{eq:AH:2D:jointly}
    (X_{\oset i j })
    \eqdist 
    (F(U_i,U_j,U_{\mset{i,j}})) \;,
  \end{equation}
  where $(U_i)_{i\in\Nats}$ and $(U_{\mset{i,j}})_{i,j\in\Nats}$ are, respectively, 
  a sequence and an array of i.i.d.\ $\Uniform[0,1]$
  random variables, which are independent of $F$.
\end{theorem}
Because the variables $U_{\mset{i,j}}$ are indexed by a set, the indices
are unordered, and we can think of the array $(U_{\mset{i,j}})$ as an
upper-triagonal matrix with \iid uniform entries.
If the function $F$ is symmetric in its first two arguments, then
$\Xinf$ is symmetric---that is, if
$F(x,y,\,.\,)=F(y,x,\,.\,)$ for all $x$ and $y$, then ${X_{ij}=X_{ji}}$ for
all $i$ and $j$. In general, however, a jointly exchangeability
matrix or $2$-array need not be symmetric.

Separately exchangeable arrays can also be given a precise characterization using \cref{theorem:AH:2D:jointly}: 

\begin{figure*}[t]
  \begin{center}
    \begin{tikzpicture}[>=stealth]%,transform canvas={xshift=-3cm,yshift=1cm}]
      \begin{scope}%[yshift=0.5cm]
        \node [mybox] (box){
          \includegraphics[width=3cm]{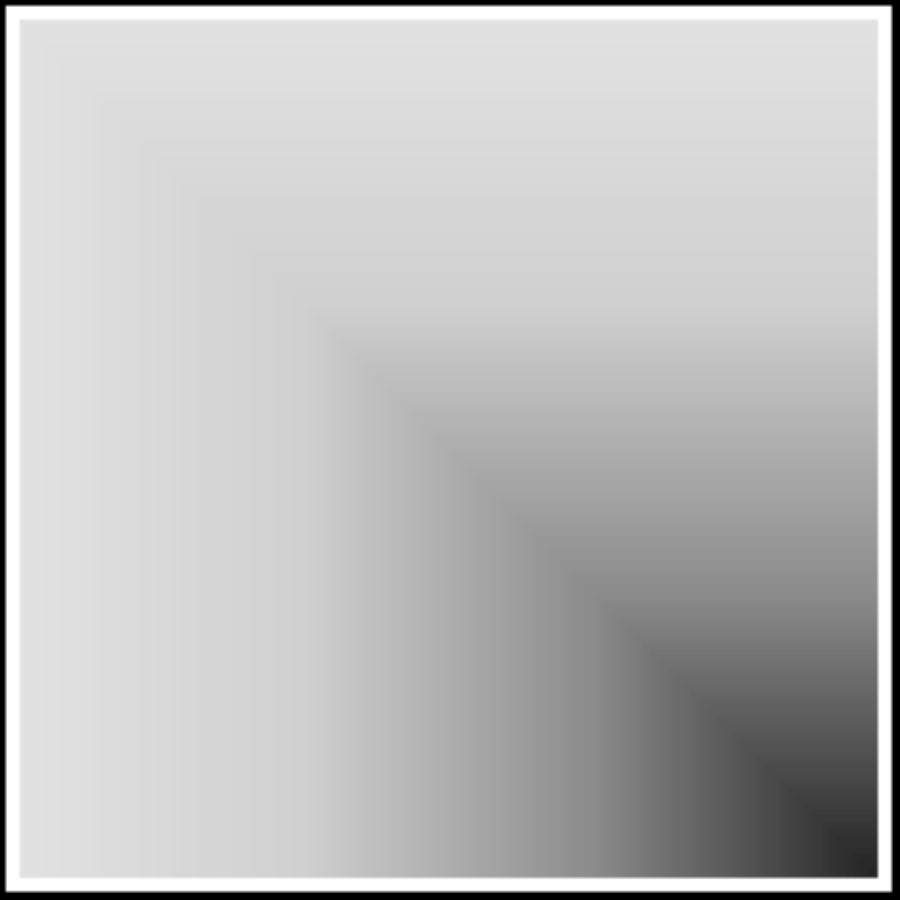} 
        };
        \node[scale=0.2,font=\scriptsize,label=left:{$0$}] at (-1.5,1.5) {};
        \node[scale=0.2,font=\scriptsize,label=above:{$0$}] at (-1.5,1.5) {};
        \node[scale=0.2,font=\scriptsize,label=right:{$1$}] at (1.5,-1.5) {};
        \node[scale=0.2,font=\scriptsize,label=below:{$1$}] at (1.5,-1.5) {};
        \node[label=above:$U_1$] at (-0.8,1.5) {};
        \draw [dashed] (0.3,0.8) -- (0.3,1.6); 
        \node[label=above:$U_2$] at (0.3,1.5) {};
        \draw [dashed] (-1.6,0.8) -- (0.3,0.8); 
        \node[label=left:$U_1$] at (-1.5,0.8) {};
        \node[label=left:$U_2$] at (-1.5,-0.3) {};
        
        \draw[thick] (-0.8,1.4)--(-0.8,1.6);
        \draw[thick] (0.3,1.4)--(0.3,1.6);
        \draw[thick] (-1.4,0.8)--(-1.6,0.8);
        \draw[thick] (-1.4,-0.3)--(-1.6,-0.3);
        
        \node[circle,fill,scale=0.4,color=red] (dot) at (0.3,0.8) {};
        
        \pgfmathsetmacro{\loc}{3.2}
        \draw (\loc,1.5)--(\loc,-1.5);
        \draw (\loc-0.05,-1.5)--(\loc+0.05,-1.5); \node[scale=0.2,label=right:$0$] at (\loc,-1.5) {};
        \draw (\loc-0.05,1.5)--(\loc+0.05,1.5); \node[scale=0.2,label=right:$1$] at (\loc,1.5) {};
        \node (P) at (\loc,-0.21) {};
        \draw[color=red,thick] ($(P)+(-0.1,0)$)--($(P)+(0.1,0)$);
        \node[font=\scriptsize,label=right:${\Graphon(U_1,U_2)}$] at (P) {};
        \node (Uij) at (\loc,0.4) {};
        \draw[color=black,thick] ($(Uij)+(-0.1,0)$)--($(Uij)+(0.1,0)$);
        \node[font=\scriptsize,label=right:${U_{\mset{1,2}}}$] at (Uij) {};
        \draw[->,color=black] (dot.east)--(P.west) node[pos=0.7,fill=white] {$\Graphon$} ;
      \end{scope}
      \begin{scope}[xshift=7.2cm]
        \node [mybox] (box){
          \includegraphics[width=3cm]{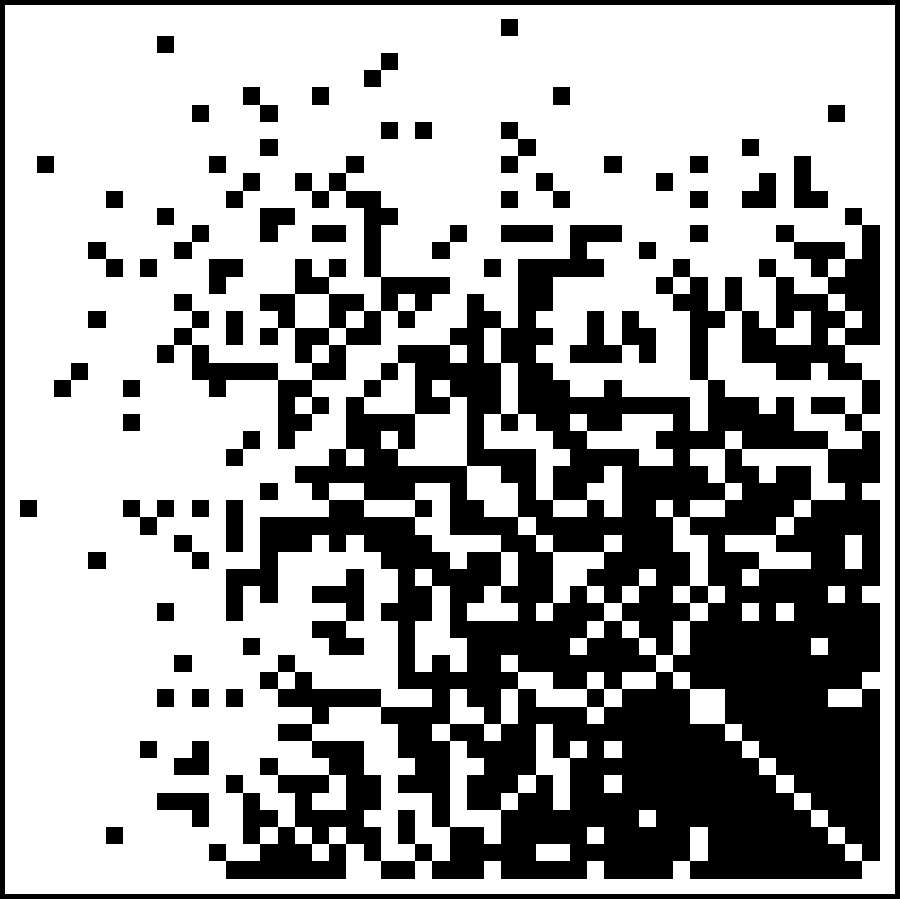}
        };
      \end{scope}
      \begin{scope}[xshift=12cm]
        \node [mybox] (box){
          \includegraphics[width=5cm]{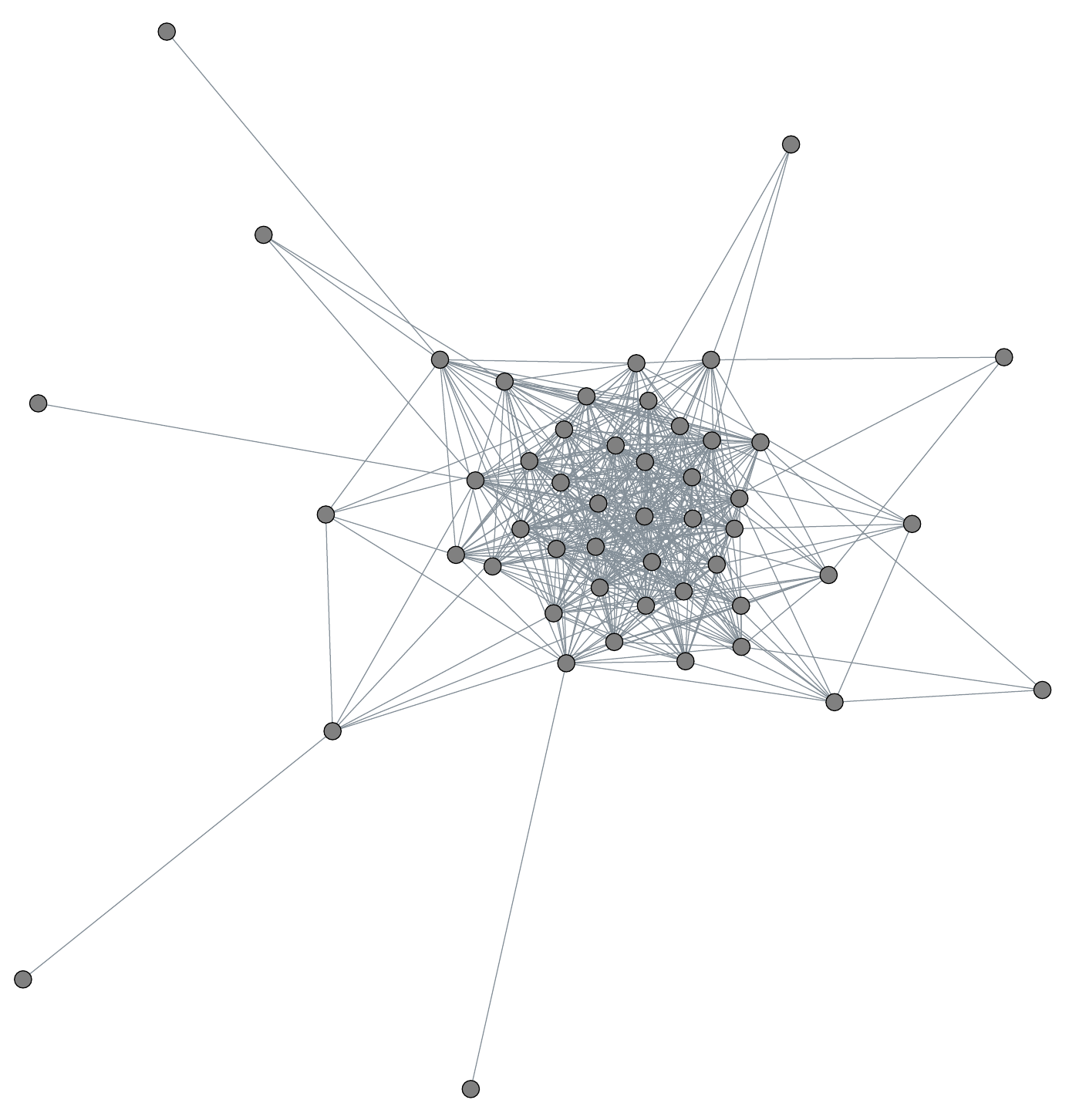}
        };
      \end{scope}
    \end{tikzpicture}
  \end{center}
  \vspace{-0.5cm}
  \caption{
    Sampling an exchangeable random graph according to \cref{eq:AH:graph}. 
    \emph{Left:} A heat-map visualization of the random function $\Graphon$ on $[0,1]^2$, given here by $\Graphon(x,y)=\min\lbrace x,y\rbrace$. Darker shades represent larger values. In the case depicted here, the edge $\{1,2\}$ is not present in the graph, because $U_{\mset{1,2}}>\Graphon(U_1,U_2)$. 
    \emph{Middle:} The adjacency matrix of a 50-vertex random graph, sampled from the function on the left. Rows (and columns) in the matrix have been ordered by their underlying $U_i$ value, resulting in a matrix resembling $\Graphon$.
    \emph{Right:} A plot of the random graph sample. The highly
    connected vertices plotted in the center correspond to values
    lower right region in $[0,1]^2$. The minimum function example, due to
    \citet{Lovasz:2013:A}, is chosen as a particularly simple
    symmetric function which is not piece-wise constant. See \cref{sec:models} for
    examples with more structure.
  }
  \label{fig:graphon}
\end{figure*}
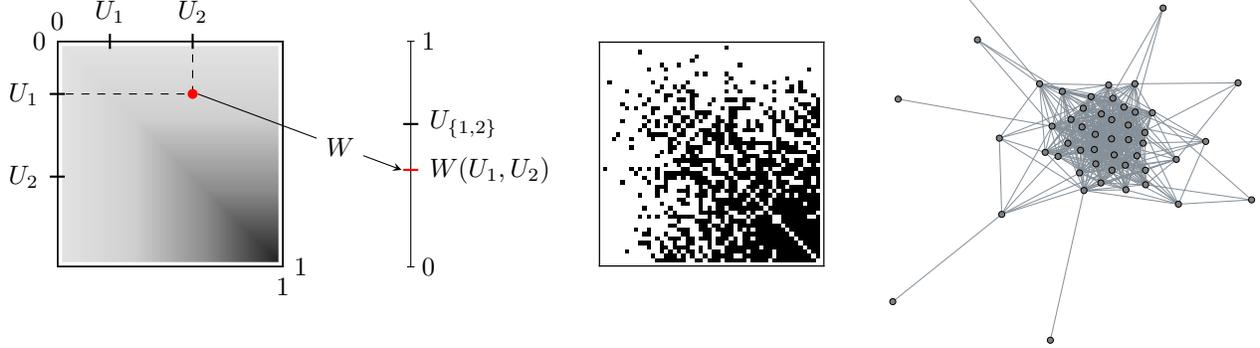

\def\Urow{U^{\mbox{\tiny row}}}
\def\Ucol{U^{\mbox{\tiny col}}}

\begin{corollary}[Aldous]%[Separately exchangeable arrays]
  \label{theorem:AH:2D:separately}
  A random array $(X_{\oset i j})$ is separately exchangeable if and only if it can be represented as follows:
  There is a random function $F:[0,1]^3\to\dataspace$ such that
  \begin{equation}
    \label{eq:AH:2D:separately}
    (X_{\oset i j})
    \eqdist 
    (F(\Urow_i,\Ucol_j,U_{\oset i j})) \;,
  \end{equation}
  where $(\Urow_i)_{i\in\Nats}$, $(\Ucol_j)$ and $(U_{\oset i j})_{i,j\in\Nats}$ are, respectively, two sequences and an array
  of i.i.d.\ $\Uniform[0,1]$ random variables, which are independent of $F$.
\end{corollary}
It is not hard to see that this follows directly from the jointly
exchangeable case: If we choose two disjoint, infinite subsets $C$ and
$R$ of $\Nats$, there exist bijections 
${r : \Nats \to R}$ and ${c : \Nats \to C}$.
Separate exchangeability of $(X_{\oset ij})$ then implies ${(X_{\oset ij}) \eqdist (X_{r_i c_j})}$.
Because a separately exchangeable array is jointly exchangeable, it can
be represented as in \cref{eq:AH:2D:jointly}, and substituting 
$(X_{r_i c_j})$, $R$ and $C$ into \eqref{eq:AH:2D:jointly} yields \eqref{eq:AH:2D:separately}.

Because separate exchangeability treats rows and columns independently, the single sequence $(U_i)$ of random variables in
\cref{eq:AH:2D:jointly} is replaced by two distinct sequences $(\Urow_i)_{i\in\Nats}$ and $(\Ucol_j)_{j\in\Nats}$,
respectively. Additionally, for each pair of distinct indices $i$ and
$j$, the single variable $U_{\nset{i,j}}$ in 
the joint case is now replaced by a pair of variables $U_{\oset ij}$
and $U_{\oset ji}$.
The index structure of the uniform random variables is the only difference between the jointly and separately 
exchangeable case.

\begin{example}[Collaborative filtering]\label{ex:collab}
  In the prototypical version of a collaborative filtering problem, users assign scores to movies. Scores may be binary
  (``like/don't like'', $X_{\oset ij}\in\lbrace 0,1\rbrace$), have a finite range (``one to five stars'', $X_{\oset ij}\in\lbrace 1,\dotsc 5\rbrace$),
  etc. Separate exchangeability then simply means that the probability of seeing any particular realization of the matrix
  does not depend on the way in which either the users or the movies are ordered.
\end{example}

Like de Finetti's and Kingman's theorem, the representation results
are complemented by a convergence result, due to 
\citet[][Theorem 3]{Kallenberg:1999}. The general result
involves some technicalities and is not stated here; the
special case for random graphs is discussed in \cref{sec:graph:limits}
(see \cref{theorem:kallenberg}).

\begin{remark}[Non-uniform sampling schemes]
  The variables $(U_i)$ and $(U_{ij})$ in the representation \eqref{eq:AH:2D:jointly} of a jointly
  exchangeable array need not be uniform: We can substitute variables with a different marginal distribution
  if we modify the random function $F$ accordingly. For \eqref{eq:AH:2D:jointly} to
  represent \emph{some} jointly exchangeable array, it is sufficient that:
  \begin{enumerate}
  \item Both $(U_i)$ and $(U_{ij})$ are \iid
  \item The variables are independent of the random function $F$.
  \end{enumerate}
  Conversely, suppose we define variables $(U_i)$ and $(U_{ij})$ satisfying 1) and 2), and
  additionally:
  \begin{enumerate}
    \setcounter{enumi}{2}
  \item Values do not re-occur, with probability one; that is, the distributions substituted for
    the uniform is atomless.
  \end{enumerate}
  Then for any jointly exchangeable array $(X_{ij})$, there is a random function $F$ such that
  the representation \eqref{eq:AH:2D:jointly} holds.
  The representation \eqref{eq:AH:2D:separately}
  of separately exchangeable arrays can be modified analogously.
  The choice of $[0,1]$ and the uniform distribution is a canonical one in probability, but one 
  could also choose $\Reals$ and any Gaussian distribution.
  The resemblance between functions  on $[0,1]^2$ and empirical graph distributions
  (see \cref{fig:graphon}) makes the unit square convenient for purposes of exposition. 
\end{remark}

\subsection{Exchangeable Graphs}

A particular important type of random structures are random graphs defined on a nonrandom vertex set.
For a graph on a countably infinite vertex set, we can consider the
vertex set to be $\Nats$ itself, without any loss of generality. A
random graph $G$ is then given by a random edge set, which is a random 
subset of ${\Nats\times\Nats}$.
A natural symmetry property of a random graph is the invariance of its distribution to a relabeling/permutation of its vertex set.  In this case, $G$ is said to be an \defn{exchangeable graph}. 
Informally, $G$ can be thought of as a random graph up to isomorphism, and so its distribution is determined by the frequency of edges, triangles, five-stars, etc., rather than by where these finite subgraphs appear.  

It is straightforward to check that $G$ is an exchangeable graph if and only if its adjacency matrix is jointly exchangeable.
More carefully, let ${(X_{ij})}$ be an array of binary random variables and put ${X_{ij}=1}$ if and only if there is an edge between vertices ${i,j \in \Nats}$ in $G$.  Then a simultaneous permutation of the rows and columns of $X$ is precisely a relabeling of the vertex set of $G$.

An important special case is when $G$ is simple---i.e., undirected and without self-loops.  In this case, $X$ is symmetric with a zero diagonal, and its representation \eqref{eq:AH:2D:jointly} 
can be simplified:  Let $F$ be a random function satisfying \eqref{eq:AH:2D:jointly}.  Without loss of generality, we may assume that $F$ is symmetric in its first two arguments.  Consider the (random) function $W$ from $[0,1]^2$ to $[0,1]$ given by $W(x,x) \defas 0$ on the diagonal and otherwise by
\[
 \Graphon(x,y) 
 &\defas \Pr [ F(x,y,U) = 1 | F] \\
 &= \int_0^1 F(x,y,u) \,\dee u, \nonumber
\]
where $U \dist \Uniform[0,1]$ is independent of $F$.
Then $W$ is random symmetric function from $[0,1]^2$ to $[0,1]$, and, by construction,
\[
  (F(U_i,U_j,U_{\nset{i,j}}))
  \eqdist
  (\Ind \lbrace U_{\nset{i,j}} < \Graphon(U_i,U_j)\rbrace)\;.
\]
We thus obtain the following specialization of the Aldous-Hoover theorem for exchangeable simple graphs:
\begin{corollary}%[Exchangeable graphs]
  Let $G$ be a random simple graph with vertex set $\Nats$ and let $X$ be its adjacency matrix.
   Then $G$ is an exchangeable graph if and only if there is a random
  function $\Graphon$ from $[0,1]^2$ to $[0,1]$
  such that
  \begin{equation}
    \label{eq:AH:graph}
    (X_{ij})
    \eqdist
    (\Ind\lbrace U_{\nset{i,j}} < \Graphon(U_i,U_j) \rbrace),
  \end{equation}
  where $U_i$ and $U_{\nset{i,j}}$ are independent i.i.d.\ uniform
  variables as in \eqref{eq:AH:2D:jointly}, which are independent of $W$.
\end{corollary}

The representation \eqref{eq:AH:graph} yields the following generative process:
\begin{enumerate}
\item 
Sample a random function $W \dist \nu$.
\item 
For every vertex $i \in \Nats$, sample an independent uniform random variable $U_i$, independent also from $W$.
\item
For every pair of vertices ${i < j \in \Nats}$, sample 
\[
  X_{ij} \given W, U_i, U_j
  \dist
  \Bernoulli(\Graphon(U_i,U_j))\;,
\]
where ${X_{ij}=1}$ indicates the edge connecting $i$ and $j$ is present;
if ${X_{ij}=0}$, it is absent.
\end{enumerate}
\cref{fig:graphon} illustrates the generation of random simple graph.

Following \citep{Lovasz:2013:A}, we call a (measurable) function from $[0,1]^2$ to $[0,1]$ a \defn{graphon}.
Thus, every exchangeable graph is 
represented by a random graphon
$\Graphon$.
In the language of integral decompositions,
 the ergodic distributions of exchangeable simple graphs
are parametrized by graphons: In \eqref{eq:integral:decomp}, we could take $\nu$ to be a distribution on the space of graphons.
We will see in \cref{sec:graph:limits} that $\theta$ has an
interpretation as a limit of the empirical adjacency matrices for larger and larger subgraphs.

\subsection{Application to Bayesian Models}

The representation results above have fundamental implications for Bayesian modeling---in fact, they provide a general
characterization of Bayesian models of array-valued data:
\begin{quote}
  {\em 
  Statistical models of exchangeable arrays can be parametrized by functions from ${[0,1]^3 \to [0,1]}$.
  Every Bayesian model of an exchangeable array is characterized by
  a prior distribution on the space of functions from
  $[0,1]^3$ to $[0,1]$.}
\end{quote}
For the special case of simple graphs, we can rephrase this idea in terms of
graphons:
\begin{quote}
  {\em 
    Statistical models of exchangeable simple graphs are parametrized by
    graphons. Every Bayesian model of
    an exchangeable simple graph is characterized by a prior distribution
    on the space of graphons.
  }
\end{quote}
In a Bayesian model of an exchangeable $2$-array,
the random function $F$ plays the role of the random parameter
$\Theta$ in \eqref{eq:hierarchy}, and the parameter space $\tspace$
is the set of measurable function ${[0,1]^3\rightarrow[0,1]}$.
Every possible value $f$ of $F$ defines an ergodic distribution 
$\kernelval \genkernel f$:
In the jointly exchangeable case, for example, 
\cref{theorem:AH:2D:jointly} shows that $\Xinf$ can be sampled
from $\kernelval \genkernel f$ by sampling
\begin{align}
  \label{eq:generative:jointly:2}
  \forall i\in\Nats:&{}
  &
  U_i&\simiid\Uniform[0,1] 
  & 
  {}&{}\\
  \label{eq:generative:jointly:3}
  \forall i,j\in\Nats:&{}
  &
  U_{\mset{i,j}}&\simiid\Uniform[0,1] 
  & 
  {}&{}
\end{align}
and computing $X$ as
\begin{align}\label{eq:sampform}
  \forall i,j\in\Nats:&{}
  &
  X_{\oset ij}& \defas f(U_i,U_j,U_{\mset{i,j}}) \;.
  & 
  {}&{}
\end{align}
Similarly, the ergodic distributions for separately exchangeable
2-arrays are given by \eqref{eq:AH:2D:separately}.
In the special case of exchangeable simple graphs, the parameter
space $\tspace$ can be reduced to the set of graphons, and the
ergodic distribution $\kernelval \genkernel w$ defined by a graphon
$w$ is given by \eqref{eq:AH:graph}.

To the best of our knowledge, \citet{Hoff:2007:1} was the first
to invoke the Aldous-Hoover theorem for statistical modeling.
The problem of estimating the distribution of an
exchangeable graph can be formulated as a regression problem on
the unknown function $w$.
This perspective was proposed in
\citep{Lloyd2012}, where the regression
problem is formulated as a Bayesian nonparametric model with a
Gaussian process prior. The regression model need not be Bayesian,
however, and recent work formulates 
the estimation of $w$ under suitable modeling conditions
as a maximum likelihood problem \citep{Wolfe:Olhede:2013:1}.

\begin{remark}[Beyond exchangeability]
  Various types of array-valued data depend on time
  or some other covariate. In this case, joint or separate
  exchangeability might be assumed to hold marginally, as
  described in \cref{sec:nonexchangeable}.
  E.g., in the case of a graph evolving over time, one could posit the existence of a graphon $\Graphon(.,.,t)$ depending also on the time $t$.
  More generally, the discussion in \ref{sec:nonexchangeable} applies to joint and separate
  exchangeability just as it does to exchangeable sequences.
  On the other hand, sometimes exchangeability will not be an appropriate assumption, even marginally.  In \cref{sec:sparsity}, we highlight some reasons why exchangeable graphs may be poor models of very large sparse graphs.
\end{remark}

\subsection{Uniqueness of representations}
\label{sec:nonuniqueness}

\begin{figure}
  \begin{center}
  \includegraphics[width=3cm]{min_function_bw.pdf}
  \qquad\qquad
  \includegraphics[width=3cm]{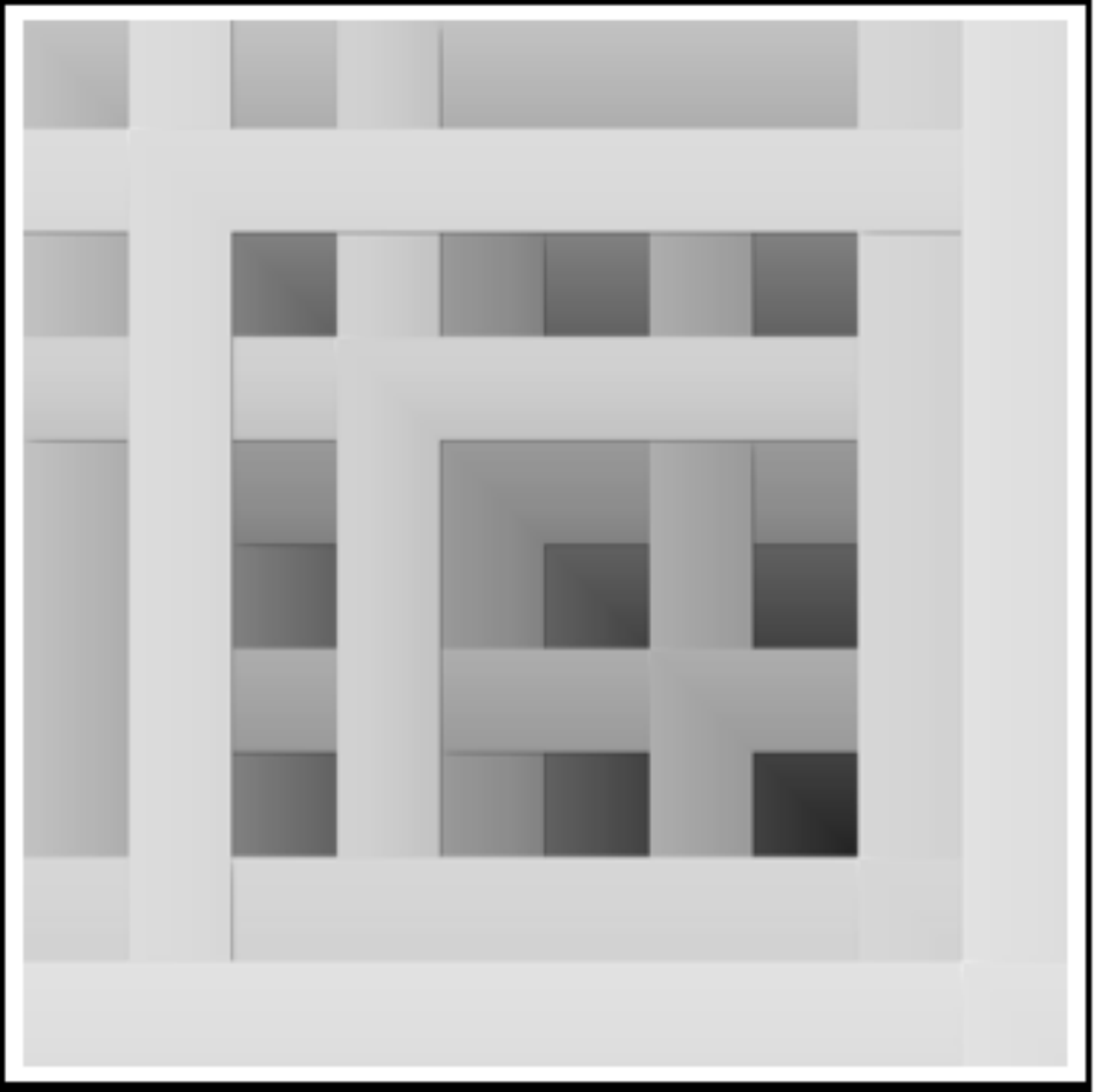}
  \end{center}
  \caption{Non-uniqueness of representations: The function on the left parametrizes
    a random graph as in \cref{fig:graphon}. On the right, this function has been modified by dividing the unit square into
    $10\times 10$ blocks and applying the same permutation of the set $\lbrace 1,\dotsc,10\rbrace$
    simultaneously to rows and columns. Since the random variables $U_i$ in \cref{eq:AH:graph}
    are i.i.d., sampling from either function defines one and the same distribution on random graphs.
  }
  \label{fig:non:uniqueness}  
\end{figure}

\begin{figure}[b]
  \caption{Two examples of graphons for which monotonization does not yield a 
    unique representation. {\it Upper row:}
    The functions $w$ and $w'$ are distinct but parametrize the same random graph.
    For both, the projection $v(x)$ (see \cref{remark:monotonization}) is the constant
    function with value $\frac{1}{2}$, which means both remain invariant and hence distinct
    under monotonization. Additionally, $w''$ also projects to ${v(x)=\frac{1}{2}}$, but
    parametrizes a different random graph, \ie projections do not
    distinguish different random graphs.
    {\it Lower row:} Another different example, where again $w$ and $w'$ are equivalent to each other, but not
    to $w''$. All three functions project to ${v(x)=\frac{1}{3}}$.}
  \begin{center}
    \begin{tikzpicture}[scale=2]
      \begin{scope}[font=\scriptsize]
        \draw (0,0)--(0,1)--(1,1)--(1,0)--(0,0);
        \foreach \x/\y in {0/.25,0/.75,.25/0,.25/.5,.5/.25,.5/.75,.75/0,.75/.5}
                 {
                   \draw[fill] ($(\x+0,\y+0)$)--($(\x+0,\y+0.25)$)--($(\x+0.25,\y+0.25)$)--($(\x+0.25,\y+0)$)--cycle;
                   \node[white] at ($(\x+.125,\y+.125)$) {$1$};
                 };
                 \foreach \x/\y in {0/0,0/.5,.25/0.25,.25/.75,.5/0,.5/.5,.75/0.25,.75/.75}
                 \node at ($(\x+.125,\y+.125)$) {$0$};
                 \node at (0.5,-.15) {\normalsize $w\phantom{'}$};
      \end{scope}
      \node at (1.15,.5) {$\equiv$};
      \begin{scope}[xshift=1.3cm]
        \draw (0,0)--(0,1)--(1,1)--(1,0)--(0,0);
        \draw[fill] (0,0.5)--(0,1)--(.5,1)--(.5,.5)--(0,0.5);
        \draw[fill] (.5,.5)--(1,.5)--(1,0)--(.5,0)--(.5,.5);
        \node at (0.25,0.25) {$0$};
        \node at (0.75,0.75) {$0$};
        \node[white] at (0.25,0.75) {$1$};
        \node[white] at (0.75,0.25) {$1$};
        \node at (0.5,-.15) {\normalsize $w'$};
      \end{scope}
      \node at (2.45,.5) {$\not\equiv$};
      \begin{scope}[xshift=2.6cm]
        \draw[fill=black!50!white] (0,0)--(0,1)--(1,1)--(1,0)--(0,0);
        \node at (0.5,0.5) {$\frac{1}{2}$};
        \node at (0.5,-.15) {\normalsize $w''$};;
      \end{scope}
      \begin{scope}[yshift=-1.45cm]
        \begin{scope}
          \draw[black!50!white,fill] (0,3/3)--(2/3,3/3)--(2/3,1/3)--(0,1/3)--cycle;
          \draw[black,fill] (2/3,1/3)--(3/3,1/3)--(3/3,0)--(2/3,0)--cycle;
          \draw (0,0)--(0,3/3)--(3/3,3/3)--(3/3,0)--cycle;
          \node at (2/6,4/6) {$\frac{1}{2}$};
          \node[white] at (5/6,1/6) {$1$};
          \node at (0.5,-.15) {\normalsize $w\phantom{'}$};
        \end{scope}
        \node at (1.15,.5) {$\equiv$};
        \begin{scope}[xshift=1.3cm]
          \draw[black!50!white,fill] (1/3,0)--(1/3,2/3)--(3/3,2/3)--(3/3,0)--cycle;
          \draw[black,fill] (0,2/3)--(0,3/3)--(1/3,3/3)--(1/3,2/3)--cycle;
          \draw (0,0)--(0,3/3)--(3/3,3/3)--(3/3,0)--cycle;
          \node[white] at (1/6,5/6) {$1$};
          \node at (2/3,1/3) {$\frac{1}{2}$};
          \node at (0.5,-.15) {\normalsize $w'$};
        \end{scope}
        \node at (2.45,.5) {$\not\equiv$};
        \begin{scope}[xshift=2.6cm]
          \draw[fill=black!33!white] (0,0)--(0,3/3)--(3/3,3/3)--(3/3,0)--cycle;
          \node at (.5,.5) {$\frac{1}{3}$};
          \node at (0.5,-.15) {\normalsize $w''$};
        \end{scope}
      \end{scope}
    \end{tikzpicture}
\end{center}
\label{fig:counterexample}
\end{figure}
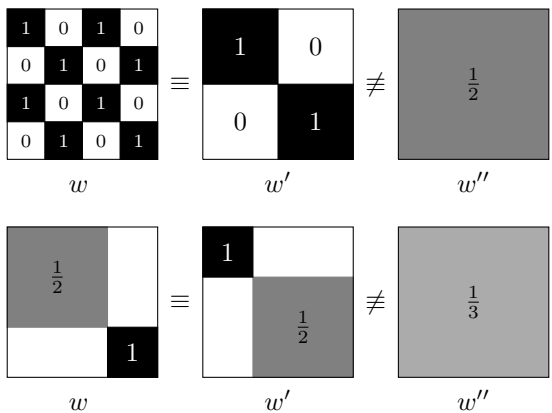

In the representation \cref{eq:AH:graph}, random graph 
distributions are parametrized by 
functions ${w:[0,1]^2\to[0,1]}$. 
This representation is not unique: Two distinct graphons may parametrize
the same random graph. In this case, the two graphons
are called \defn{weakly isomorphic} \citep{Lovasz:2013:A}.
From a statistical perspective, this means the graphon is not identifiable when regarded
as a model parameter, although it
is possible to treat the estimation problem up to equivalence of functions 
\citep[][Theorem 4]{Kallenberg:1999}.

To see that the representation by $w$ is not unique, simply note that
the graphon $w'(x,y) = w(1-x,1-y)$ is weakly isomorphic to $w$ because ${(U_i) \eqdist (V_i)}$ when ${V_i = 1-U_i}$ for $i \in \Nats$.
More generally,
let ${\phi:[0,1]\to[0,1]}$ be a \defn{measure-preserving transformation} (MPT), i.e., a map such that $\phi(U)$ is uniformly distributed when $U$ is.
By the same argument as above, the graphon $w^\phi$ given by $w^\phi(x,y) = w(\phi(x),\phi(y))$ is weakly isomorphic to $w$.
\cref{fig:non:uniqueness} shows another example of a function $w$ and its
image under a MPT.

Although any graphon $w^{\phi}$ obtained from $w$ by a MPT $\phi$ is
weakly isomorphic to $w$, the converse is \emph{not}
true: For two weakly isomorphic graphons, there need not be a MPT that
transforms one into the other \citep[see][Example 7.11]{Lovasz:2013:A}.

\begin{remark}[Monotonization does not yield canonical representations]
  \label{remark:monotonization}
  A question that often arises in this context is whether a unique
  representation can be defined through ``monotonization'': On the
  interval, every bounded real-valued function can be transformed
  into a monotone left-continuous functions by a measure-preserving transformation,
  and this left-continuous
  representation is unique \citep[e.g.][Proposition A.19]{Lovasz:2013:A}.
  It is well known in combinatorics that the
  same does \emph{not} hold on $[0,1]^2$ \citep{Borgs:Chayes:Lovasz:2010,Lovasz:2013:A}.
  More precisely, one might attempt to monotonize 
  $w$ on $[0,1]^2$ by first considering its projection 
  ${v(x):=\int_0^1 w(x,y) \,\dee y}$. The one-dimensional function
  $v$ can be transformed into a monotone representation by a unique
  MPT $\phi$, which we can then apply to both arguments of $w$ to obtain $w^\phi$. 
  Although uniqueness in the one-dimensional case suggests that 
  the resulting graphon $w^\phi$ depends only on the weak-isomorphism class of $w$,  
  this approach does not yield a canonical representation.
  \cref{fig:counterexample} shows two examples: In each example, the functions
  $w$ and $w'$ have identical projections 
  ${v=v'}$, and thus identical MPTs ${\phi=\phi'}$.  Therefore,
  we have ${w^\phi \neq {w'}^{\phi'}}$, even though $w$ and $w'$ parametrize the same random graph.
\end{remark}

\section{Models in the Machine Learning Literature}
\label{sec:models}

\begin{figure*}
\centering
\includegraphics[width=.19\linewidth]{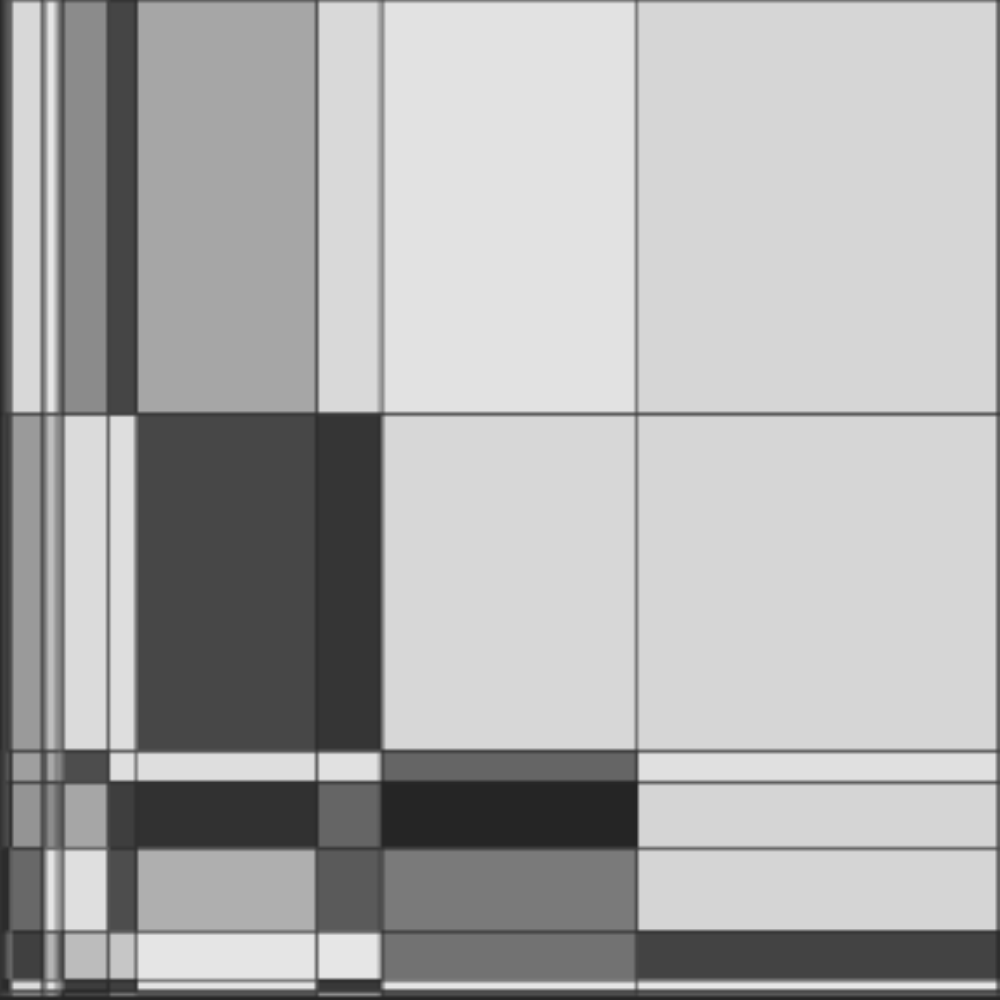}
\includegraphics[width=.19\linewidth]{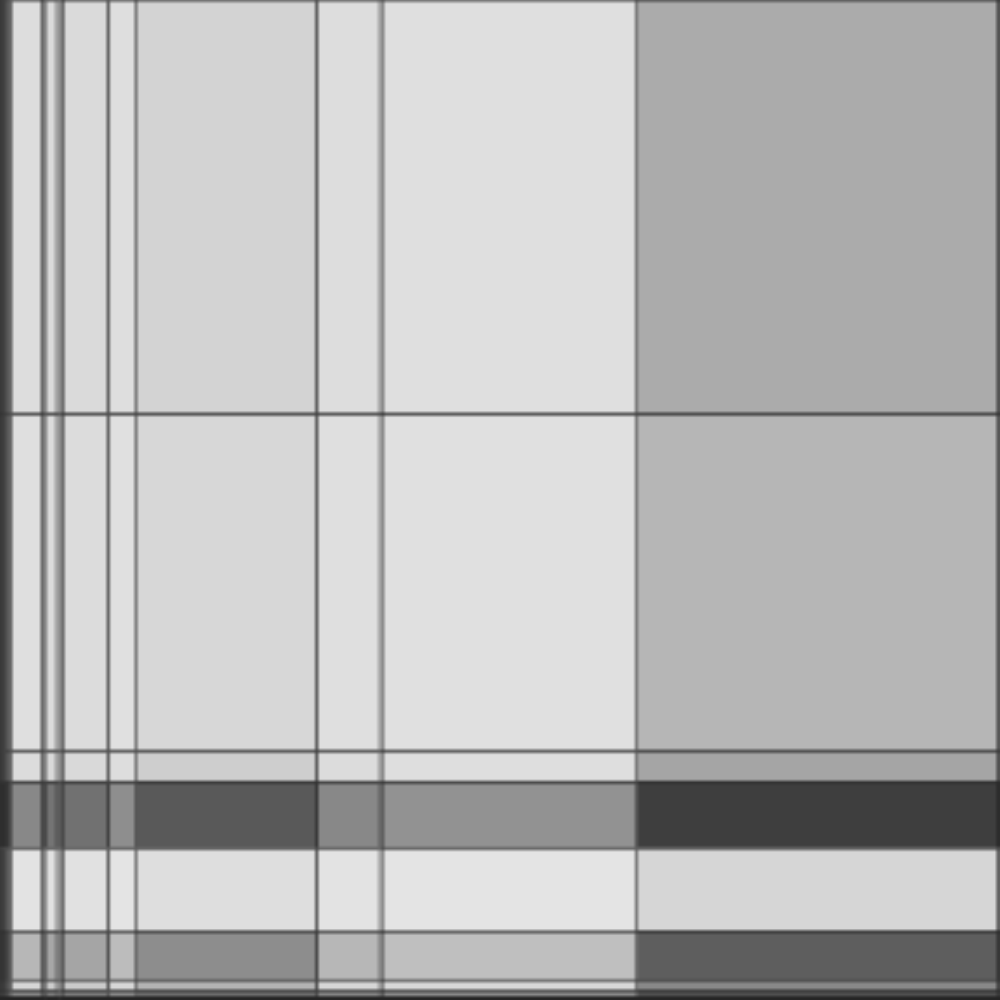}
\includegraphics[width=.19\linewidth]{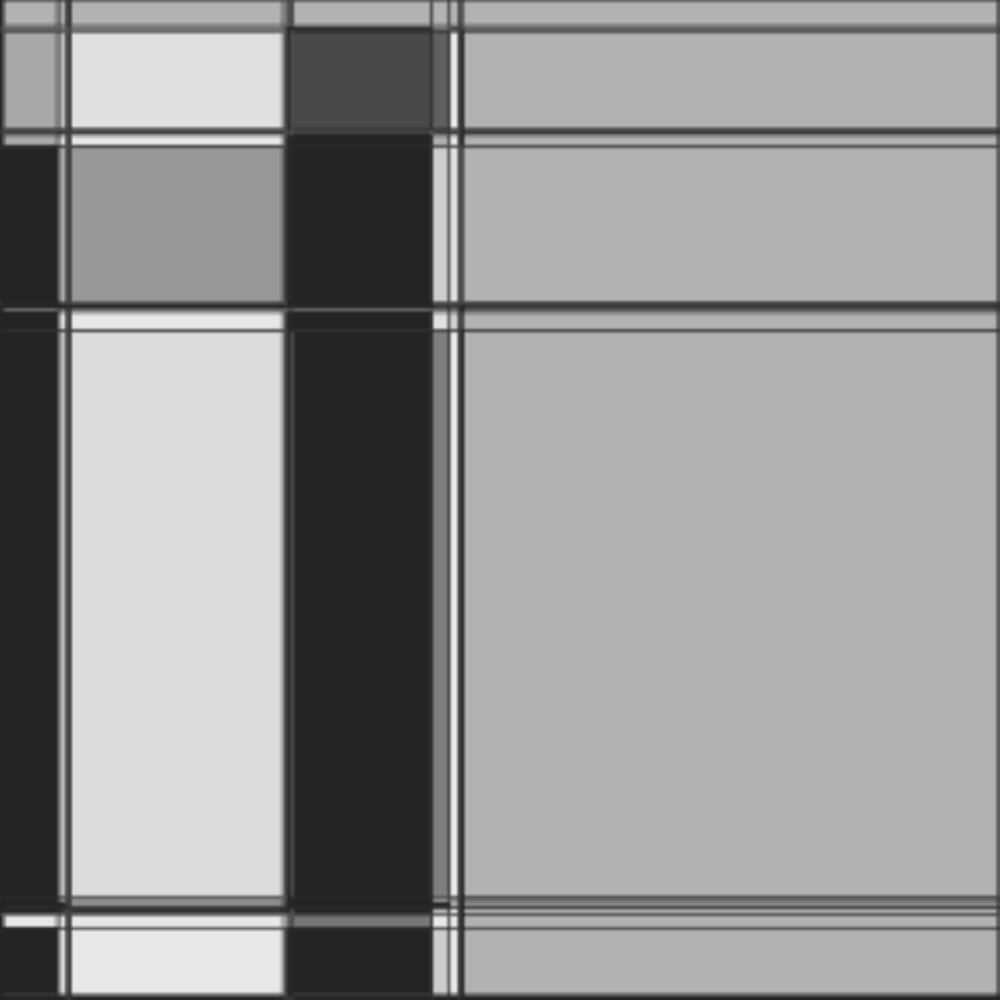}
\includegraphics[width=.19\linewidth]{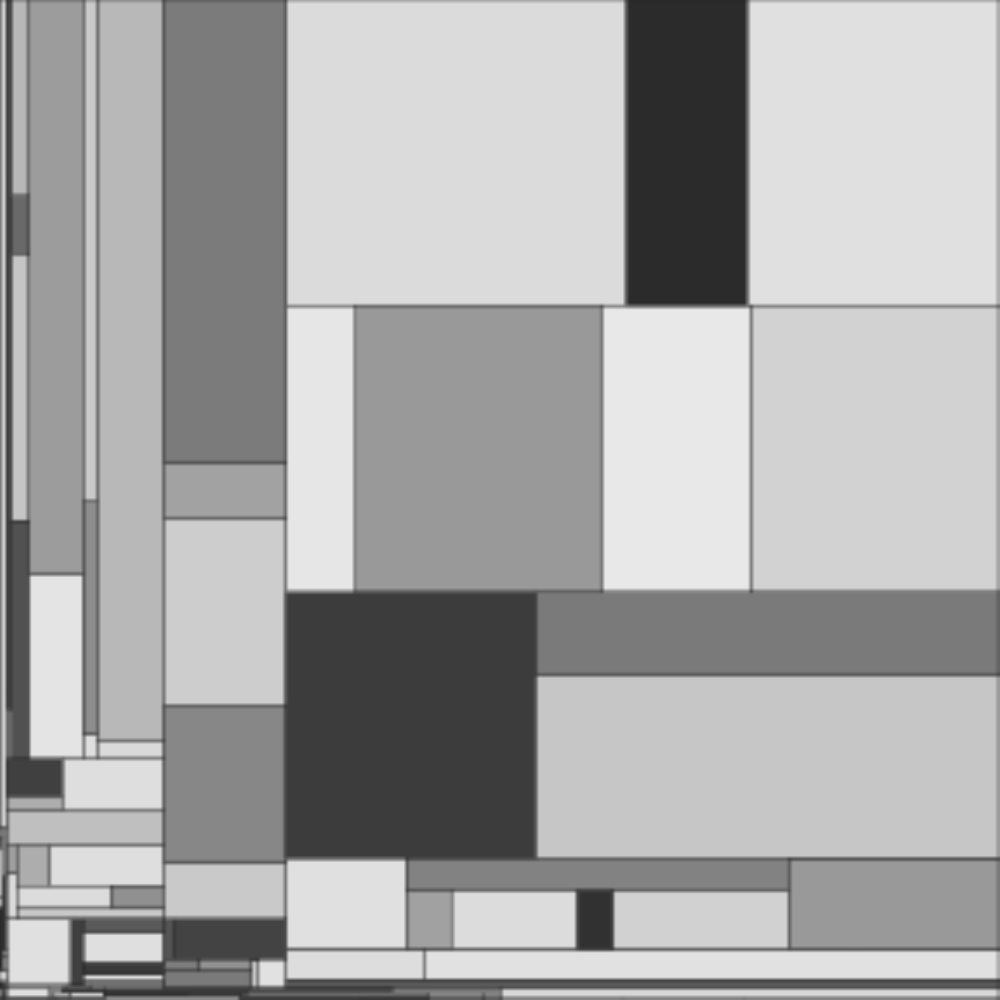}
\includegraphics[width=.19\linewidth]{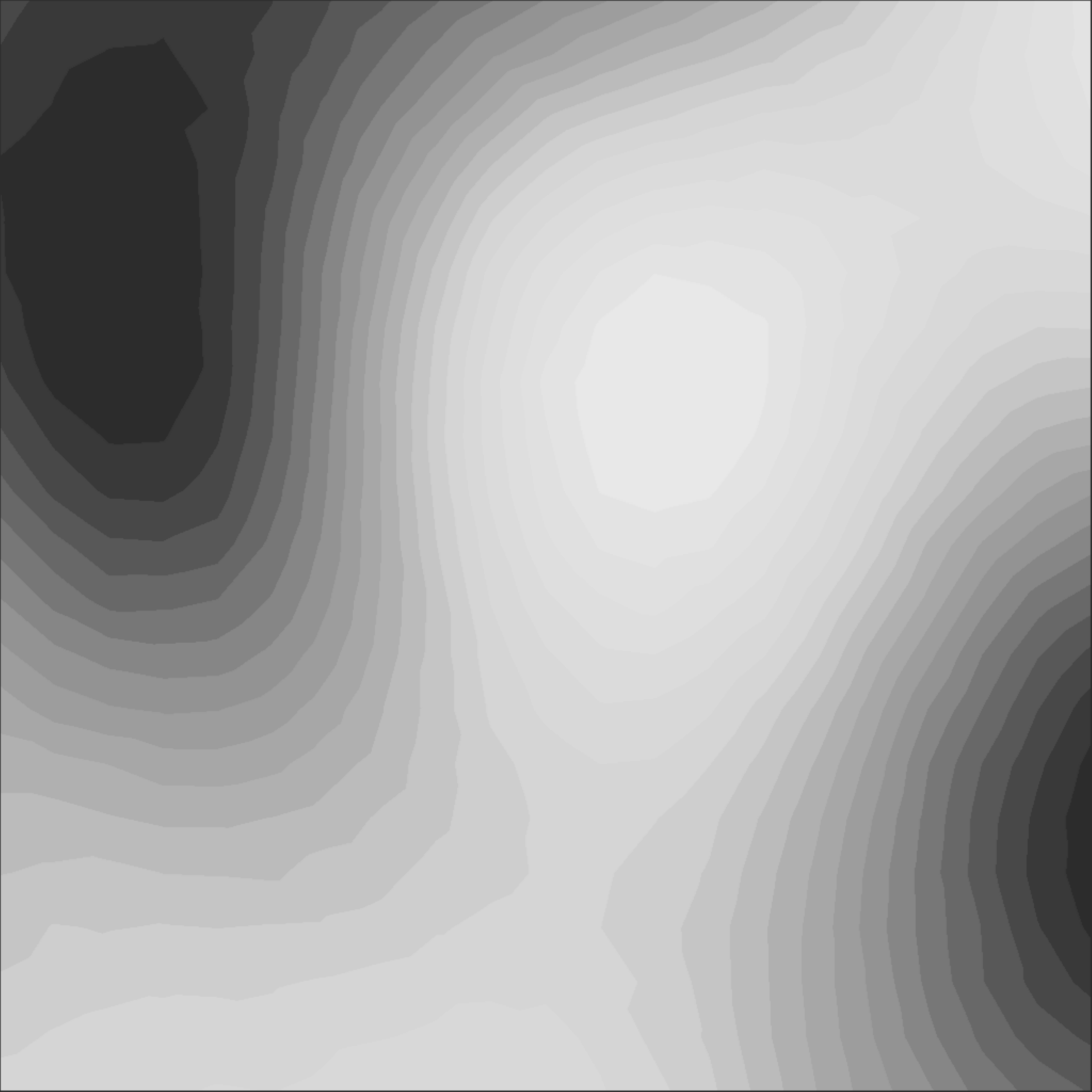}
\caption{Typical directing random functions underlying, from left to right, 1) an IRM (where partitions correspond with a Chinese restaurant process) with conditionally i.i.d.\ link probabilities; 2) a more flexible variant of the IRM with merely \emph{exchangeable} link probabilities as in~\cref{ex:elp}; 3) a LFRM (where partitions correspond with an Indian buffet process) with feature-exchangeable link probabilities as in~\cref{ex:felp}; 4) a Mondrian-process-based model with a single latent dimension; 5) a Gaussian-processed-based model with a single latent dimension.  In the first four figures, we have truncated each of the ``stick-breaking'' constructions at a finite depth.
}
\label{fig:blockstructured}
\end{figure*}

\begin{it}
  The representation theorems show that any Bayesian model of an exchangeable array can be specified by a prior on functions. 
  Models can therefore be classified according to the type of random function they employ.
  This section surveys
  several common categories of such random functions, including 
  random piece-wise constant (p.w.c.) functions, which account for the structure of models built using Chinese restaurant processes, Indian buffet processes and other combinatorial stochastic processes; 
  and random continuous functions with, e.g., Gaussian process priors. 
  Special cases of the latter include
  a range of matrix factorization and dimension reduction models proposed in the machine learning literature.
  \Cref{tab:model:classes} summarizes the classes in terms of restrictions on the random function and the  values it takes, and
  \cref{fig:blockstructured} depicts typical random functions across these classes.
\end{it}

\begin{table*}
 \centering{
    \begin{tabular}{lll}
      \toprule
      Model class & Random function $F$ & Distribution of values\\
      \midrule
      Cluster-based (\cref{sec:partition:based}) & p.w.c.\ on random product partition  & exchangeable \\
      Feature-based (\cref{sec:feature:based}) & p.w.c.\ on random product partition  & feature-exchangeable \\
      Piece-wise constant (\cref{sec:pwc:based}) & p.w.c.\ general random partition & arbitrary \\
      Gaussian process-based (\cref{sec:GP:based}) & continuous &
      Gaussian \\
      \bottomrule
    \end{tabular}
  }
 \caption{Important classes of exchangeable array models, categorized according to the type random function parametrizing the model (where p.w.c.\ stands for piece-wise constant).}
 \label{tab:model:classes}
\end{table*}

\subsection{Cluster-based models}
\label{sec:partition:based}

Cluster-based models assume that the rows and columns of the random array $X \defas (X_{\oset i j})$ can be partitioned into (disjoint) classes, such that the probabilistic structure between every row- and column-class is homogeneous.  Within social science, this idea is captured by assumptions underlying \defn{stochastic block models}~\citep{HollandLaskeyLeinhardt83,WassermanAnderson87}.

The collaborative filtering problem described in \cref{ex:collab} is a prototypical application: 
here, a cluster-based model would assume that the users can be partitioned into classes/groups/types/kinds (of users), and likewise, the movies can also be partitioned into classes/groups/types/kinds (of movies). 
Having identified the underlying partition of users and movies, each class of user would be assumed to have a prototypical preference for each class of movie.  

Because a cluster-based model is described by two partitions, the exchangeable partition models well-known from Bayesian nonparametric clustering can be used as building blocks. If the partitions are in particular generated by a 
by a Chinese restaurant process, we obtain the Infinite Relational Model (IRM), introduced in \citep{Kemp2006} and independently in \citep{Xu2006}.
The IRM can be seen as a nonparametric generalization of parametric stochastic block models \citep{HollandLaskeyLeinhardt83,WassermanAnderson87}.  In the following example, we describe the model for the special case of a $\{0,1\}$-valued array.

\begin{example}[Infinite Relational Model]
\label{example:IRM}
Under the IRM, the generative process for a finite subarray of binary random variables $X_{ij}$, $i \le n$, $j \le m$, is as follows:  To begin, we partition the rows (and then columns) into clusters according to a Chinese restaurant process: The first and second row are chosen to belong to the same cluster with probability proportional to 1 and to belong to different clusters with probability proportional to a parameter $c>0$.  Subsequently, each row is chosen to belong to an existing cluster with probability proportional to the current size of the cluster, and to a new cluster with probability proportional to $c$.  Let $\Pi \defas \{\Pi_1,\dotsc,\Pi_\kappa\}$ be the random partition of $\{1,\dotsc,n\}$ induced by this process, where $\Pi_1$ is the cluster containing 1, and $\Pi_2$ is the cluster containing the first row not belonging to $\Pi_1$, and so on.  Note that the number of clusters, $\kappa$, is also a random variable.  Let $\Pi' \defas \{\Pi'_1,\dotsc,\Pi'_{\kappa'}\}$ be the random partition of $\{1,\dotsc,m\}$ induced by this process on the \emph{columns}, possibly with a different parameter $c' >0$ determining the probability of creating new clusters.  Next, for every pair $(k,k')$ of cluster indices, $k \le \kappa$, $k' \le \kappa'$, we generate an independent beta random variable $\theta_{k,k'}$.
Finally, we generate each $X_{ij}$ independently from a Bernoulli distribution with mean $\theta_{k,k'}$, where $i \in \Pi_k$ and $j \in \Pi'_{k'}$.  As we can see, $\theta_{k,k'}$ represents the probability of links arising between elements in clusters $k$ and $k'$.

The Chinese restaurant process generating $\Pi$ and $\Pi'$ is known to be exchangeable in the sense that the distribution of $\Pi$ is invariant to a permutation of the underlying set $\{1, \dotsc, n\}$.  It is then straightforward to see that the distribution on the subarray is exchangeable.  In addition, it is straightforward to verify that, were we to have generated an $n+1 \times m+1$ array, the marginal distribution on the $n \times m$ subarray would have agreed with that of the above process.  This implies that we have defined a so-called projective family and so results from probability theory imply that there exists an infinite array and that the above process describes the distribution of every finite subarray.
\end{example}

\begin{definition}%[simple cluster-based models]
\label{def:simple:cluster:based}
  We say that a Bayesian model
  of an exchangeable array is \defn{simple cluster-based} 
  if, for some random function $F$ representing $X$,
  there are random partitions $B_1,B_2,\dotsc$ and $C_1,C_2,\dotsc$ 
  %and, for every $i,j\in \Nats$, $D^{i,j}_1,D^{(i,j)}_2,\dotsc$ 
  of the unit interval $[0,1]$ such that:
  \begin{enumerate}
  \item On each block 
      ${A_{i,j} \defas B_i\times C_j\times [0,1]}$, $F$ is constant.  Let $f_{\oset ij}$ be the value $F$ takes on block $A_{i,j}$.
  \item The block values $(f_{\oset ij})$ are themselves an exchangeable array, and independent from $(B_i)$ and $(C_j)$.
  \end{enumerate}
  We call an array simple cluster-based if its distribution is.
\end{definition}
Most examples of simple cluster-based models in the literature---including, e.g., the IRM---take the block values $f_{\oset i j}$ to be conditionally i.i.d. (and so the array $(f_{\oset ij})$ is then trivially exchangeable).  As an example of a more flexible model for $(f_{\oset ij})$, which is merely exchangeable, consider the following:
\begin{example}[exchangeable link probabilities]\label{ex:elp}
  For every block $i$ in the row partition, let $u_i$ be an independent and identically distributed Gaussian random variable.  Similarly, let $(v_j)$ be an i.i.d.\ sequence of Gaussian random variables for the column partitions.  Then, for every row and column block $i,j$, put $f_{ij} \defas \textrm{sig}(u_i + v_j)$, 
  where $\textrm{sig}\colon \Reals \to [0,1]$ is a sigmoid function.  The array $(f_{ij})$ is obviously exchangeable.
\end{example}
Like with cluster-based models of exchangeable sequences, if the number of classes in each partition is bounded, then a simple cluster-based model of an exchangeable array is a mixture of a finite-dimensional family of ergodic distributions.  
Therefore, mixtures of an infinite-dimensional family must place positive mass on partitions with arbitrarily many classes.\footnote{Exchangeable partitions may contain blocks consisting of only a single element (see \cref{sec:kingman}); these are also known as \emph{dust}. Our definitions exclude this case because it complicates presentation, but the generalization is straightforward.}

%%% Snip %%%

In order to define the more general class of cluster-based models, we relax the piece-wise constant nature of the random function $F$. We do so by starting with a simple cluster-based array ${\Theta=(\Theta_{ij})}$ with values in a space $\tspace$. Since $\Theta$ is exchangeable, it is represented by a random function $G$, and we define $F$ as
\[
F(u,u',u'') = \phi(G(u,u',u''),u'')
\]
for some function ${\phi :\tspace\times [0,1] \to \dataspace}$.

If $\theta\in\tspace$ is a fixed parameter value and $U$ a uniform random variable in $[0,1]$, then 
$\phi(\theta,U)$ is a random variable with values in $\dataspace$. Suppose $P_{\theta}$ is the distribution
of $\phi(\theta,U)$. Then the array ${X=(X_{ij})}$ represented by $F$ can be sampled by sampling ${(\Theta_{ij})}$
from $G$, and then generating $X_{ij}$ from the distribution $P_{\Theta_{ij}}$.
In other words, we can posit a family ${P=\lbrace P_{\theta} | \theta\in\tspace\rbrace}$ of distributions, 
generate an array from simple cluster-based model, and then generate $X_{ij}$, conditionally independently given
$\Theta$, from $P_{\Theta_{ij}}$. This is additional layer of randomness is
also called a \emph{randomization} of $\Theta$ \citep{Kallenberg:2001}.
\begin{definition}[randomization]
  Let $P$ be a family $\kernelfamily P \theta \tspace$ of distributions 
  on $\dataspace$,  
  and let ${\Theta \defas (\Theta_{i} \st i \in I)}$ be a collection of
  random variables taking values in $\tspace$, indexed by elements of a set $I$. 
  We say that a collection ${X \defas (X_i \st i \in I)}$ of random variables, indexed by the same set $I$, is a 
  \defn{$P$-randomization of $\Theta$} when the elements $X_{i}$ are conditionally independent 
  given $\Theta$, and
  ${X_i \given \Theta \dist \kernelval P {\Theta_i}}$ for all $i \in I$.
\end{definition}
It is straightforward to show that,  if $\Theta$ is an exchangeable array (\ie ${I=\Nats^2}$) and $X$ is a randomization of $\Theta$, then $X$ is exchangeable.  We may therefore define:

\begin{definition}[cluster-based models]
  \label{def:pbm}
  We say that a Bayesian model for an exchangeable array ${X \defas (X_{\oset i j})}$ 
  in $\dataspace$ is \defn{cluster-based}
  if $X$ is a $P$-randomization of a simple cluster-based exchangeable array 
  $\Theta \defas (\Theta_{\oset ij})$ taking values in a space $\tspace$,
  for some family $\kernelfamily P \theta \tspace$ of distributions on $\dataspace$.
  We say an array is cluster-based when its distribution is.
\end{definition}

The intuition is that the cluster memberships of every pair $i,j$ of individuals determines a 
parameter $\theta_{\oset i j}$, which in turn determines a distribution $\kernelval P {\theta_{\oset i j}}$.
The actual observed relationship $X_{\oset ij}$ is then a sample from $\kernelval P {\theta_{\oset i j}}$.

%%% Snap %%%

\begin{example}[Infinite Relational Model continued]
We may alternatively describe the IRM distribution on exchangeable arrays as follows:  Let $P$ be a family $\kernelfamily P \theta \tspace$ of distributions on $\dataspace$ (e.g., a family of Bernoulli distributions indexed by their means in $[0,1]$) and let $H$ be a prior distribution on the parameter space $[0,1]$ (e.g., a Beta distribution, so as to achieve conjugacy).  The IRM is a cluster-based model, and an array ${X\defas (X_{\oset ij})}$ distributed according to an IRM is hence a $P$-randomization of a simple cluster-based array ${\Theta \defas (\Theta_{\oset ij})}$ of parameters in $\tspace$.

In order to describe the structure of $\Theta$, it suffices to
describe the distribution of the partitions $(B_k)$ and $(C_k)$
as well as that of the block values. For the latter, the IRM simply chooses the block values to be i.i.d.\ draws from the distribution $H$.  (While the block values can be taken to be merely exchangeable, we have not seen this generalization in the literature.)  For the partitions, the IRM utilizes the stick-breaking construction of a Dirichlet process~\citep{MR1309433}.

In particular, let $W_1,W_2,\dotsc$ be an i.i.d.\ sequence of $\Beta(1,\alpha)$ random variables, for some concentration parameter $\alpha > 0$.  For every $k \in \Nats$, we then define
\[\label{crpsbp}
V_k \defas (1-W_1) \dotsm (1-W_{k-1}) W_k.
\]
With probability one, we have $V_k \ge 0$ for every $k \in \Nats$ and $\sum_{k=1}^\infty V_k = 1$ almost surely, and so the sequence $(V_k)$ characterizes a (random) probability distribution on $\Nats$. We then let $(B_k)$ be a sequence of contiguous intervals that form a partition of $[0,1]$, where $B_k$ is the half-open interval of length $V_k$.
In the jointly exchangeable case, the random partition $(C_k)$ is usually chosen either as a copy of $(B_k)$, or as partition sampled independently from the same distribution as $(B_k)$.
The partitions define a partition of $[0,1]^2$ into rectangular patches, on each of which the random function
$G$ is constant.
\end{example}
The IRM is originally defined in terms of a Chinese restaurant process 
(rather than a Dirichlet process), as we have done in \cref{example:IRM}.
In terms of the CRP, each random probability $V_k$ in \eqref{crpsbp}
is the limiting fraction of rows in the $k$th cluster $\Pi_k$ as the number of rows tends to infinity.

\subsection{Feature-based models}
\label{sec:feature:based}
\newcommand{\FA}[2]{#1_{(#2)}}

Feature-based models of exchangeable arrays have similar structure to cluster-based models.  Like cluster-based models, feature-based models partition the rows and columns into clusters, but unlike cluster-based models, feature-based models allow the rows and columns to belong to multiple clusters simultaneously.  The set of clusters that a row belongs to are then called its \defn{features}.  The interaction between row $i$ and column $j$ is then determined by the features that the row and column possess. 

The stochastic process at the heart of most existing feature-based models of exchangeable arrays is the Indian buffet process, introduced by Griffiths and Ghahramani \citep{GG06}.
The Indian buffet process (IBP) produces an allocation of features in a sequential fashion, much like the Chinese restaurant process produces a partition in a sequential fashion.  In the follow example, we will describe 
the Latent Feature Relational Model (LFRM) of Miller et al.~\citep{Miller2009},
one of the first nonparametric, feature-based models of exchangeable arrays.  For simplicity, consider the special case of a $\{0,1\}$-valued, separately-exchangeable array.
\begin{example}[Latent Feature Relational Model]
\label{example:LFRM}
Under the LFRM, the generative process for a finite subarray of binary random variables $X_{ij}$, $i \le n$, $j \le m$, is as follows: To begin, we allocate features to the rows (and then columns) according to an IBP.  In particular, the first row is allocated a Poisson number of features, with mean $\gamma > 0$.  Each subsequent row will, in general, share some features with earlier rows, and possess some features not possessed by any earlier row.  Specifically, the second row is also allocated a Poisson number of altogether new features, but with mean $\gamma/2$, and, for every feature possessed by the first row, the second row is allocated that feature, independently, with probability $1/2$.  In general, the $k$th row is allocated a Poisson number of altogether new features, with mean $\gamma/k$; and, for every subset $K \subseteq \{1,\dotsc,k-1\}$ of the previous rows, and every feature possessed by exactly those rows in $K$, is allocated that feature, independently, with probability $|K|/n$.  (We use the same process to allocate a distinct set of features to the $m$ columns, though potentially with a different constant $\gamma' >0$ governing the overall number of features.)

The observed array is generated as follows: We enumerate the row- and column- features arbitrarily, and for every row $i$ and column $j$, we let $N_i, M_j \subseteq \Nats$ be the set of features they possess, respectively.
For every pair $(k,k')$ of a row- and column- feature, we generate an independent and identically distributed Gaussian random variable $w_{k,k'}$.  Finally, we generate each $X_{i,j}$ independently from a Bernoulli distribution with mean $\textrm{sig}(\sum_{k\in N_i} \sum_{k' \in M_j} w_{k,k'})$.  Thus a row and column that possess feature $k$ and $k'$, respectively, have an increased probability of a connection as $w_{k,k'}$ becomes large and positive, and a decreased probability as $w_{k,k'}$ becomes large and negative.

The exchangeability of the subarray follows from the exchangeability of the IBP itself.  In particular, define the family of counts $\Pi_N$, $N \subseteq \{1,\dotsc,n\}$, where $\Pi_N$ is the number of features possessed by exactly those rows in $N$.  We say that $\Pi \defas (\Pi_N)$ is a \defn{random feature allocation} for $\{1,\dotsc,n\}$. (Let $\Pi'$ be the random feature allocation for the columns induced by the IBP.)
The IBP is exchangeable is the sense that 
\[
(\Pi_N) \eqdist (\Pi_{\sigma(N)})
\]
for every permutation $\pi$ of $\{1,\dots,n\}$, where $\sigma(N) \defas {\{ \sigma(n) \st n \in N \}}$.  Moreover, the conditional distribution of the subarray given the feature assignments $(N_i,M_j)$ is the same as the conditional distribution given the feature allocations $(\Pi_N,\Pi'_M)$.  It is then straightforward to verify that the subarray is itself exchangeable.  Like with the IRM example, the family of distributions on subarrays of different sizes is projective, and so there exists an infinite array and the above process describes the distribution of every subarray.
\end{example}
The LFRM is a special case of a class of models that we will call \emph{feature-based}.
From the perspective of simple cluster-based models, simple feature-based models also have a block structured representing function, but relax the assumption that values of each block form an exchangeable array.

To state the definition of this class more formally, we begin by generalizing the notion of a partition of $[0,1]$.
\footnote{
  \cref{def:feature:paintbox} was also introduced and studied in much more detail
  by \citet{Broderick:Pitman:Jordan:2013:A,Broderick:Jordan:Pitman:2013:B}---which 
  we were unaware of when this article was first submitted---and we have adjusted our terminology to match that
  of \citep{Broderick:Pitman:Jordan:2013:A}, where the term \emph{feature paint-box} was coined.
  The perspectives differ slightly: In \citep{Broderick:Pitman:Jordan:2013:A}, the feature paint-box
  is introduced based on the theory
  of exchangeable partitions, and the authors show that key aspects of this theory 
  generalize elegantly. Our perspective here is that feature models form a special type of
  exchangeable arrays---as indeed do exchangeable partitions, see \citep[][Theorem 7.38]{Kallenberg:2005}.
}
  
\begin{definition}[feature paint-box]
\label{def:feature:paintbox}
Let $U$ be a uniformly-distributed random variable and 
${E\defas(E_1, E_2,\dotsc)}$ a sequence of subsets of $[0,1]$. 
Given $E$, we say that $U$ has feature $n$ if $U \in E_n$. We call the sequence
$E$ a \defn{feature paint-box} if 
\[\textstyle
\label{eq:def:feature:paint-box}
\Pr \left \{ U \notin \bigcup_{k \ge n} E_k  \right \} \to 1 \quad \text{as} \quad n \to \infty.
\]
\end{definition}
To parse \eqref{eq:def:feature:paint-box}, first recall that
a partition is a special case of a feature paint-box. In this case, the sets
$E_n$ are disjoint and represent blocks of a partition. The relation $U\in E_k$ then indicates
that an object represented by the random variable $U$ is in block $k$ of the partition.
In a feature paint-box, the sets $E_k$ may overlap. The relation $U\in E_n$ now indicates that the object has feature
$n$. Because the sets may overlap, the object may possess multiple features. However, condition \cref{eq:def:feature:paint-box} ensures
that the number of features per object remains finite (with probability 1).  

A feature paint-box induces a partition if we equate any two objects that possess exactly
the same features.  More carefully, for every subset $N\subset\Nats$ of features, define
\[
\FA E N \defas 
\bigcap_{i \in N} E_i \cap \bigcap_{j \notin N} ([0,1] \setminus E_j) \;.
\]
Then, two objects represented by random variables $U$ and $U'$ are equivalent iff $U,U' \in \FA E N$ for some finite set $N\subset\Nats$.  As before, we could consider a simple, cluster-based representing function where the block values are given by an $(f_{N,M})$, indexed now by finite subsets $N,M \subseteq \Nats$.  Then $f_{N,M}$ would determine how two objects relate when they possess features $N$ and $M$, respectively.

However, if we want to capture the idea that the relationships between objects depend on the individual features the objects possess, we would not want to assume that the entries of $f_{N,M}$ formed an exchangeable array, as in the case of a simple, cluster-based model.
E.g., we might choose to induce more dependence between $f_{N,M}$ and $f_{N',M}$ when $N\cap N' \neq \emptyset$ than otherwise.  The following definition captures the appropriate relaxation of exchangeability:
\begin{definition}[feature-exchangeable array]
Let $Y_{N,M}$ be random variables, each indexed by a pair ${N,M \subseteq \Nats}$ of finite sets, and
consider the array ${Y\defas (Y_{N,M})}$.
If $\pi$ is a permutation of $\Nats$ and ${N \subseteq \Nats}$, denote the image set as
 ${\pi(N) \defas \{ \pi(n) \st n \in N\}}$. We say that $Y$ is \defn{feature-exchangeable} if
\[
(Y_{N,M}) \eqdist (Y_{\pi(N),\pi(M)}),
\]
for all permutations $\pi$ of $\Nats$.
\end{definition}
Informally, an array $Y$ indexed by sets of features is feature-exchangeable if its distribution is invariant to permutations of the underlying feature labels (i.e., of $\Nats$).
Here is a simple example:
\begin{example}[feature-exchangeable link probabilities]\label{ex:felp}
Let $w \defas (w_{\oset ij})$ be a conditionally i.i.d.~array of random variables in $\Reals$, 
and define ${\Theta \defas (\Theta_{N,M})}$ by
\[\textstyle
\Theta_{N,M} = \textrm{sig}( \sum_{i \in N} \sum_{j \in M} w_{\oset ij}),
\]
where $\textrm{sig}\colon \Reals \to [0,1]$ maps real values to probabilities via, e.g., the sigmoid or probit functions.  It is straightforward to verify that $\Theta$ is feature-exchangeable.
\end{example}

\begin{definition}%[simple feature-based models]
  We say that a Bayesian model of an exchangeable array $X$ is \defn{simple feature-based} when,
  for some random function $F$ representing $X$,
  there are random feature allocations $B$ and $C$
  of the unit interval $[0,1]$
  such that, for every pair $N,M \subseteq \Nats$ of finite subsets,
   $F$ takes the constant value $f_{N,M}$ on the block
  \[  
    A_{N,M} \defas \FA B N \times \FA C M \times [0,1],
  \]
  and the values $f \defas (f_{N,M})$ themselves form a feature-exchangeable array, 
  independent of $B$ and $C$.  We say an array is simple feature-based if its distribution is.
\end{definition}

Simple feature-based arrays specialize to simple cluster-based arrays if either i) the feature allocations are partitions or ii) the array $f$ is exchangeable.  The latter case highlights the fact that feature-based arrays relax the exchangeability assumption of the underlying block values.  

As in the case of simple cluster-based models,
nonparametric simple feature-based models will place positive mass on feature allocations with an arbitrary number of distinct sets.
As for cluster-based models, we define general feature-based models as randomizations of simple models:

\begin{definition}[feature-based models]\label{def:fbm}
  We say that a Bayesian model for an exchangeable array 
  $X \defas (X_{\oset i j})$ 
  in $\dataspace$ is \defn{feature-based}
   when $X$ is a $P$-randomization of a simple, feature-based, exchangeable array $\Theta \defas (\Theta_{\oset ij})$ taking values in a space $T$,
  for some family $\kernelfamily P \theta \tspace$ of distributions on $\dataspace$.
  We say an array is feature-based when its distribution is.
\end{definition}
Comparing \cref{def:pbm,def:fbm}, we see that the relationship between random functions representing $\Theta$ and $X$ are the same as in cluster-based models. The LFRM described in \cref{example:LFRM} is a special case of a feature-based model:

\begin{example}[Latent Feature Relational Model continued]
A feature-based model is determined by 
the randomization family $P$,
the distribution of the underlying feature-exchangeable array $f$ of link probabilities,
and the distribution of the random feature allocation.
In the case of the LFRM, 
$P$ is the family $\Bernoulli(p)$ distributions, for $p \in [0,1]$ (although this is easily generalized, and does not represent an important aspect of the model).
The underlying feature-exchangeable array $f$ is that described in \cref{ex:felp}.  
 
The random feature allocations underlying the LFRM can be described in terms of so-called ``stick-breaking'' constructions of the Indian buffet process.  
One of the simplest stick-breaking constructions, and the one we will use here, is due to Teh, G\"or\"ur, and Ghahramani \citep{TGG07}.  (See also \citep{Thibaux2007}, \citep{PZWGC2010} and \citep{PBJ2012}.)

Let $W_1,W_2,\dotsc$ be an i.i.d.~sequence of $\Beta(\alpha,1)$ random variables for some concentration parameter $\alpha >0$.
For every $k$, we define
$%\label{ibpsbp}
\textstyle
V_k \defas \prod_{j =1}^k W_j.
$
(The relationship between this construction and \cref{crpsbp} highlights one of several relationships between the IBP and CRP.) 
It follows that we have $1 \ge V_1 \ge V_2 \ge \dotsm \ge 0$.  
The allocation of features then proceeds as follows: for every $k \in \Nats$, we assign the feature with probability $V_k$, independently of all other features.
It can be shown that  $\sum_k V_k$ is finite with probability one, and so we have a valid feature paint-box; every object has a finite number of features with probability one.  

We can describe a feature paint-box $(B_n)$ corresponding with this stick-breaking construction of the IBP as follows:  Put $B_1 = [0,V_1)$, and then inductively, for every $n\in \Nats$, put
\[
B_{n+1} 
\defas 
\bigcup_{j=1}^{2^n-1} [b_j, (b_{j+1}-b_j)\cdot V_{n+1} ) 
\]
where $B_n = [b_1,b_2) \cup [b_3,b_4) \cup \dotsm \cup [b_{2^n-1},b_{2^n})$.  (As one can see, this representation obscures the conditional independence inherent in the feature allocation induced by the IBP.)
\end{example}

\subsection{Piece-wise constant models}
\label{sec:pwc:based}

Simple partition- and feature-based models have piece-wise constant structure, which arises because both types of models posit prototypical relationships on the basis of a \emph{discrete} set of classes or features assignments, respectively.  More concretely, a partition of $[0,1]^3$ is induced by partitions of $[0,1]$.  

An alternative approach is to consider partitions of $[0,1]^3$ directly, or partitions of $[0,1]^3$ induced by partitions of $[0,1]^2$.
Rather than attempting a definition capturing a large, natural class of such models, we present an illustrative example:

\begin{example}[Mondrian-process-based models~\citep{RT09}]
A Mondrian process is a partition-valued stochastic process introduced by Roy and Teh~\citep{RT09}.  (See also Roy~\citep[][Chp.~V]{RoyThesis} for a formal treatment.)  More specifically, a \defn{homogeneous Mondrian process on $[0,1]^2$} is a continuous-time Markov chain $(M_t\colon t \ge 0)$, where,  for every time $t \ge 0$, $M_t$ is a floorplan-partition of $[0,1]^2$---i.e., a partition of $[0,1]^2$ comprised of axis-aligned rectangles of the form $A = B \times C$, for intervals $B, C \subseteq [0,1]$.  It is assumed that $M_0$ is the trivial partition containing a single rectangle.

Every continuous-time Markov chain is characterized by 
the mean waiting times between jumps
and the discrete-time Markov process of jumps (i.e., the \emph{jump chain}) embedded in the continuous-time chain.  In the case of a Mondrian process, the mean waiting time from a partition composed of a finite set of rectangles $\{B_1\times C_1,\dotsc, B_k\times C_k \}$ is $\sum_{j=1}^k (|B_j| + |C_j|)^{-1}$.  
The jump chain of the Mondrian process is entirely characterized by its transition probability kernel, which is defined as follows:
From a partition $\{B_1\times C_1,\dotsc, B_k\times C_k \}$ of $[0,1]^2$, we choose to ``cut'' exactly one rectangle, say $B_j \times C_j$, with probability proportional to $|B_j| + |C_j|$;  Choosing $j$, we then cut the rectangle vertically with probability proportional to $|C_j|$ and horizontally with probability proportional to $|B_j|$;  Assuming the cut is horizontal, we partition $B_j$ into two intervals $B_{j,1}$ and $B_{j,2}$, uniformly at random;  The jump chain then transitions to the partition where $B_j \times C_j$ is replaced by $B_{j,1} \times C_j$ and $B_{j,2} \times C_j$;  The analogous transformation occurs in the vertical case.

As is plain to see, each partition is produced by a sequence of cuts that hierarchically partition the space.  The types of floorplan partitions of this form are called \defn{guillotine partitions}.  Guillotine partitions are precisely the partitions represented by $k$d-trees, the classical data structure used to represent hierarchical, axis-aligned partitions.

The Mondrian process possesses several invariances that allow one to define a Mondrian process $M^*_t$ on all of $\Reals^2$.  The resulting process is no longer a continuous-time Markov chain.  In particular, for all $t> 0$, $M^*_t$ has a countably infinite number of rectangles with probability one.  Roy and Teh~\citep{RT09} use this extended process to produce a nonparametric prior on random functions as follows:

Let $(\psi_n)$ be an exchangeable sequence of random variables in $\dataspace$,
let $M$ be a Mondrian process on $\Reals^2$, independent of $(\psi_n)$, 
and let $(A_n)$ be the countable set of rectangles comprising the partition of $\Reals^2$ given by $M_{c}$ for some constant $c>0$.  
Roy and Teh propose the random function $F$ from $[0,1]^3$ to $[0,1]$ given by
$F(x,y,z) = \psi_{n}$, where $n$ is such that $A_n \ni (- \log x, - \log y)$.
An interesting property of $F$ is that the partition structure along any axis-aligned slice of the random function agrees with the stick-breaking construction of the Dirichlet process, presented in the IRM example. 
Roy and Teh present results in the case where the $\psi_n$ are Beta random variables, and the data are modeled as a Bernoulli randomization of an array generated from $F$.   (See \citep{RoyThesis} and \citep{RT09} for more details.)
\end{example}

Another, very popular example of a piece-wise constant model is the mixed-membership stochastic blockmodel of 
\citet{Airoldi:Blei:Fienberg:Xing:2008}.

\subsection{Gaussian-process-based models}
\label{sec:GP:based}

All models discussed so far are characterized by a random function $F$
with piece-wise constant structure.  In this section, we briefly discuss a large and important class of models that relax this restriction by modeling the random function as a Gaussian process.  

Recall that a Gaussian process \citep[e.g.][]{CarlsBook} is a distribution on random functions:
Let ${G \defas (G_i\colon i \in I)}$ be an indexed collection of $\Reals$-valued random variables.
We say that $G$ is a \defn{Gaussian process on $I$} when, for all finite sequences of indices ${i_1,\dotsc,i_k \in I}$, the vector ${(G({i_1}),\dotsc,G({i_k}))}$ is Gaussian, where we have written ${G(i)\defas G_i}$ for notational convenience.
A Gaussian process is completely specified by two function-valued parameters: 
a \defn{mean function} $\mu : I \to \Reals$, satisfying
\[
\mu(i) = \EE \bigl( G(i) \bigr), \quad i \in I,
\]
and a positive semidefinite \defn{covariance function} $\kappa\colon I \times I \to \Reals_+$, satisfying
\[
\kappa(i,j) = \mathrm{cov}(G(i),G(j)). 
\]
If $\kappa$ is chosen appropriately, a random function sampled from the Gaussian process is continuous with
probability 1.
\begin{definition}[Gaussian-process-based exchangeable arrays]
  \label{def:GP:based}
  We say that a Bayesian model for an exchangeable array ${X \defas (X_{\oset i j})}$ 
  in $\dataspace$ is \defn{Gaussian-process-based}
   when, for some random function $F$ representing $X$, 
   the process ${F=(F_{x,y,z};\ x,y,z \in [0,1])}$ is Gaussian on $[0,1]^3$.
   We will say that an array $X$ is Gaussian-process-based when its distribution is.
\end{definition}

In terms of \cref{eq:sampform}, a Gaussian-process-based model is one where a Gaussian process prior is placed on the function $f$.
The definition is stated in terms of the space $[0,1]^3$ as domain of the uniform random variables $U$ to match our statement of the Aldous-Hoover
theorem and of previous models. In the case of Gaussian processes, however, it is arguably more natural to use the real line instead of $[0,1]$,
and we note that this is indeed possible: Given an embedding $\phi : [0,1]^3 \to J$ and a Gaussian process $G$ on $J$, 
the process $G'$ on $[0,1]^3$ given by $G'_{x,y,z} = G_{\phi(x,y,z)}$ is Gaussian.
More specifically, if the former has a mean function $\mu$ and covariance function $\kappa$, then the latter has mean $\mu \circ \phi$ and covariance $\kappa \circ (\phi \otimes \phi)$.  We can therefore talk about Gaussian processes on spaces $J$ that can be put into correspondence with the unit interval.

The random variables $X_{ij}$ in \cref{def:GP:based} are real-valued. To model ${\lbrace 0,1\rbrace}$-valued arrays,
such as random graphs, we can use a suitable randomization:
\begin{definition}[noisy sigmoidal/probit likelihood]\label{def:probitlike}
Let $\xi$ be a Gaussian random variable with mean $m \in \Reals$ and variance $v > 0$, and let 
 $\sigma: \Reals \to [0,1]$ be a sigmoidal function.
We define the family $\kernelfamily L r \Reals$ of distributions on $\{0,1\}$ by
$\kernelval L r\{1\} = \EE \bigl( \sigma(r+\xi) \bigr)$.
\end{definition}

Many of the most popular parametric models for exchangeable arrays of random variables can be constructed as (randomizations of) Gaussian-process-based arrays. For a catalog of such models and several nonparametric variants, as well as their covariance functions, see \citep{Lloyd2012}.  Here we will focus on the parametric \defn{eigenmodel}, introduced by Hoff~\citep{Hoff:2007:1,Hoff2011}, and its nonparametric cousin, introduced Xu, Yan and Qi~\citep{Xu2012}. 
To simplify the presentation, we will consider the case of a $\{0,1\}$-valued array.  

\newcommand{\innerproduct}[2]{\langle #1,#2\rangle}
\begin{example}[Eigenmodel~\citep{Hoff:2007:1,Hoff2011}]

In the case of a $\{0,1\}$-valued array, both the eigenmodel and its nonparametric extension can be interpreted as an $L$-randomizations of a Gaussian-process-based array $\Theta\defas(\Theta_{\oset ij})$, where $L$ is given as in \cref{def:probitlike} for some mean, variance and sigmoid.  To complete the description, we define the Gaussian processes underlying $\Theta$.

The eigenmodel is best understood in terms of a zero-mean Gaussian process $G$ on ${\Reals^{d} \times \Reals^d}$. (The corresponding embedding ${\phi : [0,1]^3 \to \Reals^d \times \Reals^d}$ is ${\phi(x,y,z)=\Phi^{-1}(x)\Phi^{-1}(y)}$, where $\Phi^{-1}$ is defined so that ${\Phi^{-1}(U) \in \Reals^d}$ is a vector of independent doubly-exponential (aka Laplacian) random variables when $U$ is uniformly distributed in $[0,1]$.)  The covariance function ${\kappa : \Reals^d \times \Reals^d \to \Reals_+}$ of the Gaussian process $G$ underlying the eigenmodel is simply
\[\label{eq:thekernel}
\kappa(u,v;x,y) = \innerproduct u x  \innerproduct v y, \quad u,v,x,y \in\Reals^d,
\]
where $\innerproduct . . \colon \Reals^d \times \Reals^d \to \Reals$ denotes the dot product, i.e., Euclidean inner product.  
This corresponds with a more direct description of $G$:  in particular,
\[
G(x,y) = {\innerproduct x y}_{\Lambda}
\]
where ${\Lambda \in \Reals^{d\times d}}$ is a $d \times d$ array of independent standard Gaussian random variables and ${{\innerproduct x y }_A = \sum_{n,m} x_n y_m A_{n,m}}$ is an inner product.
\end{example}

A nonparametric counterpart to the eigenmodel was introduced by \citet{Xu2012}:

\begin{example}
\newcommand{\HS}{\mathcal H}
The Infinite Tucker Decomposition model \citep{Xu2012} defines the covariance function on $\Reals^d \times \Reals^d$ to be
\[
\kappa(u,v;x,y) = \kappa'(u,x) \kappa'(v,y), \quad u,v,x,y \in\Reals^d,
\]
where $\kappa' : \Reals^d \times \Reals^d \to \Reals$ is some positive semi-definite covariance function on $\Reals^d$.  This change can be understood as generalizing the inner product in \cref{eq:thekernel} from $\Reals^d$ to a (potentially, infinite-dimensional) reproducing kernel Hilbert space (RKHS).  In particular, for every such $\kappa'$, there is an RKHS $\HS$ such that
\[
\kappa'(x,y) = {\innerproduct {\phi(x)} {\phi(y)}}_{\HS}, \quad x,y \in \Reals^d.
\]
\end{example}

A related nonparametric model for exchangeable arrays, which places fewer restrictions on the covariance structure and is derived directly from the Aldous-Hoover representation, is described in~\citep{Lloyd2012}.

\section{Convergence, Concentration and Graph Limits}
\label{sec:graph:limits}

\begin{it}
  We have already noted that the parametrization of random arrays by functions
  in the Aldous-Hoover theorem is not unique. Our statement of the theorem also
  lacks an asymptotic convergence result such as the convergence of the empirical
  measure in de Finetti's theorem. The tools to fill these gaps have only recently
  become available in a new branch of combinatorics which studies objects known
  as graph limits. This section summarizes a few elementary notions of this rapidly
  emerging field and shows how they apply to the Aldous-Hoover theorem for graphs.
\end{it}

\mbox{ }

\begin{figure*}[t]
  \begin{center}
  \includegraphics[width=2.8cm]{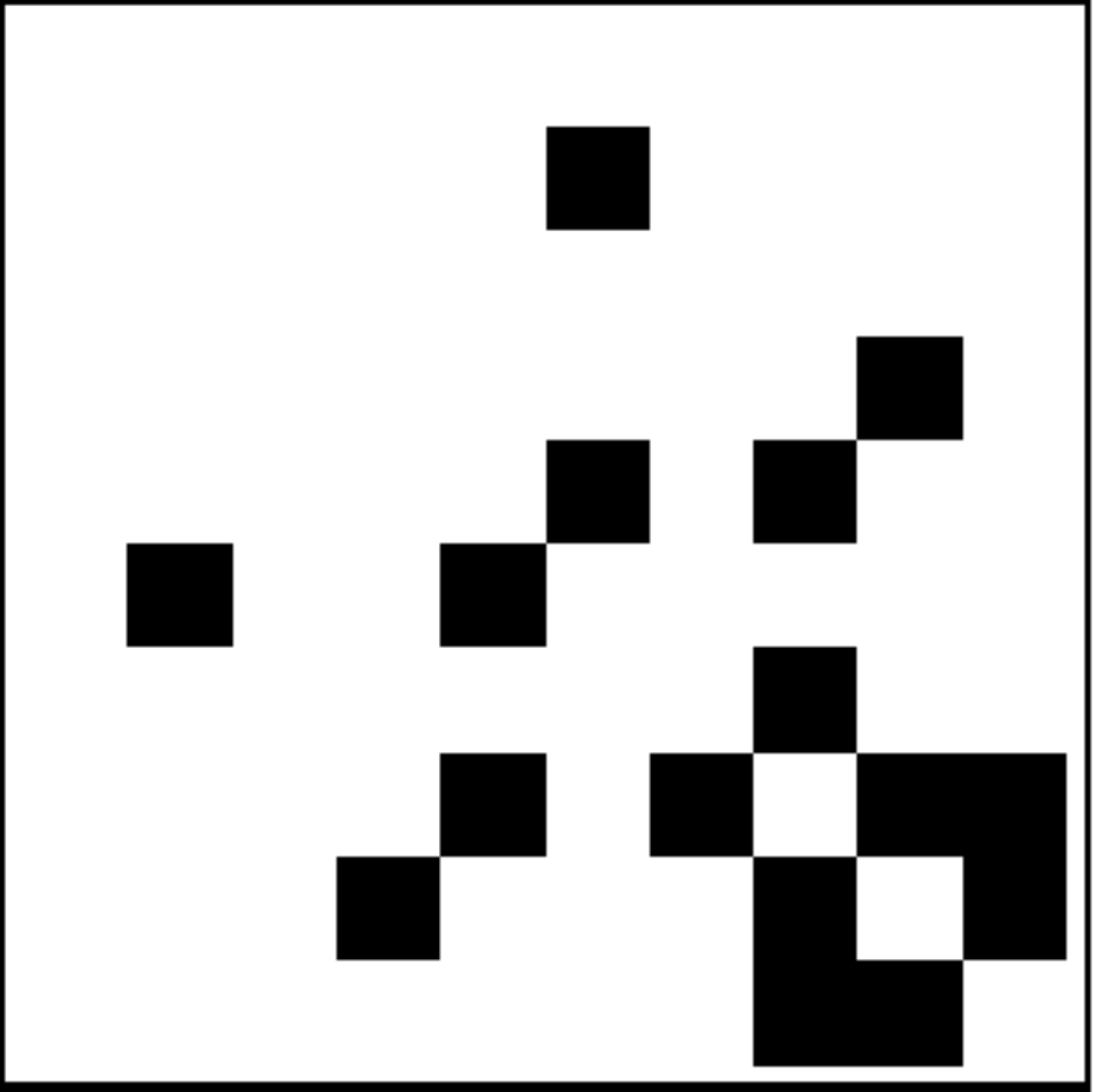}
  \includegraphics[width=2.8cm]{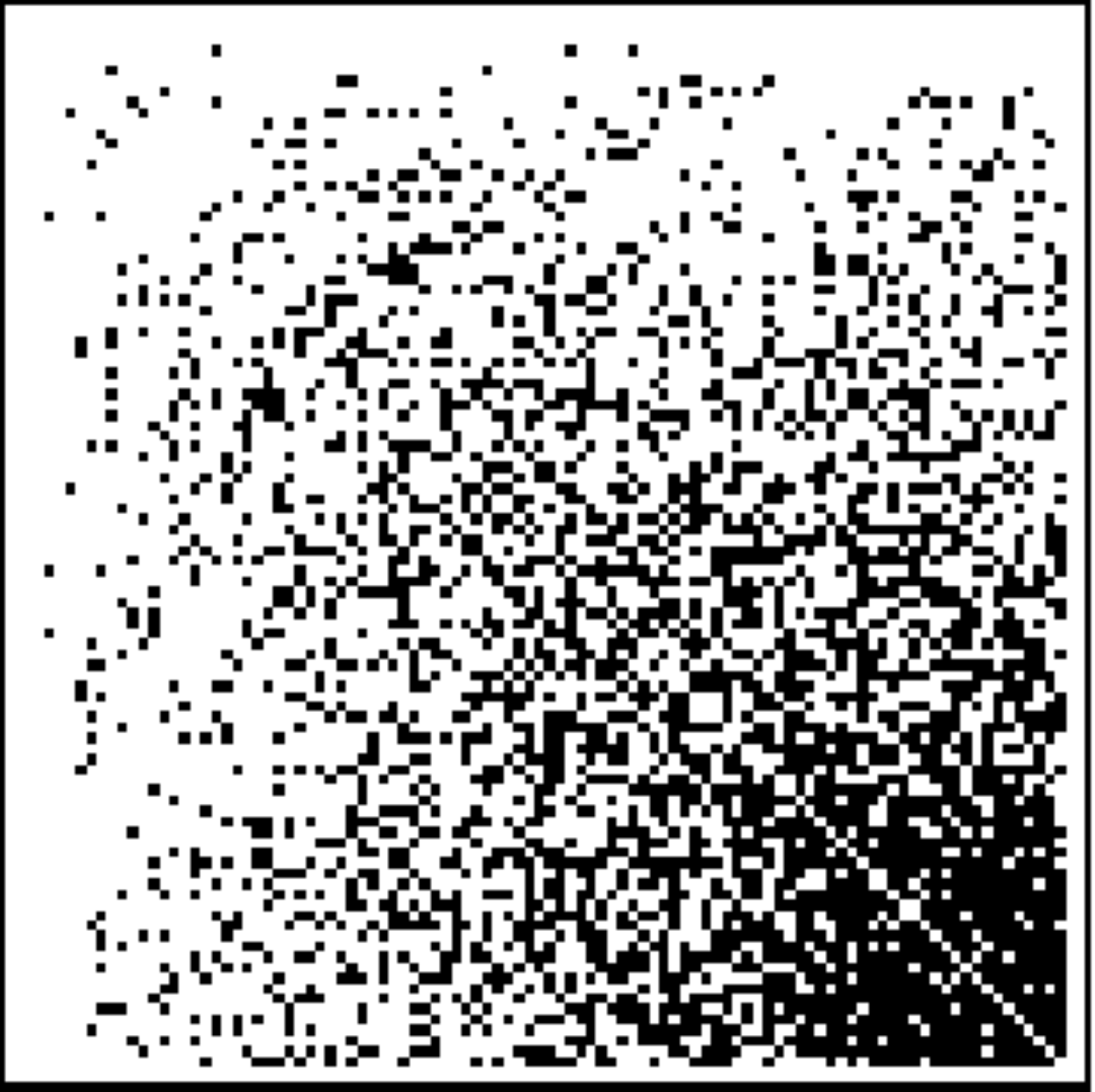}
  \includegraphics[width=2.8cm]{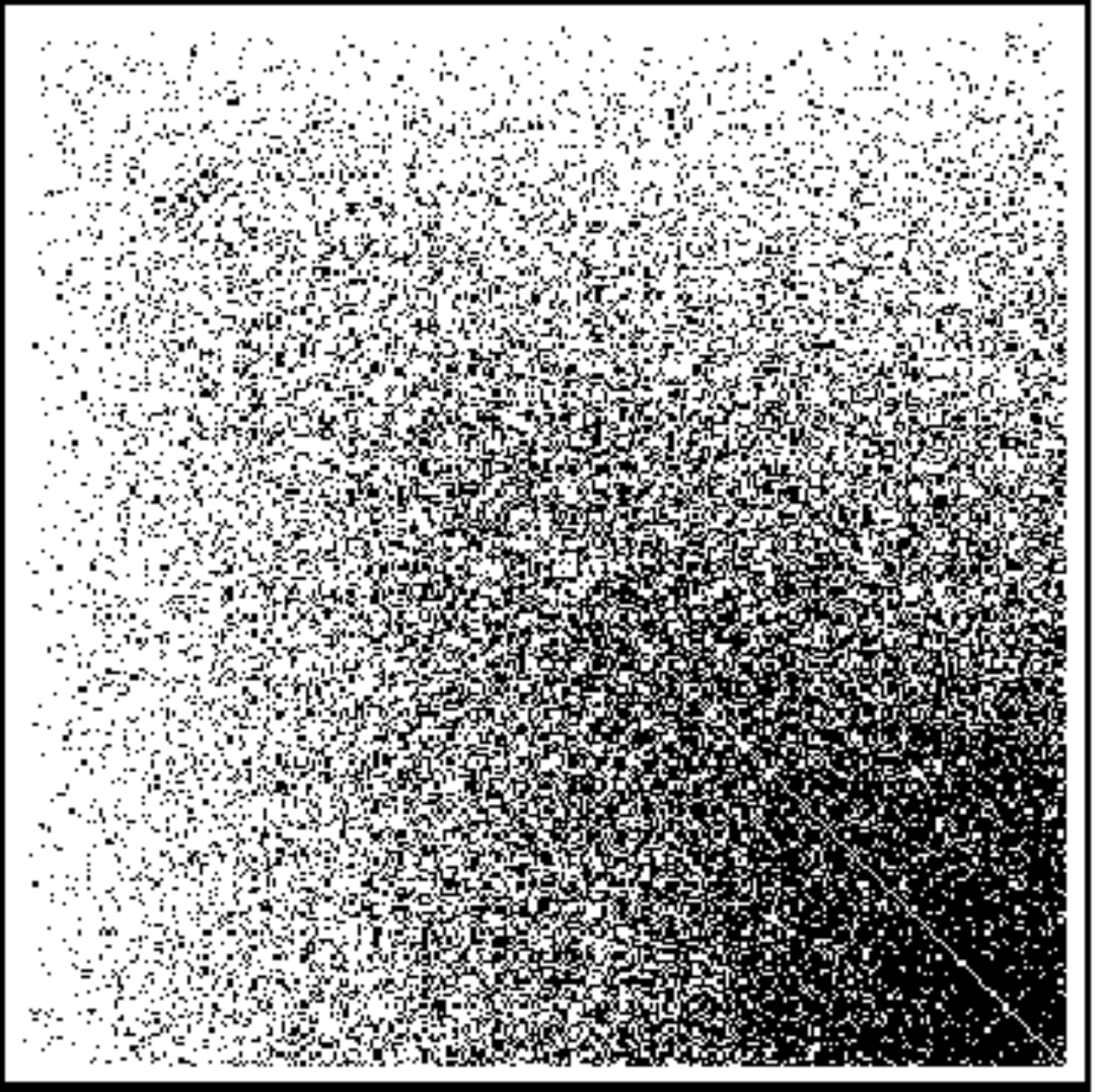}
  \qquad
  \includegraphics[width=2.8cm]{min_function_bw.pdf}\\
  \vspace{0.25cm}
  \includegraphics[width=2.8cm]{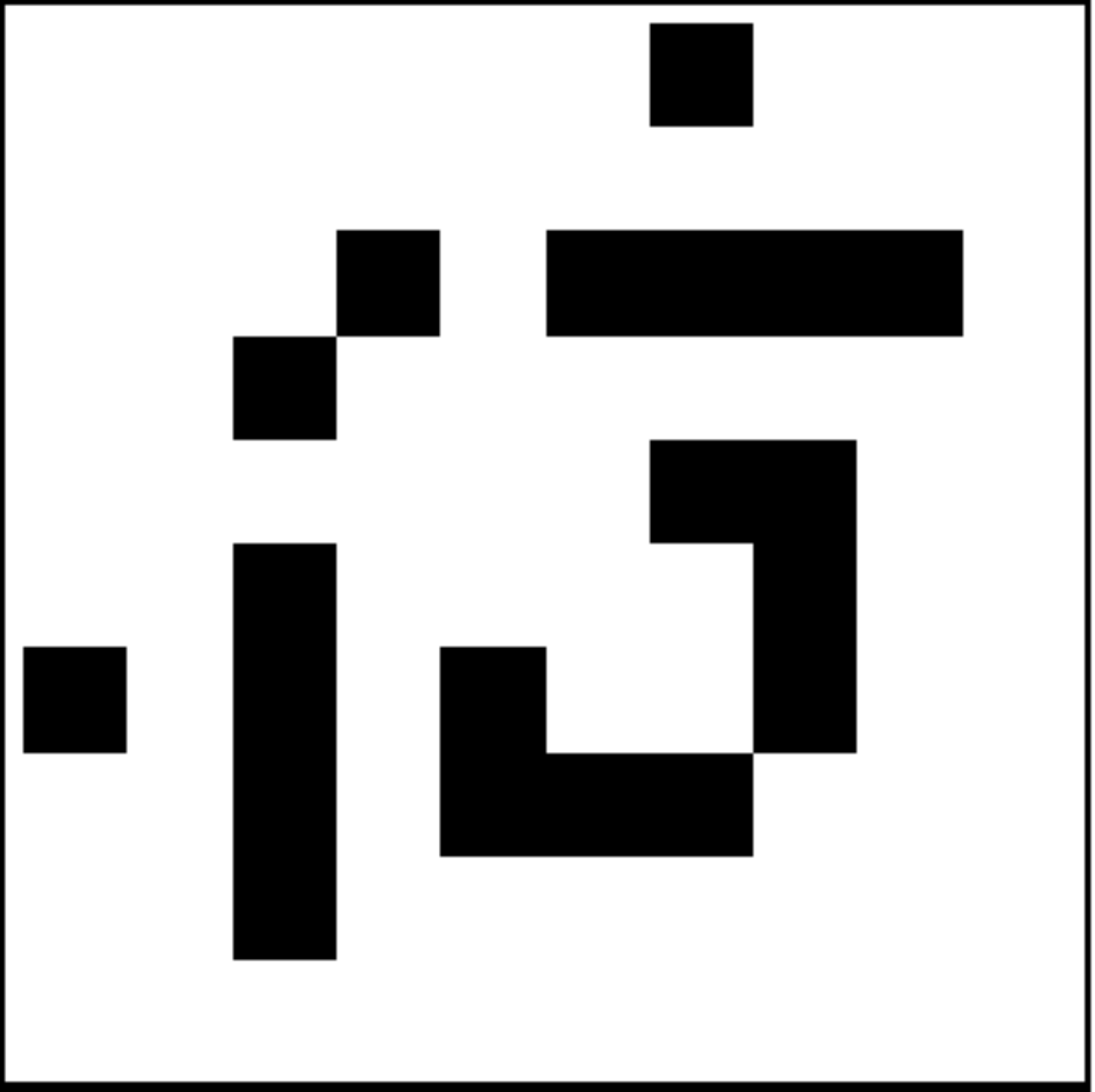}
  \includegraphics[width=2.8cm]{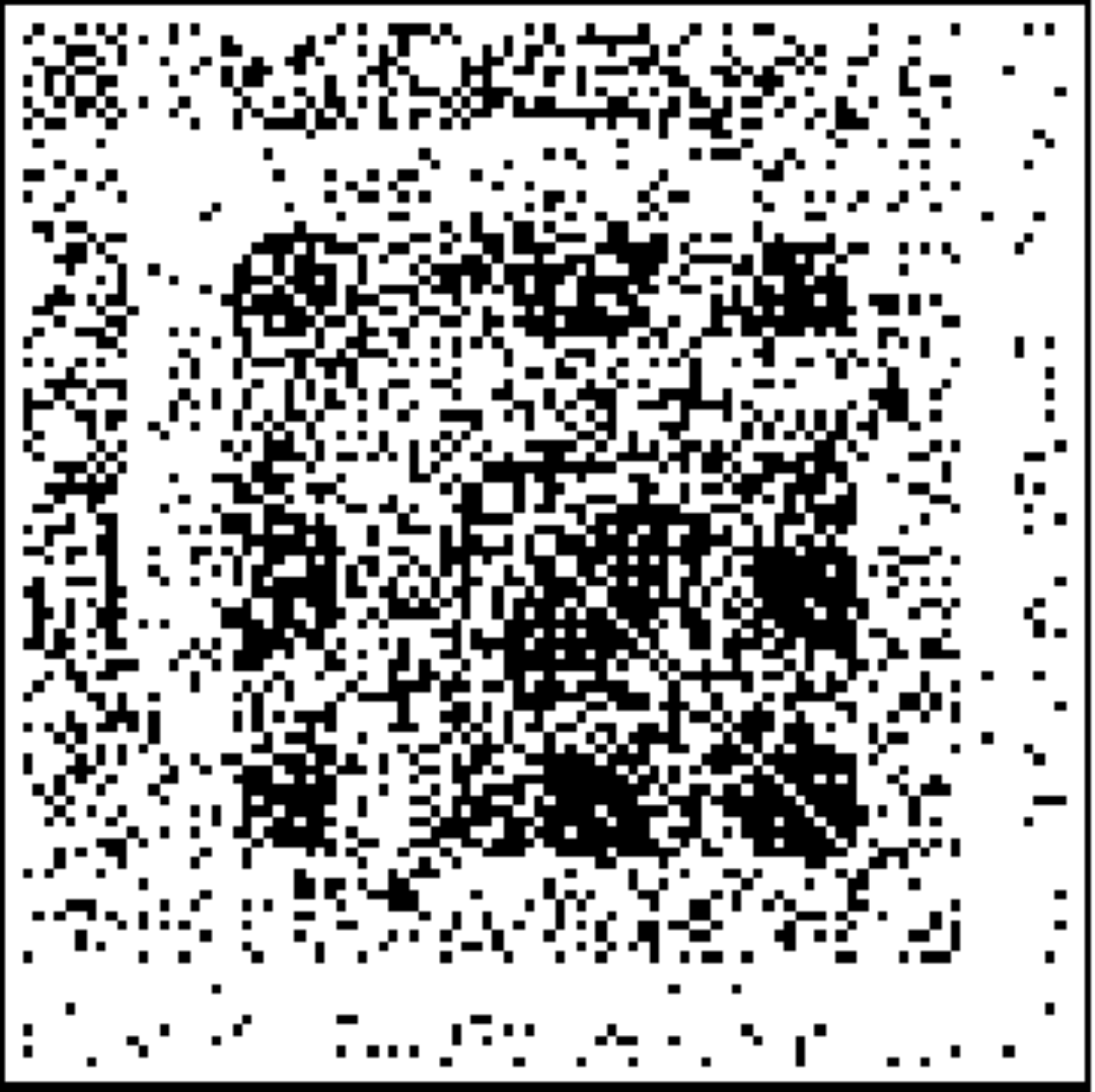}
  \includegraphics[width=2.8cm]{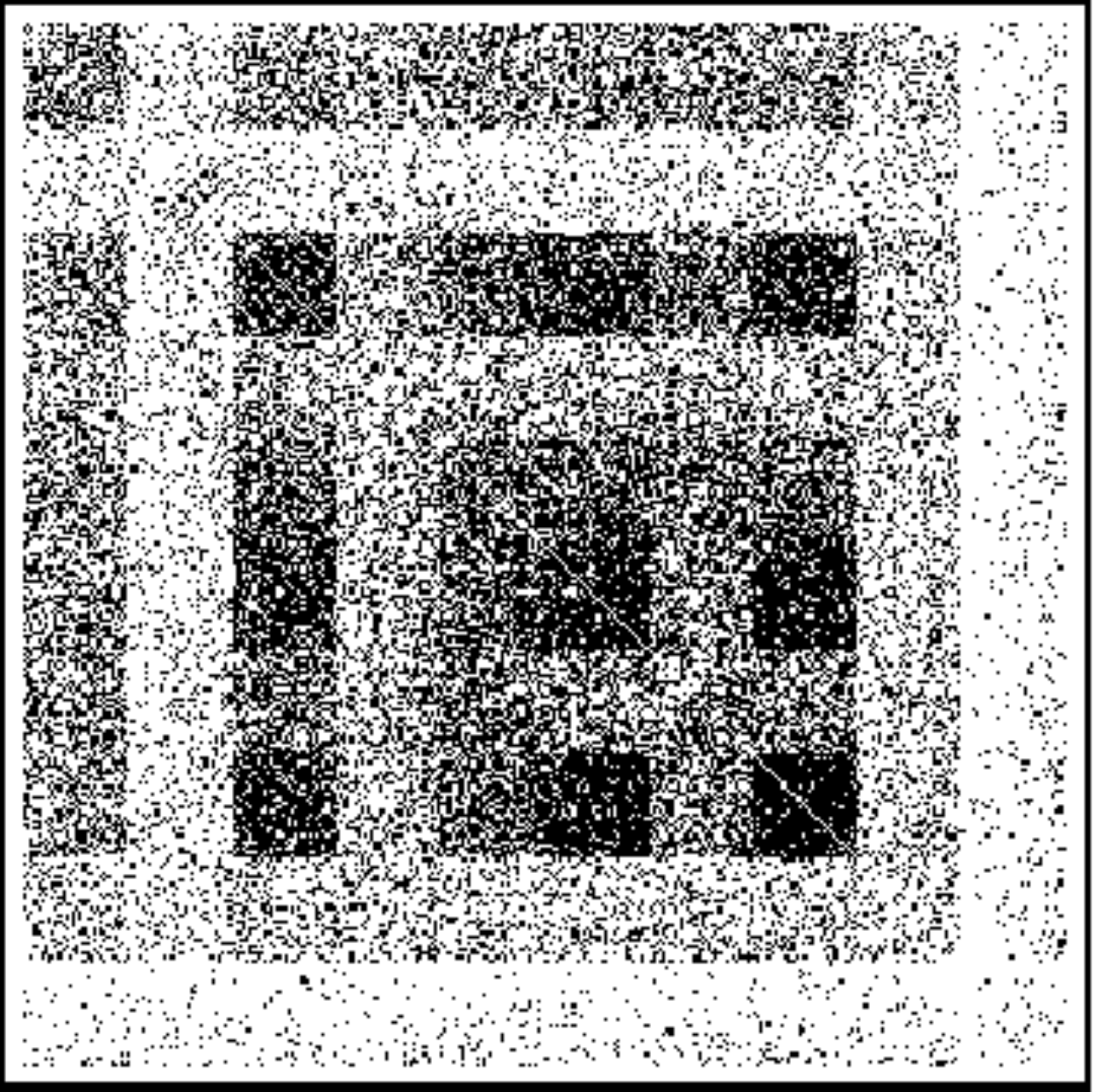}
  \qquad
  \includegraphics[width=2.8cm]{minfunction_perm_bw.pdf}
  \end{center}
  \caption{For graph-valued data, the directing random function $F$ in the Aldous-Hoover representation can be regarded as a limit of
    adjacency matrices: The adjacency matrix of a graph of size $n$ 
    can be represented as a function on $[0,1]^2$ by dividing the square into $n\times n$ patches of equal size. On each patch,
    the representing function is constant, with value equal to the corresponding entry of the adjacency matrix. (In the figure,
    a black patch indicates a value of one and hence the presence of an edge.) 
    As the size of the graph increases, the subdivision becomes finer, and converges to the function depicted on the right
    for $n\to\infty$. 
    Convergence is illustrated here for the two functions from \cref{fig:non:uniqueness}. Since the functions are equivalent,
    the two random graphs within each column are equal in distribution.
  }
  \label{fig:convergence}  
\end{figure*}

\def\vertexset{\mathbf{V}}
\def\edgeset{\mathbf{E}}
\def\cutmetric#1{\| #1 \|_{_{\square}}}
\def\dcut{d_{_{\square}}}
\def\deltacut{\delta_{_{\square}}}
\def\hatgraphon{\widehat{\graphon}}

Convergence results in statistics study the behavior of the empirical
distribution \eqref{eq:emp:measure}. The corresponding object for
exchangeable graphs is an empirical estimate of the graphon, which
is a ``checkerboard function'':
Given a finite graph $g_n$ with $n$ vertices, we subdivide $[0,1]^2$ into
${n\times n}$ square patches, resembling the ${n\times n}$ adjacency matrix. 
We then define a function $w_{g_n}$ with constant value 0 or 1 on each patch, 
equal to the corresponding entry of the adjacency matrix.
We call $w_{g_n}$ the \defn{empirical graphon} of $g_n$. 
Examples are plotted in \cref{fig:convergence}. 
Since $w_{g_n}$ is a
valid graphon, it parametrizes an infinite random graph, even
though $n$ is finite.
Aldous-Hoover theory provides a graph counterpart to the law of large numbers:
\begin{theorem}[\citet{Kallenberg:1999}]
  \label{theorem:kallenberg}
  Let $w$ be a graphon. Suppose we sample a random graph from $w$ one
  vertex at a time, and $G_n$ is the graph given by the first $n$
  vertices. Then the distributions defined by $w_{G_n}$ converge
  weakly to the distribution defined by $w$ with probability 1.
\end{theorem}

A recent development in combinatorics and graph theory is the theory
of graph limits \citep{Lovasz:Szegedy:2006,Lovasz:Szegedy:2007}.
This theory defines a distance measure $\deltacut$ on graphons
(more details below). The metric can be applied to finite graphs,
since the graphs can be represented by the empirical graphon. It
is then possible to study conditions under which sequences of graphs
converge to a limit. It turns out that limits of graphs can be
represented by graphons, and the convergence of graphs corresponds
precisely to the weak convergence of the distributions defined by
the empirical graphons. This theory refines the Aldous-Hoover theory
with a large toolbox of powerful results. We describe a few
aspects in the following. The authoritative (and very well-written)
reference is \citep{Lovasz:2013:A}.

\subsection{Metric definition of convergence}

The most convenient way to define convergence
is by defining a metric: If $d$ is a distance measure,
we can define $w$ as the limit of $w_{g_n}$ if $d(w,w_{g_n})\to 0$ as
$n\to\infty$. The metric on functions which has emerged as the
``right'' choice for graph convergence is called the \defn{cut metric},
and is defined as follows: We first define a norm as
\begin{equation}
  \label{eq:cutmetric}
  \cutmetric{w} \defas \sup_{S,T\subset[0,1]}
    \int_{S\times T} w(x,y) \, \dee x \, \dee y \,.
\end{equation}
(The integral is with respect to Lebesgue measure $\dee x$ because the variables $U_i$ are uniformly distributed.)
Intuitively---if we assume for the moment that $w$ can indeed be thought
of as a limiting adjacency matrix---$S$ and $T$ are subsets of nodes.
The integral \eqref{eq:cutmetric} measures the total number of edges
between $S$ and $T$ in the ``graph'' $w$. Since a partition of the vertices
of a graph into two sets is called a \defn{cut}, $\cutmetric{\argdot}$ is called
the \defn{cut norm}. The distance measure defined by
$\dcut(w,w'):=\cutmetric{w-w'}$ is called the \defn{cut distance}.

Suppose $w$ and $w'$ are two distinct functions which parametrize the
same distribution on graphs. The distance $\dcut$ in general perceives such functions
as different: The functions in \cref{fig:convergence}, for instance,
define the same graph, but have non-zero distance under $\dcut$.
Hence, if we were to use $\dcut$ to define convergence, the two sequences
of graphs in the figure would converge to two different limits. We therefore 
modify $\dcut$ as follows: For any given $w$, let $[w]$ be the set of
all functions $w'$ which define the same random graph. 
\begin{equation}
  \deltacut(\graphon_1,\graphon_2) \defas \inf_{\graphon'\in[\graphon_2]}
  \dcut(\graphon_1,\graphon') \;.
\end{equation}
Informally, we can think of the functions in $[w_2]$ as functions
obtained from $w_2$ by a ``rearrangement'' like the
one illustrated in \cref{fig:non:uniqueness}. 
The definition above says that,
before we measure the distance between $\graphon_1$ and
$\graphon_2$ using $\dcut$, we rearrange $\graphon_2$
in the way that aligns it most closely with $\graphon_1$.
In \cref{fig:non:uniqueness}, this closest rearrangement would
simply reverse the permutation of blocks, so that the two functions
would look identical.

The function $\deltacut$ is called the
\defn{cut pseudometric}: It is not an actual metric, since it can take
value 0 for two distinct functions. It does, however, have all other properties of a metric.
By definition, $\deltacut(w,w')=0$ holds if and only if $w$ and $w'$ parametrize
the same random graph.
\begin{definition}
  \label{def:convergence:cutmetric}
  We say that a sequence $(g_n)_{n\in\Nats}$ of graphs converges if
  $\deltacut(w_{g_n},w)\to 0$ for some measurable function
  ${w:[0,1]^2\to[0,1]}$. The function $w$ is called the limit
  of $(g_n)$, and often referred to as a \defn{graph limit}.
\end{definition}
Clearly, the graph limit is a graphon, and the two terms are used
interchangeably.
This definition indeed provides a metric counterpart to convergence
of exchangeable graph distributions:
\begin{theorem}
  A function $w$ is the graph limit of a sequence 
  $(g_n)_{n\in\Nats}$ of graphs if and only if the random graph
  distributions defined by the empirical graphons $w_{g_n}$ converge 
  weakly to the distribution defined by $w$.
\end{theorem}

\subsection{Unique parametrization in the Aldous-Hoover theorem}
\label{sec:unique:parametrization}

Recall that two graphons can be equivalent, in the sense that they are distinct
functions but define the same random graph (they are weakly isomorphic
in the language of \cref{sec:nonuniqueness}).
The equivalence classes $[w]$ form a partition of the space
$\mathbf{W}$ of all graphons, which 
motivates the definition of a ``quotient space'':
We can define a new space $\quotientspace$
by collapsing each equivalence class to a single point.
Each element ${\hatgraphon\in\quotientspace}$ corresponds to all functions in one equivalence
class, and hence to one specific random graph distribution. Since the
pseudometric $\deltacut$ only assigns distance 0 to two distinct functions
if they are equivalent, it turns into a metric on $\quotientspace$, and
$(\quotientspace,\deltacut)$ is a metric space. Although the elements
of this space are abstract objects, not actual functions, the space
has remarkable analytical properties, and is one
of the central objects of graph limit theory.

Since $\quotientspace$ 
contains precisely one element for each ergodic
distribution on exchangeable graphs, we can obtain a \emph{unique}
parametrization of exchangeable graph models
by using ${\tspace:=\quotientspace}$ as a parameter space:
If ${w\in\mathbf{W}}$ is a graphon and $\hatgraphon$ the corresponding element of 
$\quotientspace$---the element to which $w$ was collapsed in the definition 
of $\quotientspace$---we define a 
family $\kernelfamily \genkernel {\hatgraphon} {\quotientspace}$ of distributions on exchangeable arrays
by taking $\kernelval \genkernel {\hatgraphon}$ to be the distribution
induced by the uniform sampling scheme described by \cref{eq:AH:graph} when $W=w$.

Although the existence of such a probability kernel is not a trivial fact, it follows from a 
technical result of \citet{Orbanz:Szegedy:2012}.
In particular, the Aldous-Hoover theorem for an exchangeable random
graph $G$ can then be written as a unique integral decomposition
\begin{equation}
  \Pr(G\in\,.\,)=\int_{\quotientspace} \kernelvalset \genkernel \argdot {\hatgraphon} \,\nu(\dee \hatgraphon) \;,
\end{equation}
in analogy to the de Finetti representation.

\subsection{Regularity and Concentration}

\def\Szemeredi{Szemer\'edi}

All convergence results we have seen so far provide only asymptotic
guarantees of convergence, but no convergence rates. 
We give two examples of concentration
results from graph limit theory, which address similar
questions as those asked in mathematical statistics and empirical process theory:
How large a graph do we have to observe to obtain reliable estimates?

Underlying these ideas is one of the deepest results of modern graph theory,
Szemeredi's 
regularity lemma: For every very large graph $g$, there is a small, weighted
graph $\hat{g}$ that summarizes all essential structure in $g$. The only condition is that $g$
is sufficiently large. In principle, this means that $\hat{g}$ can be used as an approximation
or summary of $g$, but unfortunately, the result is only valid for graphs which are much larger
than possible in most conceivable applications. There are, however, weaker forms of this result
which hold for much smaller graphs.

To define $\hat{g}$ for a given graph $g$, we proceed as follows: 
Suppose $\Pi \defas \lbrace V_1,\dotsc,V_k\rbrace$ is a partition of $\vertexset(g)$ into $k$ sets.
For any two sets $V_i$ and $V_j$, we define $p_{\oset ij}$ as the probability that two vertices
$v\in V_i$ and $v'\in V_j$, each chosen uniformly at random from its set, are connected by an edge.
That is,
\begin{equation}
  p_{\oset ij} \defas \frac{\#\text{ edges between }V_i,V_j}{|V_i|\cdot|V_j|} \;.
\end{equation}
The graph $\hat{g}_{\Pi}$ is now defined as the weighted graph with vertex set $\lbrace 1,\dotsc, k\rbrace$
and edge weights $p_{\oset ij}$ for edge $(i,j)$.
To compare this graph to $g$, it can be helpful to blow it up to a graph $g_{\Pi}$ of the same size as $g$,
constructed as follows:
\begin{itemize}
\item Each node $i$ is replaced by a  clique of size $|V_i|$ (with all edges weighted by 1).
\item For each pair $V_i$ and $V_j$, all possible edges between the sets are inserted and weighted by $p_{\oset ij}$.
\end{itemize}
If we measure how much two graphs differ in terms of the cut distance,
$g$ can be approximated by $g_{\Pi}$ as follows:
\begin{theorem}[Weak regularity lemma \citep{Frieze:Kannan:1999}]
  \label{theorem:regularity:weak}
  Let $k\in\Nats$ and let $g$ be any graph. There is a partition $\Pi$ of $\vertexset(g)$
  into $k$ sets such that $\dcut(g,g_{\Pi})\leq 2(\sqrt{\log(k)})^{-1}$.
\end{theorem}
This form of the result is called ``weak'' since it uses a less restrictive definition of
what it means for $g$ and $g_{\Pi}$ to be close then \Szemeredi's original result. The weaker
hypothesis makes the theorem applicable to graphs that are,
by the standards of combinatorics, of modest size.

A prototypical concentration result based on \cref{theorem:regularity:weak} is the following:
\begin{theorem}[{\citep[][Theorem 8.2]{Lovasz:2009:1}}]
  Let $f$ be a real-valued function on graphs, which is smooth in the sense that ${|f(g)-f(g')|\leq\dcut(g,g')}$ for any two graphs 
  $g$ and $g'$ defined on the same vertex set. 
  Let $G(k,g)$ be a random graph of size $k$ sampled uniformly from $g$. Then the distribution of
  $f(G(k,g))$ concentrates around some value $f_0\in\Reals$, in the sense that 
  \begin{equation}
    \Pr\Bigl\lbrace |f(G(k,g))-f_0|>\frac{20}{\sqrt k}\Bigr\rbrace < 2^{-k}\;.
  \end{equation}
\end{theorem}
The relevance of such results to statistics becomes evident if we
think of $f$ as a statistic of a graph or network (such as the edge
density) which we try to estimate from an observed subgraph of size $k$.
Results of this type, for graphs and other random structures, 
are collectively known under the term \emph{property testing}, and
are covered by a sizeable literature in combinatorics and theoretical
compute science \citep{Alon:Spencer:2008,Lovasz:2013:A}.

\section{Exchangeability and higher-dimensional arrays}
\label{sec:d-arrays}

\def\mset#1{\{\!\{#1\}\!\}}

\begin{it}
The theory of exchangeable arrays extends beyond 2-dimensional arrays, and, indeed, some of the more exciting implications and applications of the theory rely on the general results.  In this section we begin by defining the natural extension of (joint) exchangeability to higher dimensions, and then give higher-dimensional analogues of the theorems of Aldous and Hoover, due to Kallenberg.  These theorems introduce exponentially-many additional random variables as the dimension increases, but only a linear number are necessary to produce an arbitrarily good approximation.  The presentation owes much to Kallenberg~\citep{Kallenberg:1999}.
\end{it}

\begin{definition}[jointly exchangeable $d$-arrays]
Let $(X_{k_1,\dotsc,k_d})$ 
be a $d$-dimensional array (or simply $d$-array) of random variables in $\dataspace$.
We say that $X$ is \defn{jointly exchangeable} when
\[
(X_{k_1,\dotsc,k_d}) \eqdist (X_{\pi(k_1),\dotsc,\pi(k_d)})
\]
for every permutation $\pi$ of $\Nats$.
\end{definition}

As in the 2-dimensional representation result,
a key ingredient in the characterization of higher-dimensional jointly exchangeable $d$-arrays will be an indexed collection $U$ of \iid latent random variables.  
In order to define the index set for $U$, 
let $\tilde \Nats^d$ be the space of multisets $J \subseteq \Nats$ of cardinality $|J| \le d$.  E.g., $\mset{1,1,3} \in \tilde \Nats^3 \subseteq \tilde \Nats^4$.
Rather than two collections---a sequence $(U_i)$ indexed by $\Nats$, and a triangular array $(U_{\mset{ i,j}})$ indexed by multisets of cardinality 2---we will use a single \iid collection $U$ indexed by elements of $\tilde \Nats^d$.
For every $I \subseteq [d] \defas \{1, \dotsc, d\}$,
we will write $\tilde k_I$ for the multiset 
\[
\mset{ k_i : i \in I }
\]
and write
\[
(U_{\tilde k_I} ;\ I \in 2^{[d]}\setminus \emptyset)
\]
for the element of the function space $[0,1]^{2^{[d]}\setminus \emptyset}$ that maps each nonempty subset $I \subseteq [d]$ to the real $U_{\tilde k_I}$, i.e., the element in the collection $U$ indexed by the multiset $\tilde k_I \in \tilde \Nats^{|I|} \subseteq \tilde \Nats^d$.

\begin{theorem}[Aldous, Hoover]
Let $U$ be an \iid collection of uniform random variables indexed by multisets $\tilde \Nats^d$.
A random $d$-array $X\defas (X_k;\ k\in\Nats^d)$ is jointly exchangeable if and only if there is random measurable function $F : [0,1]^{2^{[d]}\setminus \emptyset} \to \dataspace$ such that
\[\label{eq:jechar}
(X_k;\ k \in \Nats^d) \eqdist (F(U_{\tilde k_I} ;\ I \in 2^{[d]}\setminus \emptyset);\ k \in \Nats^d).
\]
\end{theorem}

When $d=2$, we recover \cref{theorem:AH:2D:jointly} characterizing two-dimensional exchangeable arrays.  Indeed, if we write $U_i \defas U_{\mset i}$ and $U_{ij} \defas U_{\mset{i,j}}$ for notational convenience,
then the right hand side of \cref{eq:jechar} reduces to 
\[
%(X_{\oset ij} ;\ i,j \in \Nats) \eqdist 
(F(U_i, U_j, U_{ij});\ i,j \in \Nats)
\]
for some random $F : [0,1]^3 \to \dataspace$.
When $d=3$, we instead have %that $X$ is equal in distribution to the array
\[
%(X_{\oset ijk} ;\ i,j,k \in \Nats) \eqdist 
(F(U_i, U_j, U_k, U_{ij}, U_{ik}, U_{jk}, U_{ijk});\ i,j,k \in \Nats)
\]
for some random $F : [0,1]^7 \to \dataspace$, where we have
additionally taken $U_{ijk} \defas U_{\mset {i,j,k}}$ for notational
convenience.

\subsection{Separately exchangeable $d$-arrays}

As in the two-dimensional case, arrays with certain additional symmetries can be treated as special cases.  In this section, we consider separate exchangeability in the setting of $d$-arrays, and in the next section we consider further generalizations.  We begin by defining:
\begin{definition}[separately exchangeable $d$-arrays]
We say that $d$-array $X$ is \defn{separately exchangeable}
when
\[
(X_{k_1,\dotsc,k_d}) \eqdist (X_{\pi_1(k_1),\dotsc,\pi_d(k_d)})
\]
for every collection $\pi_1,\dotsc,\pi_d$ of permutations of $\Nats$.
\end{definition}
For every $J \subseteq [d]$, let $1_J$ denote its indicator function (i.e., $1_J(x) = 1$ when $x \in J$ and $0$ otherwise), and let the vector $k_J \in \NNInts^d \defas \{0,1,2,\dotsc\}^d$ be given by
\[
k_J \defas ( k_1\, 1_J(1), \dotsc, k_d\, 1_J(d)).
\]
In order to represent separately exchangeable $d$-arrays, we will use a collection $U$ of \iid uniform random variables indexed by vectors $\NNInts^d$.
Similarly to above, we will write
\[
(U_{k_I}; \ I \in 2^{[d]}\setminus \emptyset)
\]
for the element of the function space $[0,1]^{2^{[d]}\setminus \emptyset}$ that maps each nonempty subset $I \subseteq [d]$ to the real $U_{k_I}$, i.e., the element in the collection $U$ indexed by the vector $k_I$.
Then we have:

\begin{corollary}
Let $U$ be an \iid collection of uniform random variables indexed by vectors $\NNInts^d$.
A random $d$-array $X\defas (X_k;\ k\in\Nats^d)$ is separately exchangeable if and only if there is random measurable function $F : [0,1]^{2^{[d]}\setminus \emptyset} \to \dataspace$ such that
\[
\label{eq:sechar}
(X_k;\ k \in \Nats^d) \eqdist (F(U_{k_I} ;\ I \in 2^{[d]}\setminus \emptyset);\ k \in \Nats^d).
\]
\end{corollary}

We can consider the special cases of $d=2$ and $d=3$ arrays.  Then we have, respectively,
\[
(F(U_{i0}, U_{0j}, U_{ij});\ i,j \in \Nats)
\]
for some random $F : [0,1]^3 \to \dataspace$; and 
\[
%(X_{\oset ijk} ;\ i,j,k \in \Nats) \eqdist 
(F(U_{i00}, U_{0j0}, U_{00k}, U_{ij0}, U_{i0k}, U_{0jk}, U_{ijk});\ i,j,k \in \Nats)
\]
for some random $F : [0,1]^7 \to \dataspace$. As we can see, jointly exchangeable arrays, which are required to satisfy fewer symmetries than their separately exchangeable counterparts, may take $U_{ij0} = U_{0ij} = U_{i0j} = U_{ji0} = \dotsc$.  Indeed, one can show that these additional assumptions make jointly exchangeable arrays a strict superset of separately exchangeable arrays, for $d\ge 2$.

\subsection{Further generalizations}

In applications, it is common for the distribution of an array to be invariant to permutations that act simultaneously on \emph{some but not all} of the dimensions.  E.g., if the first two dimensions of an array index into the same collection of users, and the users are \emph{a priori} exchangeable, then a sensible notion of exchangeability for the array would be one for which these first two dimensions could be permuted jointly together, but separately from the remaining dimensions.  

More generally, we consider arrays that, given a partition of the dimensions of an array into classes, are invariant to permutations that act jointly within each class and separately across classes.  More carefully:
\begin{definition}[$\pi$-exchangeable $d$-arrays]
Let $\pi = \{I_1,\dotsc,I_m\}$ be a partition of $[d]$ into disjoint classes, and let $p=(p^I;\ I \in \pi)$ be a collection of permutations of $\Nats$, indexed by the classes in $\pi$. We say that a $d$-array $X$ is $\pi$-exchangeable when
\[
(X_{k_1,\dotsc,k_d};\ k \in \Nats^d) 
\eqdist
(X_{p^{\pi_1}(k_1),\dotsc,p^{\pi_d}(k_d)};\ k \in \Nats^d),
\]
for every collection $p$ of permutations, 
where $\pi_i$ denotes the subset $I \in \pi$ containing $i$.
\end{definition}
We may now cast both jointly and separately exchangeable arrays as $\pi$-exchangeable arrays for particular choices of partitions $\pi$.  In particular, when 
$\pi = \{[d]\}$ % = \{\{1,\dotsc,d\}\}$ 
we recover joint exchangeability, 
and when $\pi = \{\{1\},\dotsc,\{d\}\}$, we recover separate exchangeability.  Just as we characterized jointly and separately exchangeable arrays, we can characterize $\pi$-exchangeable arrays.

\newcommand{\GV}[1]{\mathcal V(#1)}
Let $\pi$ be a partition of $[d]$.  In order to describe the representation of $\pi$-exchangeable $d$-arrays, we will again need a collection $U$ of \iid uniform random variables, although the index set is more complicated than before:  Let $\GV \pi \defas {\mathsf X}_{I \in \pi} \tilde \Nats^{|I|}$ denote the space of functions taking classes $I \in \pi$ to multisets $J \subseteq \Nats$ of cardinality $J \le |I|$.  We will then take $U$ to be a collection of \iid uniform random variables indexed by elements in $\GV \pi$.  

When $\pi = \{[d]\}$, $\GV \pi$ is equivalent to the space $\tilde \Nats^d$ of multisets of cardinality no more than $d$, in agreement with the index set in the jointly exchangeable case.  The separately exchangeable case is also instructive: there $\pi = \{ \{1\},\dotsc,\{d\}\}$ and so $\GV \pi$ is equivalent to the space of functions from $[d]$ to $\tilde \Nats^1$, which may again be seen to be equivalent to the space $\NNInts^d$ of vectors, where $0$ encodes the empty set $\emptyset \in \tilde \Nats^1 \cap \tilde \Nats^0$.  For a general partition $\pi$ of $[d]$, an element in $\GV \pi$ is a type of generalized vector, where, for each class $I \in \pi$ of dimensions that are jointly exchangeable, we are given a multiset of indices.

For every $I \subseteq [d]$, let $\tilde k_{\pi I} \in \GV \pi$ be given by 
\[
\tilde k_{\pi I}(J) = \tilde k_{I \cap J},\quad J \in \pi,
\]
where $\tilde k_J$ is defined as above for jointly exchangeable arrays.
We will write 
\[
(U_{\tilde k_{\pi I}};\ I \in 2^{[d]}\setminus \emptyset)
\]
for the element of the function space $[0,1]^{2^{[d]}\setminus \emptyset}$ that maps each nonempty subset $I \subseteq [d]$ to the real $U_{\tilde k_{\pi I}}$, i.e., the element in the collection $U$ indexed by the generalized vector $\tilde k_{\pi I}$.
Then we have:

\begin{corollary}[Kallenberg {\citep{Kallenberg:1999}}]
Let $\pi$ be a partition of $[d]$, and 
let $U$ be an \iid collection of uniform random variables indexed by generalized vectors $\GV \pi$.
A random $d$-array $X\defas (X_k;\ k\in\Nats^d)$ is $\pi$-exchangeable if and only if there is random measurable function $F : [0,1]^{2^{[d]}\setminus \emptyset} \to \dataspace$ such that
\[
(X_k;\ k \in \Nats^d) \eqdist (F(U_{\tilde k_{\pi I}} ;\ I \in 2^{[d]}\setminus \emptyset);\ k \in \Nats^d).
\]
\end{corollary}

\subsection{Approximations by simple arrays}

These representational results require a number of latent random variables exponential in the dimension of the array, i.e., roughly twice as many latent variables are needed as the entries generated in some subarray.  
Even if a $d$-array is sparsely observed, each observation requires the introduction of potentially $2^d$ variables.  (In a densely observed array, there will be overlap, and most  latent variables will be reused.) 

Regardless of whether this blowup poses a problem for a particular application, it is interesting to note that exchangeable $d$-arrays can be approximated by arrays with much simpler structure, known as \defn{simple arrays}.

\begin{definition}[simple $d$-arrays]
Let $U = (U^I_k;\ I \in \pi, k\in \Nats)$ be an \iid collection of uniform random variables.
We say that a $\pi$-exchangeable $d$-array $X$ is \defn{simple} when 
there is a random function $F\colon [0,1]^{[d]} \to \dataspace$ such that 
\[
(X_k;\ k \in \Nats^d)
\eqdist
(F(U^{\pi_1}_{k_1},\dots,U^{\pi_d}_{k_d});\ k \in \Nats^d),
\]
where $\pi_j$ is defined as above.
\end{definition}

Again, it is instructive to study special cases: in the jointly exchangeable case, 
taking $U_j \defas U^{\{[d]\}}_j$,
we get
\[
(F(U_{k_1},\dotsc,U_{k_d}); k \in \Nats^d)
\]
and, in the separately exchangeable case, we get
\[
(F(U^1_{k_1},\dotsc,U^d_{k_d});\ k \in \Nats^d),
\]
taking $U^i_j \defas U^{\{i\}}_j$.
We may now state the relationship between general arrays and simple arrays:

\begin{theorem}[simple approximations, Kallenberg {\citep[][Thm.~2]{Kallenberg:1999}}]
Let $X$ be a $\pi$-exchangeable $d$-array.  Then there exists a sequence of simple $\pi$-exchangeable arrays $X^1, X^2,\dotsc$ such that, for all finite sub-arrays $X_J \defas (X_k; k \in J)$, $J \subseteq \Nats^d$, the distributions of $X_J$ and $X^n_J$ are mutually absolutely continuous, and the associated densities tend \emph{uniformly to} $1$ as $n \to \infty$ for fixed $J$.
\end{theorem}

\section{Sparse random structures and networks}
\label{sec:sparsity}

\begin{it}
  Exchangeable random structures are not ``sparse''. 
  In an exchangeable infinite graph, for example, 
  the expected number of edges attached to each node is either infinite or zero.
  In contrast, graphs representing network data typically have a finite number
  of edges per vertex, and exhibit properties like power-laws
  and ``small-world phenomena'', which can only occur in sparse graphs. 
  Hence, even though exchangeable graph models are widely used in network analysis,
  they are inherently misspecified. 
  Since the lack of sparseness is a direct mathematical consequence
  of exchangeability, networks and sparse random structures pose a problem that seems
  to require genuinely non-exchangeable models. The development of a coherent theory
  is, despite intense efforts in mathematics, a largely unsolved problem.
  In this section, we make the problem more precise and describe how, at least in principle,
  exchangeability might be substituted by other symmetry properties.
  The topic raises a host of challenging questions to which, in most cases, we have no answers.
\end{it}
\mbox{ }

\subsection{Dense vs Sparse Random Structures}

In an exchangeable structure, events either never occur, or they occur
infinitely often with a fixed (though unknown) probability. The
simplest example is an exchangeable binary sequence: 
By de Finetti's theorem, the binary variables are conditionally \iid 
with a Bernoulli distribution. If we sample infinitely often,
conditionally on the random Bernoulli parameter 
taking value ${p\in[0,1]}$, 
the fraction of ones in the sequence
will be precisely $p$. Therefore, we either observe a constant
proportion of ones (if ${p>0}$) or no ones at all (if ${p=0}$). 

In an exchangeable graph, rather than ones and zeros, we have to
consider the possible subgraphs (single edges, triangles, five-stars,
etc). Each possible subgraph occurs either never, or infinitely
often. 
Since an infinite graph may have infinitely many edges even if it is sparsely connected, the number of edges is best
quantified in terms of a rate: 
\begin{definition}
  Let $(g_n)$ be a sequence of graphs ${g_n=(v_n,e_n)}$, where $g_n$
  has $n$ vertices.
  We say that the sequence is \defn{sparse} if, as $n$ increases,
  $|e_n|$ is of size $O(n)$
  (is upper-bounded by $c\cdot n$ for a 
  constant $c$).
  It is called \defn{dense} if ${|e_n|=\Omega(n^2)}$ (lower-bounded by ${c\cdot n^2}$ for a constant $c$).
\end{definition}
If a random graph is sampled step-wise one vertex at a time, the
partial graphs at each step also form a sequence, and we can refer
to the random graph as dense or sparse, depending on whether the
sequence is dense or sparse. (This definition has to be used with
caution, since changing the order in which vertices are generated
may affect the rate.) 
A typical example of dense random graphs are
infinite random graphs in which each vertex has infinite
degree.
Random graphs with bounded
degrees are sparse.
Many important types of graph and array data are inherently sparse: In
a social 
network with billions of users, individual users do not, 
on average, have billions of friends.
\begin{fact}
  Exchangeable graphs are not sparse.
  If a random graph is exchangeable, it is either dense or empty.
\end{fact}
The argument is simple: Let $G_n$ be an $n$-vertex random undirected graph sampled according to \cref{eq:AH:graph}. 
The expected proportion of edges in present in $G_n$, out of all ${{n \choose 2} = \frac{n(n-2)}{2}}$ possible edges, is independent of $n$ and given by 
${\varepsilon\defas\frac{1}{2}\int_{[0,1]^2} W(x,y) \, \dee x \, \dee y}$.
(The factor $\frac{1}{2}$ occurs since $W$ is symmetric.) 
If ${\varepsilon = 0}$, it follows that $G_n$ is empty with probability one and therefore trivially sparse.  On the other hand, if ${\varepsilon > 0}$, 
we have ${\varepsilon\cdot {n \choose 2} =\Theta(n^2)}$ edges in expectation and so, by the law of large numbers, $G_n$ is dense with probability one.

\begin{remark}[Graph limits are dense]
  The theory of graph limits described in \cref{sec:graph:limits} 
  is closely related to exchangeability, and is
  inherently a theory of dense graphs: If we construct a sequence
  of graphs with sparsely growing edge sets, convergence in cut metric
  is still well-defined, but the limit object is always the empty graphon,
  \ie a function on $[0,1]^2$ which vanishes almost everywhere.
\end{remark}
One possible way to generate sparse graphs is of course to modify
the sampling scheme for exchangeable graphs to generate fewer edges.
\begin{example}[The BJR model \citep{Bollobas:Janson:Riordan:2007}]
  There is a very simple way to modify the Aldous-Hoover approach into
  one that generates sparse random graphs:
  Suppose we sample a finite graph with a fixed number
  $n$ of vertices. We simply multiply the
  probability in our usual sampling scheme by $1/n$:
  \begin{equation*}
    \label{eq:bollobas:riordan}
    X_{ij}\sim\mbox{Bernoulli}\Bigl(\frac{1}{n}w(U_i,U_j)\Bigr) 
    \qquad\text{ for }\qquad
    i,j\leq n\;.
  \end{equation*}
  Comparison with our argument why exchangeable graphs are dense immediately
  shows that a graph sampled this way is sparse.
  More generally, we can multiply $w$ by some other rate function
  $\rho_n$ (instead of specifically ${\rho_n=1/n}$), and ask how this
  model behaves for ${n\to\infty}$.
  Statistical properties of this model are studied by
  \citet*{Bickel:Chen:Levina:2011}, who consider the behavior of moment
  estimators for the edge density, triangle density and other subgraph
  densities.
\end{example}
An obvious limitation of the BJR model
is that it does not actually attempt to model network structure:
It can equivalently
be sampled by sampling from a graphon as in \eqref{eq:AH:graph}
and then deleting each edge independently at random, with probability $(1-\rho_n)$.
(In the  parlance of random graph theory, this is exchangeable
sampling followed by \iid bond percolation.)
In other words, the BJR model modifies an exchangeable graph to 
fit a first-order statistic of a network (the number of edges), but
it cannot generate typical network structures, such as power
laws. 

\subsection{Beyond exchangeability: Symmetry and ergodic theory}
\label{sec:symmetry}

The example of networks and sparse structures shows that there are
important random structures which are not exchangeable. This raises
the question whether integral decompositions and statistical models,
which we have throughout derived from exchangeability, can be obtained
in a similar manner for structures that are not exchangeable.
In principle, that is possible: Exchangeability is a special
case of a probabilistic symmetry. It turns
out that integral decompositions can be derived from
much more general symmetries than exchangeability.

A \defn{probabilistic symmetry} is defined by choosing
a group $\mathbb{G}$ of transformations ${g: \xspace_\infty
  \to \xspace_\infty}$.
A random structure $\Xinf$ is called \defn{invariant to $\mathbb G$}
or \defn{$\mathbb G$-symmetric} 
${g(X) \eqdist X}$ for all ${g \in \mathbb G}$.
If so, we also say that the distribution of $\Xinf$ is $\mathbb{G}$-invariant.
For example, a sequence of random variables is exchangeable if and
only if the distribution of the sequence is invariant under the
group of permutations of $\Nats$ acting on the indices of the sequence.
Exchangeability of arrays (as in the Aldous-Hoover theorem)
corresponds with a subgroup generated by row and column permutations.
Invariant measures play a key role in several fields of mathematics, especially ergodic theory.

A very general result, the \emph{ergodic decomposition theorem}, shows that integral decompositions
of the form \eqref{eq:integral:decomp} are a general consequence of probabilistic symmetries,
rather than specifically of exchangeability. The general theme is that
there is some correspondence of the form
\begin{equation*}
  \text{ invariance property } \longleftrightarrow \text{ integral decomposition } \;.
\end{equation*}
In principle, Bayesian models can be constructed based on any type of symmetry,
as long as this symmetry defines a useful set of ergodic distributions.
The following statement of the
ergodic decomposition theorem glosses over various technical details;
for a precise statement, see \citep[][Theorem A1.4]{Kallenberg:2005}.
\begin{theorem}[{Varadarajan \citep{Varadarajan:1963}}]
  \label{theorem:varadarajan}
  Let $\mathbb G$ be a ``nice'' group acting on a space $\structspace$ of infinite structures.
  Then there exists a family $\mathcal E \defas \kernelfamily \genkernel \theta \tspace$ of distributions on $\structspace$ such that,  
  if the distribution of a random structure $X_{\infty}$ is $\mathbb G$-invariant, 
  it has a representation
  of the form
  \begin{equation} %xxx
    \label{eq:varadarajan}
    \Pr(X_{\infty}\in\,.\,)=\int_{\tspace} \kernelvalset \genkernel \argdot \theta  \,\nu(\dee \theta) \,
  \end{equation}
  for a \emph{unique} distribution $\nu$ on $\tspace$.
  The distributions $\kernelval \genkernel \theta \in \mathcal E$ are the so-called \emph{ergodic distributions} associated with $\mathbb G$. 
\end{theorem}

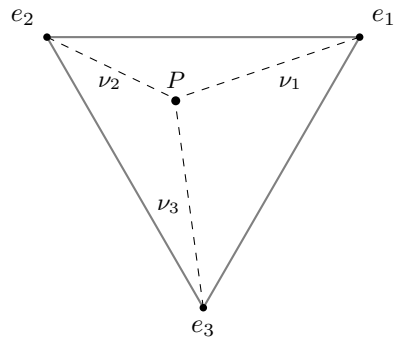
\begin{figure}
  \centering{
    \begin{tikzpicture}[scale=0.8]%[transform canvas={scale=0.7,xshift=13cm,yshift=-1.9cm}]
      \begin{scope}%[scale=0.8]
        \coordinate (e1) at (30:3cm); 
        \coordinate (e2) at (150:3cm); 
        \coordinate (e3) at (-90:3cm);
        
        \draw [thick,gray] (e1.south) -- (e2.north east) -- (e3.north west) -- cycle;
        
        \node [scale=0.3,circle,fill,label=above right:{$e_1$}] at (e1) {};
        \node [scale=0.3,circle,fill,label=above left:{$e_2$}] at (e2) {};
        \node [scale=0.3,circle,fill,label=below:{$e_3$}] at (e3) {};
        
        \small \node[scale=0.4,circle,fill,label=above:{$P$}] (P) at (barycentric cs:e1=0.5,e2=0.8,e3=0.4) {};
        
        \draw [dashed] (e1.south) -- (P.north east) node[pos=0.5,scale=0.1,label=below right:$\nu_1$] {};
        \draw [dashed] (e2.north east) -- (P.north west) node[pos=0.5,scale=0.1,label=below:$\nu_2$] {};
        \draw [dashed] (e3.north west) -- (P.south) node[pos=0.5,scale=0.1,label=left:$\nu_3$] {};
      \end{scope}
    \end{tikzpicture}
  }
  \caption{If $\mathcal{E}$ is finite, the de~Finetti mixture representation \cref{eq:de:finetti} and the more general representation
  \cref{eq:varadarajan} reduce to a finite convex combination. The points inside the set---\ie the distributions $P$ with the symmetry
  property defined by the group $\mathbb{G}$---can be represented as convex combinations $P=\sum_{e_i\in\mathcal{E}}\nu_i e_i$, with coefficients
  $\nu_i\geq 0$ satisfying $\sum_i\nu_i=1$. When $\mathcal{E}$ is infinite, an integral is substituted for the sum.
  }
  \label{fig:polytope}
\end{figure}
We have already encountered the components of \eqref{eq:varadarajan}
in \cref{sec:structures}:
In Bayesian terms, $\kernelval \genkernel \theta$ again corresponds to the observation distribution 
and $\nu$ to the prior.
Geometrically, the integral representation \cref{eq:varadarajan} can be
regarded as convex combination. \cref{fig:polytope} illustrates this
idea for a toy example with three ergodic measures.
A special case of the ergodic decomposition theorem 
is well-known in Bayesian theory as a result of Freedman
\citep{Freedman:1962:1,Freedman:1963:1}:
\begin{example}[Freedman's theorem]
  Consider a sequence ${X_1,X_2,\dotsc}$ as in de~Finetti's theorem. Now replace invariance under permutations
  by a stronger condition: Let $\mathrm O(n)$ be the orthogonal group of rotations and reflections
  on $\Reals^n$, \ie the set of
  ${n\times n}$ orthogonal matrices. We now demand that, if we regard any initial sequence of $n$ variables
  as a random vector in $\Reals^n$, then rotating and/or reflecting this vector does not change the distribution of the sequence:
  That is, for every ${n\in\Nats}$ and ${M\in \mathrm O(n)}$, 
  \begin{equation}
    \label{eq:rotation:invariance}
    (X_1,X_2,\dotsc) \eqdist (M(X_1,\dotsc,X_n), X_{n+1}, \dotsc)\,.
  \end{equation}
  In the language of \cref{theorem:varadarajan}, the group $G$ is the set of all rotations and reflections acting on all finite  prefixes of a sequence.
%  ${G=\cup_{n\in\Nats} \mathrm O(n)}$.
  For every $\sigma > 0$, let $\kernelval {\mathcal N} \sigma $ be the distribution of zero-mean Gaussian random variable with standard deviation $\sigma$.
  Freedman showed that, if $X^{\infty}$ satisfies \cref{eq:rotation:invariance}, 
  then its distribution is a scale mixture of Gaussians:
  \begin{equation}
    \label{eq:freedman}
    \Pr(X^{\infty}\in\argdot) 
      = 
      \int_{\NNReals} \kernelvalset {\mathcal N^\infty} \argdot \sigma \,\nu_{\NNReals}(\dee\sigma)\,.
  \end{equation}
  Thus, $\mathcal{E}$ contains all factorial distributions of zero-mean Gaussian distributions on $\Reals$,
  $\tspace$ is the set $\PosReals$ of variances, and $\nu$ a distribution on $\PosReals$.
\end{example}

Compared to de Finetti's theorem, where $\mathbb{G}$ is the group of
permutations, Freedman's theorem increases the size of $\mathbb{G}$:
Any permutation can be represented as an orthogonal matrix, but here rotations have been added
as well. In other words, we are strengthening the hypothesis by
imposing more constraints on the distribution of $X^{\infty}$. As a result,
the set $\mathcal{E}$ of ergodic measures
shrinks from all factorial measures to the set of factorials of zero-mean Gaussians. This is again 
an example of a general theme:
\begin{equation*}
  \text{larger group}
  \quad\longleftrightarrow\quad
  \text{more specific representation}
\end{equation*}
In contrast, the Aldous-Hoover theorem \emph{weakens} the hypothesis of de Finetti's theorem---in
the matrix case, for instance, the set of all permutations of the index set $\Nats^2$ is restricted
to those which preserve rows and columns---and hence yields a more general representation.

\begin{remark}[Symmetry and sufficiency]
  An alternative way to define symmetry in statistical models is through sufficient statistics:
  Intuitively, a symmetry property identifies information which is not relevant to the statistical
  problem; so does a sufficient statistic. For example, the empirical distribution retains
  all information about a sample except for the order in which observations were recorded.
  The empirical distribution is hence a sufficient statistic for the set distributions of exchangeable
  sequences.
  In an exchangeable graph model, the empirical graphon (the checkerboard function in
  \cref{fig:convergence}) is a sufficient statistic.
  If the sufficient statistic is finite-dimensional 
  and computes an average $\frac{1}{n}\sum_{i}S_0(x_i)$ over observations for 
  some function $S_0$, the ergodic distributions
  are exponential family models \cite{Kuechler:Lauritzen:1989}.
  A readable introduction to this
  topic is given by \citet{Diaconis:1992}.
  The definitive reference is the monograph of \citet{Lauritzen:1988}, who
  refers to the set $\mathcal{E}$ of ergodic distributions as an \emph{extremal family}.
\end{remark}

Not every probabilistic symmetry is applicable in statistics in the
same way as exchangeability is.
To be useful to statistics, the symmetry must satisfy two conditions:
\begin{enumerate}
\item The set $\mathcal{E}$ of ergodic measures should be a ``small'' subset of the set
  of symmetric measures.
\item The measures $\kernelval \genkernel \theta$ should have a tractable representation, such
  as Kingman's paint-box or the Aldous-Hoover sampling scheme.
\end{enumerate}
\cref{theorem:varadarajan} guarantees neither.
If (1) is not satisfied, the representation is useless for statistical purposes:
The integral representation \cref{eq:varadarajan} means that the information in
$X_{\infty}$ is split into two parts, the information contained in the parameter
value $\theta$ (which a statistical procedure tries to extract) and the randomness
represented by $\kernelval \genkernel \theta$ (which the statistical procedure discards).
If the set $\mathcal{E}$ is too large, $\Theta$ contains almost all the information
in $X_{\infty}$, and the decomposition becomes meaningless. We will encounter an
appealing notion of symmetry for sparse networks in the next section---which, however,
seems to satisfy neither condition (1) or (2).
It is not clear at present whether there are useful types of
symmetries which
do not imply some form of invariance to a group of permutations.
Although the question is abstract, the incompatibility of
sparseness and exchangeability means it is directly relevant to Bayesian
statistics.

\subsection{Stationary networks and involution invariance}

Is there a form of invariance that yields statistical models for
network data? There is indeed a very natural notion of invariance in
networks, called involution invariance, which we describe in more
detail below. This property has interesting mathematical properties
and admits an ergodic decomposition as in \cref{theorem:varadarajan},
but it seems to be too weak for applications in statistics.

A crucial difference between network structures and exchangeable graphs is that, in 
most networks, location in the graph matters. If conditioning on location is informative,
exchangeability is broken. Probabilistically, location is modeled by marking a distinguished
vertex in the graph.
A \defn{rooted graph} $(g,v)$ is simply a graph $g$ in which a particular vertex $v$ has been 
marked.
A very natural notion of invariance for networks is called
\emph{involution invariance} \citep{Aldous:Lyons:2007:1} 
or \emph{unimodularity}
\citep{Benjamini:Schramm:2001:1}, and
can be thought of as a form of stationarity:
\begin{definition}
  Let $P$ be the distribution of a random rooted graph, and 
  define a distribution $\tilde{P}$ as follows: A sample ${(G,w)\sim\tilde{P}}$
  is generated by sampling ${(G,v)\sim P}$, and then sampling
  $w$ uniformly from the neighbors of $v$ in $G$. 
  The distribution $P$ is called \defn{involution invariant} if
  $P=\tilde{P}$.
\end{definition}
The definition says that, if an observer randomly walks along the graph
$G$ by moving to a uniformly selected neighbor in each step,
the \emph{distribution} of the network around the observer remains
unchanged (although the actual neighborhoods in a sampled graph may vary).

Involution invariance is a symmetry property which admits an 
ergodic decomposition, and \citet*{Aldous:Lyons:2007:1} have characterized the
ergodic measures. 
This characterization is abstract, however, and there
is no known ``nice'' representation resembling, for example, the sampling
scheme for exchangeable graphs. Thus, of the two desiderata described
in \cref{sec:symmetry}, property (2) does not seem to hold.
We believe that property (1) does not hold either:
Although we have no proof at present, it seems that involution
invariance is too weak a constraint to yield interesting statistical
models (\ie the set of ergodic distributions is a ``large''
subset of the involution invariant distributions). 

Since exchangeability and involution invariance are the only well-studied
probabilistic symmetries for random graphs, the question how
statistical models of networks can be characterized is an open problem:
\begin{center}
  \emph{ 
    Is there a notion of probabilistic symmetry whose ergodic measures
    in \eqref{eq:varadarajan} describe useful statistical models for
    sparse graphs with network properties?
  }
\end{center}
There is a sizeable literature on sparse random graph models which
can model power laws and other network properties; see, for example
\citep{Durrett:2006}. These are probability models and can be
simulated, but estimation from data is often intractable, due
to stochastic dependencies between the edges in the random graph.
On the other hand, \emph{some} dependence between edges is necessary to obtain a power
law and similar properties. Hence, a suitable notion of symmetry would
have to restrict dependencies between edges sufficiently to permit
statistical inference, but not to the full conditional independence
characteristic of the exchangeable case.

\section{Further References}
\label{sec:references}

Excellent non-technical references on the general theory of exchangeable arrays and other exchangeable random structures
are two recent surveys by \citet{Aldous:2009,Aldous:2010:1}. His well-known lecture notes \citep{Aldous:1983}
also cover exchangeable arrays.
The most comprehensive available reference on the general theory is the monograph by \citet{Kallenberg:2005} (which
presupposes in-depth knowledge of measure-theoretic probability).
Kingman's original article \citep{Kingman:1978:2} provides a concise reference on exchangeable random partitions.
A thorough, more technical treatment of exchangeable partitions can be found in \citep{Bertoin:2006}.

\citet{Schervish:1995} gives an insightful discussion of the application of exchangeability to Bayesian statistics.
There is a close connection between probabilistic symmetries (such as exchangeability) and sufficient statistics, which is covered
by a substantial literature. See \citet{Diaconis:1992} for an introduction and further references. 
For applications of
exchangeability results to machine learning models, see \citep{Fortini:Petrone:2012}, who discuss applications 
of the partial exchangeability result of \citet{Diaconis:Freedman:1980:1} to the infinite hidden Markov model
\citep{Beal:Ghahramani:Rasmussen:2002}.

The theory of graph limits in its current form was initiated by Lov\'asz and Szegedy \citep{Lovasz:Szegedy:2006,Lovasz:Szegedy:2007} and
\citet{Borgs:Chayes:Lovasz:Sos:Szegedy:Vesztergombi:2005}. It builds on work of \citet{Frieze:Kannan:1999}, who introduced both the
weak regularity lemma (\cref{theorem:regularity:weak}) and the cut norm $\dcut$.
In the framework of this theory, the Aldous-Hoover representation of exchangeable graphs can be derived by purely analytic
means \citep[][Theorem 2.7]{Lovasz:Szegedy:2006}. The connection between graph limits and Aldous-Hoover theory was established,
independently of each other, by \citet{Diaconis:Janson:2007} and by \citet{Austin:2008}. 
An accessible introduction to the analytic perspective is the survey \citep{Lovasz:2009:1}, which assumes basic familiarity with measure-theoretic
probability and functional analysis, but is largely non-technical. The
monograph \citep{Lovasz:2013:A} gives a comprehensive account.

Historically, the Aldous-Hoover representation was established in
independent works of Aldous and of Hoover in the late 1970s.
Aldous' proof used probability-theoretic methods, whereas Hoover, a logician, leveraged techniques from model theory. 
In 1979, Kingman \citep{Kingman:1979:1} writes
\begin{quote}
...a general solution has now been supplied by Dr David Aldous of Cambridge. [...]
The proof is at present very complicated, but there is reason to hope that the techniques developed can be applied to more general experimental designs.
\end{quote}
Aldous' paper \citep{Aldous:1981}, published in 1981, attributes the idea of the published version of the proof to Kingman.
The results were later generalized considerably by Kallenberg \citep{Kallenberg:1999}.

\mbox{ }\\
{\noindent\bf Acknowledgments.} We have learned much about random graphs from Cameron Freer, and are also indebted to James Robert Lloyd for many useful discussions. Two anonymous reviewers have provided very detailed feedback, which in our opinion has greatly improved the article. We thank Karolina Dziugaite, Creighton Heaukalani, Jonathan Huggins and Christian Steinr\"ucken for helpful comments on the manuscript.

\bibliographystyle{natbib}
\bibliography{bibtex-preamble,references}

\end{document}